\documentclass[12pt]{amsart}
\usepackage[]{amscd,amsfonts,amsmath,amsxtra,amssymb}
\usepackage{graphics}
\usepackage[all]{xy}
\usepackage{hyperref}
\usepackage{eucal}
\usepackage{rotating}
\usepackage[T1]{fontenc}
\newtheorem{sub}{}[section]
\newtheorem{subsub}{}[sub]
\usepackage{rotating}
\usepackage[active]{srcltx}
\def\ov#1{\overline{#1}}

\def\coker{\mathop{\rm coker}\nolimits}
\def\Hom{\mathop{\rm Hom}\nolimits}
\def\Surj{\mathop{\rm Surj}\nolimits}
\def\Inj{\mathop{\rm Inj}\nolimits}
\def\HHom{\mathop{\mathcal Hom}\nolimits}

\def\Ext{\mathop{\rm Ext}\nolimits}
\def\EExt{\mathop{\mathcal Ext}\nolimits}
\def\Tor{\mathop{\rm Tor}\nolimits}

\def\Pic{\mathop{\rm Pic}\nolimits}

\def\Aut{\mathop{\rm Aut}\nolimits}

\def\End{\mathop{\rm End}\nolimits}

\def\EEnd{\mathop{\mathcal End}\nolimits}

\def\Gr{\mathop{\rm Gr}\nolimits}

\def\imm{\mathop{\rm im}\nolimits}
\def\deg{\mathop{\rm deg}\nolimits}

\def\rk{\mathop{\rm rk}\nolimits}
\def\Deg{\mathop{\rm Deg}\nolimits}

\def\spec{\mathop{\rm spec}\nolimits}

\def\lra{\longrightarrow}

\def\sigg{\mathop{\hbox{$\displaystyle\sum$}}\limits}

\def\hfl#1#2{\smash{\mathop{\ \hbox to 12mm{\rightarrowfill}}
\limits^{\scriptstyle#1}_{\scriptstyle#2} \ }}
\def\hflb#1#2{\smash{\mathop{\hbox to 12mm{\leftarrowfill}}
\limits^{\scriptstyle#1}_{\scriptstyle#2}}}

\def\m#1{{\hbox{$#1$}}}

\def\ot{\otimes}
\def\og{\leavevmode\raise.3ex\hbox{$\scriptscriptstyle\langle\!\langle$}}
\def\fg{\leavevmode\raise.3ex\hbox{$\scriptscriptstyle\,\rangle\!\rangle$}}
\def\span#1{\langle#1\rangle}

\def\nsp{\lbrace 0\rbrace}

\def\dsp{\displaystyle}
\def\Ssect#1#2{\pagebreak[3]\begin{sub}\label{#2}{\sc\small\small
#1}\rm\medskip}

\def\sepsec{\vskip 1.8cm}
\def\sepsub{\vskip 1cm}
\def\sepsubsub{\vskip 0.75cm}
\def\sepprop{\vskip 0.5cm}

\def\xmat#1{\[\xymatrix{#1}\]}
\def\flinc{\ar@{^{(}->}}
\def\fleq{\ar@{=}}
\def\flon{\ar@{->>}}
\def\fmaps{\ar@{|-{>}}}
\def\fflat{\ar@{-}}
\def\fimpl{\ar@{=>}}
\def\Nligne{\hfil\break}

\def\ED{\vskip 1cm\end{document}}

\def\REF#1{\m{(\ref{#1})}}

\newcommand{\B}{{\mathbb B}}

\newcommand{\C}{{\mathbb C}}

\newcommand{\D}{{\mathbb D}}
\newcommand{\F}{{\mathbb F}}
\newcommand{\E}{{\mathbb E}}
\newcommand{\G}{{\mathbb G}}

\newcommand{\Y}{{\mathbb Y}}
\newcommand{\U}{{\mathbb U}}
\newcommand{\V}{{\mathbb V}}

\renewcommand{\L}{{\mathbb L}}

\def\T{{\mathbb T}}

\newcommand{\kc}{{\mathcal C}}
\newcommand{\kd}{{\mathcal D}}
\newcommand{\ke}{{\mathcal E}}
\newcommand{\kf}{{\mathcal F}}
\newcommand{\kg}{{\mathcal G}}

\newcommand{\ki}{{\mathcal I}}

\newcommand{\kl}{{\mathcal L}}
\newcommand{\km}{{\mathcal M}}
\newcommand{\kn}{{\mathcal N}}
\newcommand{\ko}{{\mathcal O}}

\newcommand{\ku}{{\mathcal U}}
\newcommand{\kv}{{\mathcal V}}

\newcommand{\kx}{{\mathcal X}}
\newcommand{\ky}{{\mathcal Y}}
\newcommand{\kz}{{\mathcal Z}}
\begin{document}

\def\refname{References}
\def\contentsname{Summary}
\def\proofname{Proof}
\def\abstractname{Resume}

\author{Jean--Marc Dr\'{e}zet}
\thanks{This paper was partially written during
a stay in Tata Institute of Fundamental Research (Mumbai)}
\address{
Institut de Math\'ematiques de Jussieu - Paris Rive Gauche\\
Case 247\\
4 place Jussieu\\
F-75252 Paris, France}
\email{jean-marc.drezet@imj-prg.fr}
\title[{Deformations of moduli spaces}] {Reachable sheaves on ribbons
and deformations of moduli spaces of sheaves}

\begin{abstract} A {\em primitive multiple curve} is a Cohen-Macaulay
irreducible projective curve $Y$ that can be locally embedded in a smooth
surface, and such that $C=Y_{red}$ is smooth. In this case, $L=\ki_C/\ki_C^2$ is
a line bundle on $C$. If $Y$ is of multiplicity 2, i.e. if \ \m{\ki_C^2=0}, $Y$
is called a {\em ribbon}. If $Y$ is a ribbon and \ \m{h^0(L^{-2})\not=0}, then
$Y$ can be deformed to smooth curves, but in general a coherent sheaf on $Y$
cannot be deformed in coherent sheaves on the smooth curves. 

It has been proved in \cite{dr8} that a ribbon with associated
line bundle $L$ such that \m{\deg(L)=-d<0} can be deformed to reduced curves
having 2 irreducible components if $L$ can be written as
\[L \ = \ \ko_C(-P_1-\cdots-P_d) \ , \]
where \m{P_1,\cdots,P_d} are distinct points of $C$. In this case we prove
that quasi locally free sheaves on $Y$ can be deformed to torsion free sheaves
on the reducible curves with two components. This has some consequences on
the structure and deformations of the moduli spaces of semi-stable sheaves on
$Y$.
\end{abstract}

\maketitle
\tableofcontents


\section{Introduction}

A {\em primitive multiple curve} is an algebraic variety $Y$ over $\C$ which is
Cohen-Macaulay, such that the induced reduced variety \m{C=Y_{red}} is a
smooth projective irreducible curve, and that every closed point of $Y$ has a
neighbourhood that can be embedded in a smooth surface. These curves have been
defined and studied by C.~B\u anic\u a and O.~Forster in \cite{ba_fo}. The
simplest examples are infinitesimal neighbourhoods of projective smooth curves
embedded in a smooth surface (but most primitive multiple curves cannot be
globally embedded in smooth surfaces).

Let $Y$ be a primitive multiple curve with associated reduced curve $C$, and
suppose that \m{Y\not=C}. Let \m{\ki_C} be the ideal sheaf of $C$ in $Y$. The
{\em multiplicity} of $Y$ is the smallest integer $n$ such that \m{\ki_C^n=0}.
The sheaf \ \m{L=\ki_C/\ki_C^2} \ is a line bundle on $C$, called the {\em line
bundle on $C$ associated to $Y$}.

Primitive multiple curves of multiplicity 2 are called {\em ribbons}. They have
been parametrised in \cite{ba_ei}. Primitive multiple curves of any multiplicity
and the coherent sheaves on them have been studied in \cite{dr1}, \cite{dr2},
\cite{dr4} and \cite{dr5}.

The deformations of ribbons to smooth projective curves have been studied by
M.~Gonz\'alez in \cite{gon}: he proved a ribbon $Y$, with associated smooth
curve $C$ and associated line bundle $L$ on $C$ is smoothable if \
\m{h^0(L^{-2})\not=0}.

\sepsub

\Ssect{Deformations to reduced reducible curves}{def1}

Deformations of primitive multiple curves \m{Y=C_n} of any multiplicity
\m{n\geq 2} to reduced curves having multiple components which are smooth,
intersecting transversally, have been studied in \cite{dr7}: we consider flat
morphisms \m{\pi:\kc\to S}, where $S$ is a smooth connected curve, for some
\m{P\in S}, \m{\kc_P} is isomorphic to $Y$ and for \m{s\not=P}, \m{\kc_s} is a
curve with multiple components which are smooth and intersect transversally.
Then $n$ is the maximal number of components of such deformations of $Y$. 
In this case we say that the deformation is {\em maximal}, the number of
intersection points of two components is exactly \m{-\deg(L)} and the genus of
the components is the one of $C$ (these deformations are called {\em maximal
reducible deformations}).
If \m{\deg(L)=0} we call $\pi$ a {\em fragmented deformation}. This case
\m{\deg(L)=0} has been completely treated in \cite{dr7}: a primitive
multiple curve of multiplicity $n$ can be deformed in disjoint unions of $n$
smooth curves if and only if \m{\ki_C} is isomorphic to the trivial bundle on
\m{C_{n-1}}. In \cite{dr6} it has been proved that this last condition is
equivalent to the following: there exists a flat family of smooth curves
\m{\kc\to S}, parametrised by a smooth curve $S$, \m{s_0\in S} such that
\m{\kc_{s_0}=C}, such that $Y$ is isomorphic to the $n$-th infinitesimal
neighbourhood of $C$ in $\kc$.

The problem of determining which primitive multiple curves of multiplicity $n$
can be deformed to reduced curves having exactly $n$ components, allowing
intersections of the components, is more difficult. A necessary condition is \ 
\m{h^0(L^*)>0}. It has been proved in \cite{dr8} that a ribbon with associated
line bundle $L$ such that \m{\deg(L)=-d<0} can be deformed to reduced curves
having 2 irreducible components if $L$ can be written as
\[L \ = \ \ko_C(-P_1-\cdots-P_d) \ , \]
where \m{P_1,\cdots,P_d} are distinct points of $C$.

\end{sub}

\sepsub

\Ssect{Deformations of sheaves}{def2}

Let $Y$ be a ribbon, \m{C=Y_{\text{red}}}, and  $L$ the associated line bundle
on $C$. We suppose that $L$ can be written as
\[L \ = \ \ko_C(-P_1-\cdots-P_d) \ , \]
where \m{P_1,\cdots,P_d} are distinct points of $C$, or that \m{L=\ko_C}.
Let $\ke$ be a coherent sheaf on $Y$. Let \ \m{p:\ke\to\ke_{|C}} \ be the
restriction morphism and \ \m{\ke_1=\ker(p)}, which is concentrated on $C$. We
call \ \m{R(\ke)=\rk(\ke_{|C})+\rk(\ke_1)} \ the {\em generalised rank} of
$\ke$, and \ \m{\Deg(\ke)=\deg(\ke_{|C})+\deg(\ke_1)} \ the {\em generalised
degree} of $\ke$. It has been proved in \cite{dr2} that there exist unique
integers \m{a,b\geq 0} and a nonempty open subset \m{U\subset C} such that for
every \m{x\in U} we have \ \m{\ke_x\simeq a\ko_{Y,x}\oplus b\ko_{C,x}}. If this
is true for every \m{x\in C}, $\ke$ is called a {\em quasi locally free sheaf},
and this is the case if and only if \m{\ke_{|C}} and \m{\ke_1} are vector
bundles on $C$.

Let $S$ be a smooth connected curve, \m{P\in S} and \m{\rho:\kx\to S} a
deformation of $Y$ to smooth curves, i.e. $\rho$ is flat, \m{\rho^{-1}(P)\simeq
Y}, and for every \m{s\in S\backslash\{P\}}, \m{\kx_s=\rho^{-1}(s)} \ is a
smooth projective irreducible curve. It is easy to see (using Hilbert
polynomials with respect to an ample line bundle on $\kx$) that if $\kf$ is a
coherent sheaf on $\kx$ flat on $S$, then \m{R(\kf_{|Y})} must be even. This
means that not all sheaves on $Y$ can be deformed to sheaves on the smooth
fibres. Moreover if \m{\kf_{|Y}} is locally free, \m{\deg(\kf_s)} is even.

Now let \m{\pi:\kc\to S} be a maximal reducible deformation of $Y$: we have
\m{\pi^{-1}(P)\simeq Y}, and for every \m{s\in S\backslash\{P\}},
\m{\kc_s=\pi^{-1}(s)} \ is a reduced projective curve having 2 irreducible
components, smooth and intersecting transversally in \m{-\deg(L)} points (see
\ref{prel2} for the precise definition). Let $E$ be a coherent sheaf on $Y$. We
say that {\em $E$ is reachable} (with respect to $\pi$) if there exists a
smooth connected curve \m{S'}, \m{P'\in S'}, a non constant morphism
\m{f:S'\to S} such that \m{f(P')=P} and a coherent sheaf $\E$ on
$f^*(\kc)=\kc\times_SS'$, flat on \m{S'}, such that \ \m{\E_{P'}\simeq E}.

The main result of this paper is theorem \ref{theo6}

\sepprop

\begin{subsub}{\bf Theorem: } Every quasi locally free sheaf on 
$Y$ is reachable.
\end{subsub}

\sepprop

To prove this result we need to study the torsion free sheaves on reducible
curves with several components, and the moduli spaces of such semi-stable
sheaves. Some work has been done on this subject by M.~Teixidor~i~Bigas in
\cite{tei1}, \cite{tei2}, \cite{tei3}, \cite{tei4}, but only for vector
bundles. We will need here to consider sheaves which have not the same rank on
all the components. We will also use blowing-up of ribbons, which are also used
in \cite{ch-ka} in another context.

\end{sub}

\sepsub

\Ssect{Motivation}{moti}

\begin{subsub} Study of the non reduced structure of some moduli spaces of 
stable sheaves on $Y$ -- \rm
Some coherent sheaves on projective varieties have a non reduced versal
deformation space. In particular, some moduli spaces of stable sheaves are non 
reduced. 

We treat here the case of some sheaves on ribbons: the {\em quasi locally free 
sheaves of rigid type} (see. \ref{mod_rig_r} and \ref{rel_m} for more details). 
Let $Y$ be a ribbon and \m{\pi:\kc\to S} a deformation of $Y$ to reduced curves 
with two components as in \ref{def1}. 

Let $E$ be a stable quasi locally free sheaf of rigid type on $Y$. From 
\cite{dr4}, $E$ belongs to an unique irreducible component of the corresponding 
moduli space of stable sheaves. Let $\bf M$ be the open subset of this 
component corresponding to quasi locally free sheaves of rigid type. It is 
proved in \cite{dr4} that ${\bf M}_{reg}$ is smooth, and that the tangent sheaf 
\m{T{\bf M}} of $\bf M$ is locally free. Hence we obtain a vector bundle  
\m{T{\bf M}/T({\bf M}_{reg})} on ${\bf M}_{reg}$, its rank is \m{h^0(L^*)}.

We find, using theorem 1.2.1, that $E$ can be deformed in two distinct ways to 
sheaves on the reduced curves. In particular $\bf M$ deforms to two components 
of the moduli spaces of sheaves on the reduced curves, and $\bf M$ appears as 
the ``limit'' of varieties with two components,  whence the non reduced 
structure of $\bf M$.

In \cite{dr9} we show that this limit is responsible of a sub-line bundle $\kl$ 
of \m{T{\bf M}/T({\bf M}_{reg})}, and that $\kl$ is closely related to the 
deformation $\pi$.
\end{subsub}

\medskip

\begin{subsub} Moduli spaces of vector bundles -- \rm
If instead of a deformation of $Y$ to reducible curves on considers a 
deformation \m{\rho:\kd\to S} to {\em smooth} curves, among the quasi locally 
free sheaves on $Y$ only some of them can be deformed to sheaves on the smooth 
fibers $\kd_s$. For example if $Y$ can be deformed to reduced curves
having 2 irreducible components, then it can be deformed to smooth curves 
(\cite{dr8}, prop. 2.5.1), and it follows easily from theorem 1.2.1 that the 
vector bundles on $Y$ can be deformed to vector bundles on the smooth fibers.
The moduli spaces of stable vector bundles on $Y$ (which are smooth) deform 
then to moduli spaces of vector bundles on the smooth fibers $\kd_s$. It would 
then be possible to deduce properties of moduli spaces of vector bundles on the 
smooth fibers from the study of the moduli spaces of stable bundles on the 
special fiber $Y$. The study of the moduli spaces of stable bundles $Y$ can 
also use that of moduli spaces of more complicated sheaves as those of quasi 
locally free sheaves of rigid type (if, as in \cite{dr4} some proofs by 
induction are needed).
\end{subsub}

\end{sub}

\sepsub

\Ssect{Moduli spaces of (semi-)stable sheaves on primitive reducible
curves}{rel_mod0}

\begin{subsub}\label{primred2} Primitive reducible curves -- \rm (cf.
\ref{prim_red}-). A {\em primitive reducible curve} is an algebraic curve $X$
such that
\begin{enumerate}
\item[--] $X$ is connected.
\item[--] The irreducible components of $X$ are smooth projective curves.
\item[--] Any three irreducible components of $X$  have no common point, and 
any two components are transverse.
\end{enumerate}
Suppose that $X$ has two components \m{D_1}, \m{D_2}. Let $\ke$ be a coherent
sheaf on $X$. Then $\ke$ is torsion free if and only if there exist vector
bundles \m{E_i} on \m{D_i}, \m{i=1,2}, and for every \m{x\in D_1\cap D_2} a
vector space \m{W_x} and surjective maps \ \m{f_{i,x}:E_{i,x}\to W_x}, such that
$\ke$ is isomorphic to the kernel of the surjective morphism of sheaves on $X$
\xmat{E_1\oplus E_2\ar[rrr]^-{\bigoplus(f_{1,x},-f_{2,x})} & & & \bigoplus_{x\in
D_1\cap D_2}\ov{W_x}}
(where \m{\ov{W_x}} is the skyscraper sheaf concentrated on $x$ with fibre
\m{W_x} at this point). For \m{i=1,2}, we have \ \m{E_i=\ke_{|D_i}/T_i}, where
\m{T_i} is the torsion subsheaf. We call
\[\tau(\ke) \ = \ (\rk(E_1),\deg(E_1),\rk(E_2),\deg(E_2))\]
the {\em type} of $\ke$. We have, for every \m{x\in D_1\cap D_2},
\m{\dim(W_x)\leq\inf(\rk(E_1),\rk(E_2))}. If for every \m{x\in D_1\cap D_2},
\m{\dim(W_x)=\inf(\rk(E_1),\rk(E_2))}, we say that {\em $\ke$ is linked}. For
example, if \m{\rk(E_1)=\rk(E_2)}, $\ke$ is linked if and only if $\ke$ is
locally free. We prove in \ref{prim_red}- that deformations of linked sheaves
are linked, and that \m{\deg(E_1)} and \m{\deg(E_2)} (and hence also
\m{\tau(\ke)}) are invariant by deformations. Moreover, any two linked sheaves
with the same type can be put in a flat family of linked sheaves on $X$
parametrised by an integral variety (theorem \ref{lem1m}). We prove also that
every torsion free sheaf on $X$ can be deformed to linked sheaves (theorem
\ref{theo5}). 
\end{subsub}

\sepprop

\begin{subsub}Moduli spaces of sheaves -- \rm Let \m{\ko_X(1)} be an ample line
bundle on $X$ and $H$ a polynomial in one variable with rational coefficients.
Let \m{\km_X(H)} be the moduli space of semi-stable torsion free sheaves on $X$
with Hilbert polynomial $H$ with respect to \m{\ko_X(1)}. It follows from
\ref{primred2} that every irreducible component $\kn$ of \m{\km_X(H)} has a
dense open subset corresponding to linked sheaves (having all the same type
$\tau$), and every semi-stable linked sheaf of type $\tau$ defines a point of
$\kn$. It follows that the components of \m{\km_X(H)} are indexed by the types
of the linked semi-stable sheaves of Hilbert polynomial $H$. We will
denote by \m{\kn_X(r_1,d_1,r_2,d_2)} the component containing points
corresponding to linked sheaves of type \m{(r_1,d_1,r_2,d_2)}.
We prove that \m{\kn_X(r_1,d_1,r_2,d_2)} is smooth at every point corresponding
to a stable linked sheaf $\ke$ such that \m{E_1} and \m{E_2} are simple vector
bundles (theorem \ref{theo4}). But in general, we have \
\m{\Ext^2_{\ko_X}(\ke,\ke)\not=\nsp}.
\end{subsub}

\end{sub}

\sepsub

\Ssect{Moduli spaces of sheaves of rigid type on ribbons}{mod_rig_r}

Let $Y$ be a ribbon, \m{C=Y_{\text{red}}} and $L$ the associated line bundle on
$C$. A quasi locally free sheaf $\ke$ on $Y$ is called {\em of rigid type} if it
is locally free, or locally isomorphic to \m{a\ko_Y\oplus\ko_C}, for some
integer \m{a\geq 0}. If $\ke$ is of rigid type, then the deformations of $\ke$
are also of rigid type, and \m{\deg(\ke_{|C})}, \m{\deg(\ke_1)} are invariant by
deformations. Moreover, any two quasi locally free sheaves of rigid type with
the same \m{\rk(\ke_{|C})}, \m{\deg(\ke_{|C})}, \m{\deg(\ke_1)} can be put in a
flat family of rigid sheaves with the same invariants parametrised by an
integral variety (cf. \cite{dr4}).

The semi-stability conditions on $Y$ are the same for every choice of an ample
line bundle on $Y$ (cf. \ref{coh_prim}-). The Hilbert polynomial of a coherent
sheaf $\ke$ on $Y$ depends only on the two invariants \m{R(\ke)} and
\m{\Deg(\ke)}. Let \m{\km_Y(R,D)} be the moduli space of semi-stable torsion
free sheaves on $Y$ with generalised rank $R$ and generalised degree $D$.
Suppose that \m{R=2a+1} and \m{d_0}, \m{d_1} are integers such that \
\m{D=d_0+d_1+a.\deg(L)}. Let \m{\kn(R,d_0,d_1)} denote the open irreducible
subset of \m{\km_Y(R,D)} corresponding to stable sheaves of rigid type $\ke$
such that \m{\deg(\ke_{|C})=d_0} and \m{\deg(\ke_1\ot L^*)=d_1} (it is denoted
by \m{\kn(a,1,d_0,d_1)} in \cite{dr4} and \cite{dr5}). We give in \cite{dr5}
sufficient conditions for the non-emptiness of \m{\kn(R,d_0,d_1)} (cf.
\ref{rig_t}). If it is non empty, then \m{\kn(R,d_0,d_1)_{\text{red}}} is
smooth, and at the point corresponding to the stable sheaf $\ke$, its tangent
space is canonically isomorphic to \m{H^1(\End(\ke))}, whereas the tangent
space of \m{\kn(R,d_0,d_1)} is isomorphic to \m{\Ext^1_{\ko_X}(\ke,\ke)}. The
quotient
\[\Ext^1_{\ko_X}(\ke,\ke)/H^1(\End(\ke)) \ \simeq \ H^0(\EExt^1(\ke,\ke))\]
is canonically isomorphic to \m{H^0(L^*)} (cf. \cite{dr2}, cor. 6.2.2).

\end{sub}

\sepsub

\Ssect{Relative moduli spaces of semi-stable sheaves}{rel_m}

We use the notations of \ref{def2}. Let $\ke$ be a quasi locally free
sheaf of rigid type on $Y$, locally isomorphic to \m{a\ko_Y\oplus\ko_C}. Let
\m{\pi:\kc\to S} be a maximal reducible deformation of $Y$. Then by the main
result, in \ref{def2}, \m{\ke} can be deformed to torsion free sheaves on the
primitive reducible curves \m{\kc_s}, \m{s\not=P}. 

Let \m{\ko_\kc(1)} be an ample line bundle on $\kc$. Then we have, for every
\m{s\in S\backslash\{P\}},\Nligne
\m{\deg(\ko_{\kc_{1,s}}(1))=\deg(\ko_{\kc_{2,s}}(1))}, and so the semi-stability
of sheaves on the curves \m{\kc_s} does not depend on the choice of
\m{\ko_\kc(1)}. Let $H$ be the Hilbert polynomial of $\ke$. Let
\[\tau:{\bf M}(\kc/S,H)\lra S\]
be the relative moduli space of sheaves on the fibres of $\pi$, so that for
every closed point \m{s\in S}, \m{\tau^{-1}(s)} is the moduli space
\m{\km_{\kc_s}(H)} of semi-stable sheaves on \m{\kc_s} with Hilbert polynomial
$H$.

We have \m{\rk(\ke_{|C})=a+1} and \m{\rk(\ke_1)=a}. Let
\m{d_0=\deg(\ke_{|C})}, \m{d_1=\deg(\ke_1)}. Let \m{S'} be a smooth connected
curve, \m{P'\in S'}, \m{f:S'\to S} a non constant morphism such that \m{f(P')=P}
and $\E$ a coherent sheaf on $f^*(\kc)=\kc\times_SS'$, flat on \m{S'}, such that
\ \m{\E_{P'}\simeq\ke}. We prove in proposition \ref{prop14} that for
\m{t\not=P'} in a suitable neighbourhood of \m{P'}, \m{\E_t} is a linked sheaf
such that, if \ \m{E_i=\E_{t|\kc_{i,f(t)}}/T_i} (where \m{T_i} is the
torsion subsheaf), then either
\[\rk(E_1) \ = a \ , \quad \deg(E_1) \ = \ d_1 \ , \quad \rk(E_2) \ = a+1 \ ,
\quad \deg(E_2) \ = \ d_0 \ , \]
or
\[\rk(E_2) \ = a \ , \quad \deg(E_2) \ = \ d_1 \ , \quad \rk(E_1) \ = a+1 \ ,
\quad \deg(E_1) \ = \ d_0 \ . \]
In other words, in the moduli space \m{{\bf M}(\kc/S,H)},
\m{\kn(2a+1,d_0,d_1)\subset\km_Y(H)} is the ''limit'' of the only components
\m{\kn_{\kc_s}(a+1,d_0,a,d_1)}, \m{\kn_{\kc_s}(a,d_1,a+1,d_0)} of
\m{\km_{\kc_s}(H)}, \m{s\not=P}. In particular this explains (cf. \ref{moti}) 
why \m{\kn(2a+1,d_0,d_1)} cannot be reduced.

\end{sub}

\sepsub

\Ssect{The proof of the main theorem}{pro}

We use the notations of \ref{def2}. Let \m{\kz\subset\kc} be the closure in
$\kc$ of the locus of the intersection points of the components of
\m{\pi^{-1}(s)}, \m{s\not=P}. It consists of \m{d=-\deg(L)} smooth curves
intersecting $C$ in distinct points \m{P_1,\cdots,P_d} such that \
\m{L=\ko_C(-P_1-\cdots-P_d)}. The blowing-up \m{\kd\to\kc} of $\kz$ is a
fragmented deformation of the blowing up of \m{P_1,\cdots,P_d} in $Y$, which is
a ribbon with associated smooth curve $C$ and associated line bundle \m{\ko_C}.
We first prove the theorem for $\kd$. By studying the relations between sheaves
on $\kc$ and sheaves on $\kd$ we can prove the theorem for $\kc$.

\sepsub

\end{sub}

\Ssect{Outline of the paper}{outl}

{\bf Section 2} contains reminders of definitions and results, and some
technical lemmas used in the rest of the paper.

{\bf Section 3} contains mainly reminders of the main properties of coherent
sheaves on ribbons, with some useful lemmas.

{\bf Section 4} is devoted to the study of torsion free sheaves on primitive
reducible curves (cf. \ref{rel_mod0}). The results obtained here could probably
be useful if one wants to study more precisely the moduli spaces of sheaves on
these curves, without restricting to the locally free ones.

In {\bf Section 5} we consider the blowing-up of a finite number of points
of a ribbon $Y$. It is again a ribbon $\widetilde{Y}$. We study the relations
between torsion free sheaves on $Y$ and torsion free sheaves on $\widetilde{Y}$.

In {\bf Section 6} we study the coherent sheaves $\ke$ on $\kc$ such that
\m{\ke_{|Y}} is a quasi locally free sheaf. 

In {\bf Section 7} we prove the main theorem, using the results of the
preceding sections.

\end{sub}

\sepsub

{\bf Notations -- } If $X$ is a smooth variety and $D$ a divisor of $X$, 
\m{\ko_X(D)} is the line bundle associated to $D$, i.e. the dual of the ideal
sheaf \m{\ko_X(-D)} of $D$. If $\ke$ is a coherent sheaf on $X$, let \
\m{\ke(D)=\ke\ot\ko_X(D)}.

-- In this paper, an {\em algebraic variety} is a quasi-projective scheme over 
$\C$.

-- If $X$ is an algebraic variety and \m{Y\subset X} a closed subvariety, let 
\m{\ki_{Y,X}} denote the ideal sheaf of $Y$ in $X$. If there is no ambiguity,
we will denote it also by \m{\ki_Y}.

-- if \m{x\in X} is a closed point, \m{\C_x} will denote the skyscraper sheaf 
concentrated on $x$ with fibre $\C$ at $x$.

-- if \m{f:\kx\to S} is a flat morphism of algebraic varieties, and \m{f:T\to S}
another morphism, we shall sometimes use the notation \
\m{f^*(\kx)=\kx\times_ST}.

\sepsec

\section{Preliminaries}

\Ssect{Relative moduli spaces of (semi-)stable sheaves}{rel_mod}

(cf. \cite{ma2}, \cite{si})

Let $X$, $S$ be algebraic varieties and \ \m{X\to S} \ a flat projective 
morphism. Let \m{\sigma:S'\to S} be an $S$-variety. Let \ \m{X'=X\times_SS'}. A
{\em family of sheaves} on \m{X'/S'} is a coherent sheaf on \m{X\times_SS'},
flat on \m{S'}.

Let \m{\ko_X(1)} be an ample line bundle on $X$ and $H$ a polynomial in one 
variable with rational coefficients. A {\em family of semi-stable sheaves} with 
Hilbert polynomial $H$ on \m{X'/S'} is a family of sheaves $\ke$ on \m{X'/S'} 
such that for every closed point \m{s'\in S'}, the restriction \m{\ke_{s'}} of 
$\ke$ to \ \m{X_{\sigma(s')}=X'_{s'}} is a semi-stable sheaf with Hilbert 
polynomial $H$ with respect to 
\m{\ko_{X_{\sigma(s')}}(1)=\ko_X(1)_{|X_{\sigma(s')}}}. According to \cite{si},
there exists a moduli space
\[\tau:{\bf M}(X/S,H)\lra S\]
such that for every closed point \m{s\in S}, \m{\tau^{-1}(s)={\bf M}(X_s,H)} \ 
is the moduli space of semi-stable sheaves on \m{X_s} with Hilbert polynomial 
$H$. To every family $\ke$ of semi-stable sheaves with Hilbert polynomial $H$ on
\m{X'/S'}, one associates a canonical morphism
\[f_\ke:S'\lra{\bf M}(X/S,H)\]
such that for every \m{s'\in S'}, \m{f_\ke(s')} is the point of \m{{\bf 
M}(X_{\sigma(s')},H)} corresponding to \m{\ke_{s'}}. The morphism $\tau$ needs 
not to be flat (cf. \cite{in2}).
\end{sub}

\sepsub

\Ssect{Restrictions of sheaves on hypersurfaces}{res_hyp}

Let $X$ be a smooth algebraic variety, and $Y$, $Z$ hypersurfaces such that 
\m{Y\subset Z}. Let \ \m{Z=Y\cup T}, with $T$ a minimal hypersurface. Then we 
have \ \m{\ko_X(-Z)\subset\ko_X(-Y)}. Let $E$ be a vector bundle on $X$. Then 
we have a commutative diagram with exact rows
\xmat{0\ar[r] & E\ar[r]\fleq[d] & E(Y)\ar[r]\flinc[d] & 
E(Y)_{|Y}\ar[d]^{\phi_E}\\
 0\ar[r] & E\ar[r] & E(Z)\ar[r] & E(Z)_{|Z}}

\sepprop

\begin{subsub}\label{lem1d}{\bf Lemma: } Let \m{z\in Z}. Then
\[\imm(\phi_E) \ = \ \ki_{T,Z}\big[E(Z)_{|Z}\big] \ . \]
\end{subsub}

(recall that \m{\ki_{T,Z}} is the ideal sheaf of $T$ in $Z$)

\begin{proof} Let \m{f,g\in\ko_{X,z}} be equations of $Y$, $T$ respectively. 
Then \m{fg} is an equation of $Z$. Then \m{\ko_X(Y)_z} is the dual of \ 
\m{\ko_X(-Y)_z=(f)\subset\ko_{X,z}}, and is generated by
\[\xymatrix@R=5pt{(f)\ar[r] & \ko_{X,z}\\ \alpha.f\fmaps[r] & \alpha}\]
Similarly, \m{\ko_X(Z)_z} is the dual of \ 
\m{\ko_X(-Z)_z=(fg)\subset\ko_{X,z}}, and is generated by
\[\xymatrix@R=5pt{(fg)\ar[r] & \ko_{X,z}\\ \alpha.fg\fmaps[r] & \alpha}\]
It is clear that \ \m{\imm(\phi_{\ko_Y,z})=g.\ko_X(Z)_z}. The lemma follows 
immediately.
\end{proof}

\sepprop

\begin{subsub}\label{lem1f}{\bf Lemma: } Let $\E$ be a vector bundle on $X$ 
and $F$ a vector bundle on $Y$. Then there is a canonical isomorphism
\[\Ext^1_{\ko_X}(F,\E) \ \simeq \ \Hom(F,\E(Y)_{|Y}) \ . \]
\end{subsub}

\begin{proof} We have \ \m{\HHom(F,\E)=0}, hence, from the Ext spectral 
sequence (cf. \cite{go}, 7.3) it suffices to prove that \ 
\m{\EExt^1_{\ko_X}(F,\E)\simeq\HHom(F,\E(Y)_{|Y})}.

Let \m{r=\rk(F)}. Locally we have a locally free resolution of $F$:
\xmat{0\ar[r] & \E_1\ar[r]\fleq[d] & \E_0\ar[r]\fleq[d] & F\ar[r] & 0\\
 & \E_0(-Y) & r\ko_X}
hence we have an exact sequence
\[0\lra\HHom(\E_0,\E)\lra\HHom(\E_0(-Y),\E)\lra\EExt^1_{\ko_X}(F,\E)\lra 0 \ . 
\]
The result follows immediately.
\end{proof}

\end{sub}

\sepsub

\Ssect{A lemma on extensions}{le_ex}

Let $X$ be an algebraic variety, and $A$, $B$, $C$, $D$, $\ke$, coherent 
sheaves, such that we have a commutative diagram with exact rows
\xmat{0\ar[r] & C\ar[r]\ar[d]^\alpha & \ke\ar[r]\fleq[d] & D\ar[r]\ar[d]^\beta 
& 0 \\ 0\ar[r] & A\ar[r] & \ke\ar[r] & B\ar[r] & 0}
Suppose that $\alpha$ is injective, and let \ \m{N=\coker(\alpha)}. Then we 
have a canonical isomorphism \ \m{\ker(\beta)\simeq N}. Let \ \m{i:N\to D} be 
the inclusion.

Let \m{q\geq 0} be an integer and
\[\delta:H^q(N)\lra H^{q+1}(C) \ , \qquad \delta':H^q(D)\to H^{q+1}(C)\]
be the mappings coming from the exact sequences
\[0\lra C\lra A\lra N\lra 0 \quad \text{and} \quad 0\lra C\lra\ke \lra D\lra 0\]
respectively. The following lemma is easily proved (by using \v Cech 
cohomology):

\sepprop

\begin{subsub}\label{lem1e}{\bf Lemma: } The following triangle is commutative
\xmat{H^q(N)\ar[rr]^\delta\ar[dd]^{H^q(i)} & & H^{q+1}(C)\\ \\
H^q(D)\ar[uurr]_{\delta'}}
\end{subsub}

\end{sub}

\sepsub

\Ssect{Generalised elementary modifications of vector bundles on curves}{el_mod}

(cf. \cite{ma3})

Let $C$ be a smooth algebraic curve, \m{x\in C}, $E$ a vector bundle on $C$ and 
\ \m{N\subset E_x} \ a linear subspace. If $W$ is a finite dimensional vector 
space, let \m{W_x} denote the skyscraper sheaf concentrated on $x$ with fibre 
$W$ at $x$. Let \m{E_N} denote the kernel of the canonical morphism
\[E\lra (E_x/N)_x \ . \]
It is a vector bundle. We call \m{E_N} the {\em elementary modification of $E$} 
defined by $N$ (elementary modifications are well known for rank 2 vector
bundles).

We can also define elementary deformations from a finite set of distinct points 
of $C$: \m{x_1,\ldots,x_n}, and for \m{1\leq i\leq n}, a linear subspace 
\m{N_i\subset E_{x_i}}. We obtain the vector bundle \m{E_{N_1,\ldots,N_n}}, 
which is also obtained by successive elementary transformations involving one 
point, i.e. we define \m{E_{N_1,\ldots,N_i}} inductively by the relation \ 
\m{E_{N_1,\ldots,N_{i+1}}=(E_{N_1,\ldots,N_i})_{N_{i+1}}}.

\sepprop

\begin{subsub}\label{perio} Periodicity -- \rm
We have \ \m{E(-x)\subset E_N}, and \m{E(-x)} is itself an elementary
modification of \m{E_N}. We have a canonical exact sequence
\[0\lra E(-x)\lra E_N\lra N_x\lra 0 \ . \]
The kernel of the induced map \ \m{E_{N,x}\to N} \ is canonically isomorphic to
\Nligne\m{E(-x)_x/(N\ot\ko_C(-x)_x))}. i.e we have
\[E(-x) \ = \ (E_N)_{E(-x)_x/(N\ot\ko_C(-x)_x)}\]
(this can be seen for example by taking a local trivialisation of $E$).
\end{subsub}

\sepprop

\begin{subsub}\label{GEM_dual} Duality -- \rm We can see \m{(E_x/N)^*} as a 
linear subspace of \m{E_x^*}, and we have
\[(E^*)_{(E_x/N)^*} \ = \ (E_N)^*(-x) \ . \]
\end{subsub}

\end{sub}

\sepsub

\Ssect{The Ext spectral sequence}{koda0}

(cf. \cite{go})

Let $\ke$, $\kf$ be coherent sheaves on an algebraic variety $X$. Then there 
exists a spectral sequence such that
\[E_2^{pq} \ = \ H^p(X,\EExt^q(\ke,\kf)) \ \Longrightarrow \ 
\Ext^{p+q}(\ke,\kf) \ . \]
It follows that we have an exact sequence
\xmat{0\ar[r] & H^1(X,\HHom(\ke,\kf))\ar[r]^-\tau & 
\Ext^1_{\ko_X}(\ke,\kf)\ar[r] & H^0(X,\EExt^1_{\ko_X}(\ke,\kf))\ar[r] & 
H^2(X,\HHom(\ke,\kf)).}
In particular, if \ \m{H^2(X,\HHom(\ke,\kf))=\nsp}, for example if 
\m{\dim(X)=1}, then we have an exact sequence
\xmat{0\ar[r] & H^1(X,\HHom(\ke,\kf))\ar[r]^-\tau & 
\Ext^1_{\ko_X}(\ke,\kf)\ar[r] & H^0(X,\EExt^1_{\ko_X}(\ke,\kf))\ar[r] & 0.}

The morphism $\tau$ can be described using \v Cech cohomology: let \m{(U_i)} be 
an open cover of $X$. Let \m{\alpha\in H^1(X,\HHom(\ke,\kf))} corresponding to 
a cocycle \m{(\alpha_{ij})}, where \ 
\m{\alpha_{ij}\in\Hom(\ke_{|U_{ij}},\kf_{|U_{ij}})}. Then \m{\tau(\alpha)} 
corresponds to the extension \ \m{0\to\kf\to\boldsymbol{\kg}\to\ke\to 0}, where 
$\boldsymbol{\kg}$ is obtained by gluing the sheaves \m{(\kf\oplus\ke)_{|U_i}} 
with the automorphisms of \m{(\kf\oplus\ke)_{|U_{ij}}} defined by the matrices 
\m{\begin{pmatrix} I_\kf & \alpha_{ij}\\ 0 & I_\ke\end{pmatrix}}.

Let \m{Y\subset X} be a closed subvariety. Then if  
\begin{equation}\label{equ15}0\lra\kf\lra\boldsymbol{\kg}\lra\ke\lra 0
\end{equation}
is an extension coming from an element of \m{\imm(\tau)}, the restriction 
\begin{equation}\label{equ16}0\lra\kf_{|Y}\lra\boldsymbol{\kg}_{|Y}\lra\ke_{|Y}
\lra 0\end{equation}
is exact, because locally on $X$, the equation \m{(\ref{equ15})} is split 
The map
\[H^1(X,\EEnd(E))\lra H^1(X,\EEnd(\ke_{|Y}))\]
induced by the canonical morphism \ \m{\EEnd(E)\to\EEnd(\ke_{|Y})} sends the 
element corresponding to \REF{equ16} to the one corresponding to \REF{equ15}.
\end{sub}

\sepsub

\Ssect{Infinitesimal deformations of coherent sheaves}{Koda1}

\begin{subsub}\label{def_sh} Deformations of sheaves -- \rm Let $X$ be a 
projective algebraic variety and $E$ a coherent sheaf on $X$. A {\em 
deformation} of $E$ is a quadruplet \ \m{\kd=(S,s_0,\ke,\alpha)}, where 
\m{(S,s_0)} is the germ of an analytic variety, $\ke$ is a coherent sheaf on 
\m{S\times X}, flat on $S$, and $\alpha$ an isomorphism \m{\ke_{s_0}\simeq E}.
If there is no risk of confusion, we also say that $\ke$ is an infinitesimal 
deformation of $E$. Let \m{Z_2=\spec(\C|t]/(t^2))}. When \m{S=Z_2} and \m{s_0} 
is the closed point $*$ of \m{Z_2}, we say that $\kd$ is an {\em infinitesimal 
deformation} of $E$. Isomorphisms of deformations of $E$ are defined in an 
obvious way. If \ \m{f:(S',s'_0)\to(S,s_0)} \ is a morphism of germs, the 
deformation \m{f^\#(\kd)} is defined as well. A deformation \ 
\m{\kd=(S,s_0,\ke,\alpha)} \ is called {\em semi-universal} if for every 
deformation \ \m{\kd'=(S',s'_0,\ke',\alpha')} \ of $E$, there exists a morphism 
\ \m{f:(S',s'_0)\to(S,s_0)} \ such that \ \m{f^\#(\kd)\simeq\kd'}, and if the 
tangent map \ \m{T_{s'_0}S'\to T_{s_0}S} \ is uniquely determined. There always 
exists a semi-universal deformation of $E$ (cf. \cite{si_tr}, theorem I).

Let $\ke$ be an infinitesimal deformation of $E$. Let \m{p_X} denote the 
projection \ \m{Z_2\times X\to X}. Then there is a canonical exact sequence
\[0\lra E\lra p_{X*}(\ke)\lra E\lra 0 \ , \]
i.e. an extension of $E$ by itself. In fact, by associating this extension to 
$\ke$ one defines a bijection between the set of isomorphism classes of 
infinitesimal deformations of $E$ and the set of isomorphism classes of 
extensions of $E$ by itself, i.e. \m{\Ext^1_{\ko_X}(E,E)}.
\end{subsub}

\sepprop

\begin{subsub}\label{ko_da} Koda\"\i ra-Spencer morphism -- \rm Let \ 
\m{\kd=(S,s_0,\ke,\alpha)} \ be a deformation of $E$, and \m{X_{s_0}^{(2)}} the 
infinitesimal neighbourhood of order 2 of \ \m{X_{s_0}=\{s_0\}\times X} \ in 
\m{S\times X}. Then we have an exact sequence on \m{X_{s_0}^{(2)}}
\[0\lra T_{s_0}S\ot E\lra\ke/m_s^2\ke=\ke_{|X_{s_0}^{(2)}}\lra E\lra 0 \ . \]
By taking the direct image by \m{p_X} we obtain the exact sequence on $X$
\[0\lra T_{s_0}S\ot E\lra p_{X*}(\ke/m_s^2\ke)\lra E\lra 0 \ , \]
hence a linear map
\[\omega_{s_0}:T_{s_0}S\lra\Ext^1_{\ko_X}(E,E) \ , \]
which is called the {\em Koda\"\i ra-Spencer morphism} of $\ke$ at \m{s_0}.

We say that $\ke$ is a {\em complete deformation} if \m{\omega_{s_0}} is 
surjective. If $\kd$ is a semi-universal deformation, \m{\omega_{s_0}} is an 
isomorphism.
\end{subsub}

\end{sub}

\sepsub

\Ssect{The Koda\"\i ra-Spencer map for vector bundles on families of 
curves}{Koda}

\begin{subsub}\label{Koda_1} Vector bundles on families of smooth curves --
\rm A {\em family of smooth curves} parametrised by an algebraic variety $U$
is a flat morphism \ \m{\rho:\kx\to U} \ such that for every closed point
\m{x\in U}, \m{\kx_x=\rho^{-1}(x)} is a smooth projective connected curve. 
If \ \m{f:Y\to U} \ is a morphism of algebraic varieties,
\m{f^*(\kx)=\kx\times_UY} \ is a family of smooth curves parametrised by $Y$.
If $\E$ is a vector bundle on $\kx$, let \ \m{f^\#(\E)=p_\kx^*(\E)} (where
\m{p_\kx} is the projection \ \m{\kx\times_UY\to\kx}).
\end{subsub}

\sepprop

Let $S$ be a smooth curve, \m{s_0\in S}, \m{\pi:\kx\to S} \ a flat family of 
smooth projective curves and \ \m{C=\kx_{s_0}=\pi^{-1}(s_0)}. When dealing with 
morphisms of smooth curves \m{\phi:T\to S}, we will always assume that there is
exactly one closed point of $T$ over \m{s_0} (if $\phi$ is not constant and 
\m{s_0\in\imm(\phi)}, it is always possible to remove from $T$ some points in 
the finite set \m{\phi^{-1}(s_0)}).

Let $T$, \m{T'} be smooth curves, \m{t_0\in T}, \m{t'_0\in T'}, 
\m{\alpha:T\to S}, \m{\alpha':T'\to S} morphisms such that \m{\alpha(t_0)=s_0},
\m{\alpha'(t'_0)=s_0}. Suppose that the tangent maps \
\m{T_{t_0}\alpha:T_{t_0}T\to T_{s_0}S},\Nligne \m{T_{t_0}\alpha':T_{t'_0}T'\to
T_{s_0}S} \ are injective. Let \ \m{Z=T\times_ST'}, \m{p_T}, \m{p_{T'}}, 
$\beta$ the projections \ \m{Z\to T}, \m{Z\to T'}, \m{Z\to S} \ respectively. 
Let \ \m{\ky=\beta^*(\kx)}. Let \ \m{z_0=(t_0,t'_0)\in Z}.

\sepprop

\begin{subsub}\label{Koda_2} Semi-universal families -- \rm
Let $E$ be a vector bundle on $C$.
There exists a smooth connected variety $R$, a flat morphism \m{\eta:R\to Z},
a vector bundle $\E$ on \m{\eta^*(\ky)} such that:
\begin{enumerate}
\item[--] There exists $r_0\in R_{z_0}$ such that \ $\E_{r_0}\simeq E$.
\item[--] For every $z\in Z$, $\E_{|R_z\times\kx_{\beta(z)}}$ is a complete
deformation, i.e. for every \m{r\in R_z}, the Koda\"\i ra-Spencer map
\[\omega_r:(TR_z)_r\lra\Ext^1_{\ko_{\kx_{\beta(z)}}}(\E_r,\E_r)\]
is surjective.
\item[--] $R$ has the following local universal property: for every \m{z\in
Z}, every $r\in R_z$ every neighbourhood $U$ of $z$, every vector bundle
$\kf$ on $\beta^{-1}(U)$ such that \ $\kf_z\simeq\E_r$, there exists a
neighbourhood $U'\subset U$ of $z$ and a morphism \ $\phi:U'\to R$ \ such that
\ $\phi(z)=r$, $\eta\circ\phi=I$ \ and \ $\phi^\#(\ke)\simeq\kf$.
\end{enumerate}
(for example $R$ can be constructed by using relative Quot schemes). A useful 
consequence of the existence of $R$ is that {\em there exists a smooth curve 
\m{S'}, \m{s'_0\in S'}, a morphism \m{\phi:S'\to S} such that 
\m{\phi(s'_0)=s_0}, and a vector bundle $\ke$ on \m{\phi^*(\kx)} such that \ 
\m{\ke_{s'_0}\simeq E}} (one can take for \m{S'} an appropriate curve in $R$ 
through \m{r_0}).
\end{subsub}

\sepprop

\begin{subsub}\label{Koda_3} The Koda\"\i ra-Spencer map -- \rm Let $T$, \m{T'} 
be smooth curves, \m{t_0\in T}, \m{t'_0\in T'}, 
\m{\alpha:T\to S}, \m{\alpha':T'\to S} morphisms such that \m{\alpha(t_0)=s_0},
\m{\alpha'(t'_0)=s_0}. Suppose that the tangent maps \Nligne
\m{T_{t_0}\alpha:T_{t_0}T\to T_{s_0}S}, \m{T_{t_0}\alpha':T_{t'_0}T'\to
T_{s_0}S} \ are injective. Let \ \m{Z=T\times_ST'}, \m{p_T}, \m{p_{T'}}, 
$\beta$ the projections \ \m{Z\to T}, \m{Z\to T'}, \m{Z\to S} \ respectively. 
Let \ \m{\ky=\beta^*(\kx)}. Let \ \m{z_0=(t_0,t'_0)\in Z}.
\rm Let $E$ be a
vector bundle on $C$, $\ke$ a vector bundle on \m{\alpha^*(\kx)}, \m{\ke'} a
vector bundle on \m{\alpha'^*(\kx)}, such that there are isomorphisms \
\m{\ke_{s_0}\simeq E}, \m{\ke'_{s_0}\simeq E}. We have two vector bundles
\m{p_T^\#(\ke)}, \m{p_{T'}^\#(\ke')} on \m{\beta^*(\kx)}. There exist a
neighbourhood $U$ of \m{z_0} and morphisms \ \m{f:U\to R}, \m{f':U\to R} \ such
that \ \m{\eta\circ f=\eta\circ f'=I} \ and \ \m{f^\#(\E)\simeq p_T^\#(\ke)},
\m{{f'}^\#(\E)\simeq p_{T'}^\#(\ke')}. We have
\[T\eta_{r_0}\circ Tf_{z_0} \ = \ T\eta_{r_0}\circ Tf'_{z_0} \ = \ I_{TZ_{z_0}}
\ ,\]
hence \ \m{\imm(Tf_{z_0}-Tf'_{z_0})\subset TR_{z_0,r_0}}. Now we define the
linear map
\[\omega_{\ke,\ke'}:TS_{s_0}\lra\Ext^1_{\ko_C}(E,E)\]
by:
\[\omega_{\ke,\ke'} \ = \
\omega_{r_0}\circ(Tf_{z_0}-Tf'_{z_0})\circ(T\beta_{z_0})^{-1}.\]
It is easy to see that this definition is independent of the choice of $R$
satisfying the above properties.

If $\kx$ is the trivial family, i.e. \m{\kx=C\times S} and \m{T=T'=S}, then
there is a canonical choice for \m{\ke'}: \m{\ke_0=p_C^*(E)} (where
\m{p_C:C\times S\to C} is the projection). We can then define
\[\omega_\ke \ = \ \omega_{\ke,\ke_0}:TS_{s_0}\lra\Ext^1_{\ko_C}(E,E) \ , \]
This is the usual Koda\"\i ra-Spencer morphism.
\end{subsub}

\sepprop

\begin{subsub}\label{Koda_4} Other equivalent definition -- \rm It uses the
properties of coherent sheaves on primitive double curves (cf.
\ref{coh_prim}). We assume for simplicity that \ \m{T=T'=Z}. Let \m{{\mathfrak
m}_{z_0}\subset\ko_{Z,z_0}}\ be the maximal ideal, and
\[Y \ = \ C^{(2)} \ = \ \pi^{-1}(\spec(\ko_{S,s_0}/{\mathfrak m}_{s_0}^2)) \ = \
\beta^{-1}(\spec(\ko_{Z,z_0}/{\mathfrak m}_{z_0}^2))
\]
the second infinitesimal neighbourhood of $C$ in $\kx$ or \m{\beta^*(\kx)}, 
which is a primitive double curve with associated smooth curve $C$ and 
associated line bundle 
\[\ki_{C}/\ki_{C}^2 \ \simeq TS_{s_0}^*\ot\ko_C \ \simeq \ \ko_C \ . \]
(the isomorphism \m{\ki_{C}/\ki_{C}^2\simeq\ko_C} depends on the choice of a 
generator of \m{T_{s_0}S}). We have an exact sequence
\[0\lra E\ot TZ_{z_0}^*\lra\ke_{|Y}\lra E\lra 0\]
(given by the canonical filtration of \m{\ke_{|Y}}). Let \ 
\m{\sigma_\ke\in\Ext^1_{\ko_Y}(E\ot TZ_{z_0}^*,E)} \ be the element 
associated to this exact sequence. The associated morphism \ \m{E\ot 
TZ_{z_0}^*\to E\ot TZ_{z_0}^*} (cf. \ref{cstr}) is of course the identity.
We define similarly \ \m{\sigma_{\ke'}\in\Ext^1_{\ko_Y}(E\ot TZ_{z_0}^*,E)}.
Then we have
\[\omega_{\ke,\ke'} \ = \ (\sigma_\ke-\sigma_{\ke'})\circ(T\beta_{z_0})^{-1} \ 
. \]
\end{subsub}

\sepprop

\begin{subsub}\label{Koda_5} Properties of the Koda\"\i ra-Spencer map -- \rm
Let $W$, \m{W'} be smooth curves, \m{w_0\in W}, \m{w'_0\in W'}, \m{\lambda:W\to 
T}, \m{\lambda':W'\to T'} morphisms such that \m{\lambda(w_0)=t_0} and 
\m{\lambda'(w'_0)=t'_0}, and such that the tangent maps \m{T\lambda_{w_0}} and
\m{T\lambda_{w'_0}} are isomorphisms. Then we have
\[\omega_{\lambda^\#(\ke),{\lambda'}^\#(\ke')} \ = \ \omega_{\ke,\ke'} \ . \]
\end{subsub}

\sepprop

\begin{subsub}\label{prop1b}{\bf Proposition: } Let $u$ be a non zero element 
of \m{TS_{s_0}}. Let \ \m{\sigma\in\Ext^1_{\ko_C}(E,E)}. Then there exists a 
smooth curve $\bf T$, \m{{\bf t}_0\in T}, a morphism 
\m{\boldsymbol{\alpha}:{\bf T}\to S} such that \m{\boldsymbol{\alpha}({\bf 
t}_0)=s_0}, and a vector bundle $\boldsymbol{\ke}$ on 
\m{\boldsymbol{\alpha}^*(\kx)} such that \m{\boldsymbol{\ke}_{{\bf t}_0}\simeq 
E} and that
\[\omega_{\boldsymbol{\ke},\ke}(u) \ = \ \sigma \ . \]
\end{subsub}

\begin{proof} It suffices to take appropriate curves in $R$. \end{proof}

\end{sub}

%\sepsec
\newpage

\section{Coherent sheaves on primitive double curves}\label{coh_prim}

\Ssect{Primitive double curves}{pdc}

(cf. \cite{ba_fo}, \cite{ba_ei}, \cite{dr2}, \cite{dr1}, \cite{dr4},
\cite{dr5},
\cite{dr6}, \cite{ei_gr}).

Let $C$ be a smooth connected projective curve. A {\em multiple curve with
support $C$} is a Cohen-Macaulay scheme $Y$ such that \m{Y_{red}=C}.

Let $n$ be the smallest integer such that \m{Y\subset C^{(n-1)}}, \m{C^{(k-1)}}
being the $k$-th infinitesimal neighbourhood of $C$, i.e. \
\m{\ki_{C^{(k-1)}}=\ki_C^{k}} . We have a filtration \ \m{C=C_1\subset
C_2\subset\cdots\subset C_{n}=Y} \ where $C_i$ is the biggest Cohen-Macaulay
subscheme contained in \m{Y\cap C^{(i-1)}}. We call $n$ the {\em multiplicity}
of $Y$.

We say that $Y$ is {\em primitive} if, for every closed point $x$ of $C$,
there exists a smooth surface $S$, containing a neighbourhood of $x$ in $Y$ as
a locally closed subvariety. In this case, \m{L=\ki_C/\ki_{C_2}} is a line
bundle on $C$ and we have \ \m{\ki_{C_i}=\ki_C^i},
\m{\ki_{C_{i}}/\ki_{C_{i+1}}=L^i} \ for \m{1\leq i<n}. We call $L$ the line
bundle on $C$ {\em associated} to $Y$. Let \m{P\in C}. Then there exist
elements $z$, $x$ of \m{m_{S,P}} (the maximal ideal of \m{\ko_{S,P}}) whose
images in \m{m_{S,P}/m_{S,P}^2} form a basis, and such that for \m{1\leq i<n}
we have \ \m{\ki_{C_i,P}=(z^{i})}. The image of $x$ in \m{\ko_{C,P}} is then a
generator of the maximal ideal.

The simplest case is when $Y$ is contained in a smooth surface $S$. Suppose
that $Y$ has multiplicity $n$. Let \m{P\in C} and \m{f\in\ko_{S,P}}  a local
equation of $C$. Then we have \ \m{\ki_{C_i,P}=(f^{i})} \ for \m{1<j\leq n},
in particular \m{I_{Y,P}=(f^n)}, and \ \m{L=\ko_C(-C)}.

For any \m{L\in\Pic(C)}, the {\em trivial primitive curve} of multiplicity $n$,
with induced smooth curve $C$ and associated line bundle $L$ on $C$ is the
$n$-th infinitesimal neighbourhood of $C$, embedded by the zero section in the
dual bundle $L^*$, seen as a surface.

In this paper we are interested only in the case \m{n=2}. Primitive multiple 
curves of multiplicity 2 are called {\em primitive double curves} or {\em 
ribbons}.

\end{sub}

\sepsub

\Ssect{Canonical filtrations and invariants of coherent sheaves}{cspdc}

Let $Y$ be a ribbon, \m{C=Y_{\text{red}}}, and $L$ the associated line bundle
on $C$. Let \m{P\in C} be a closed point, \m{z\in\ko_{Y,P}} an equation of $C$ 
and $M$ a \m{\ko_{Y,P}}-module of finite type. Let $\ke$ be a coherent sheaf on 
$Y$. We denote by \ \m{\tau_\ke:\ke\ot L\to\ke} \ the canonical morphism.

\sepprop

\begin{subsub}\label{QLL-def1}First canonical filtration -- \rm
The {\em first canonical filtration of $M$} is
\[\nsp\subset G_1(M)\subset M\]
where \m{G_1(M)} is the kernel of the canonical surjective morphism \ 
\m{M\to M\ot_{\ko_{Y,P}}\ko_{C,P}} .
So we have \ \m{G_1(M)=zM}. We will note \ \m{G_0(M)=M/G_1(M)}. The direct sum
\[{\rm Gr}(M) \ = \ G_1(M)\oplus G_0(M)\]
is a \m{\ko_{C,P}}-module. The following properties are obvious: \Nligne
-- we have $G_1(M)=\nsp$ if and only if $M$ is a $\ko_{C,P}$-module,\Nligne
-- a morphism of \m{\ko_{Y,P}}-modules sends the first canonical filtration of 
the first module to the one of the second. 

We define similarly the {\em first canonical filtration of $\ke$}:
\[0\subset G_1(\ke)\subset\ke\]
where \m{G_1(\ke)} is the kernel of the canonical surjective morphisms \ 
\m{\ke\to\ke_{\mid C}}. We will note \ \m{G_0(\ke)=\ke/G_1(\ke)=\ke_{|C}}. We
will also use the notation \m{\ke_1=G_1(\ke)}. Let
\[{\rm Gr}(\ke) \ = \ G_1(\ke)\oplus G_0(\ke) \ . \]
It is a \m{\ko_{C}}-module. The following properties are obvious: \Nligne
- we have \m{G_1(\ke)=\ki_C\ke},\Nligne
- we have $G_1(\ke)=0$ if and only if $\ke$ is a sheaf on $C$,\Nligne
- every morphism of coherent sheaves on \m{Y} sends the first canonical 
filtration of the first to that of the second.
\end{subsub}

\sepsubsub

\begin{subsub}\label{2-fc}Second canonical filtration -- \rm
The {\em second canonical filtration of $M$} is
\[0\subset G^{(1)}(M)\subset M \ , \]
with \ \m{G^{(1)}(M)=\big\lbrace u\in M ; z^iu=0\big\rbrace}. We have \ 
\m{G_1(M)\subset G^{(1)}(M)}. We will note \Nligne \m{G^{(2)}(M)=M/G^{(1)}(M)}. 
The direct sum
\[\Gr_2(M) \ = \ G^{(1)}(M)\oplus G^{(2)}(M) \]
is a \m{\ko_{C,P}}-module. Every morphism of \m{\ko_{Y,P}}-modules sends the 
second canonical filtration of the first module to that of the second. 

We define similarly the {\em second canonical filtration of $\ke$}:
\[0\subset G^{(1)}(\ke)\subset \ \ke \ , \qquad G^{(2)}(\ke) \ = \ 
\ke/G^{(1)}(\ke) \ . \]
We will also use the notation \ \m{\ke^{(1)}=G^{(1)}(\ke)}. The direct sum
\[\Gr_2(\ke) \ = \ G^{(1)}(\ke)\oplus G^{(2)}(\ke) \]
is a \m{\ko_{C}}-module. We have \ \m{G_1(\ke)\subset G^{(1)}(\ke)}, and every
morphism of coherent sheaves on $Y$ sends the second canonical filtration of 
the first to that of the second.
\end{subsub}

\sepprop

\begin{subsub}\label{rel_filt} Relations between the two filtrations - \rm
we have a canonical isomorphism
\[G_1(\ke) \ \simeq \ G^{(2)}(\ke)\ot L \ . \]
The canonical morphism \m{\ke\ot\L\to\ke} induces two morphisms of sheaves on 
$C$:
\[\lambda_\ke=:G^{(2)}(\ke)\ot L\lra G^{(1)}(\ke) \ , \qquad 
\mu_\ke=:G_0(\ke)\ot L\lra G_1(\ke) \ . \]
The first morphism of sheaves is injective, the second is surjective, and we
have
\[\ker(\mu_\ke) \ \simeq \ \coker(\lambda_\ke)\ot L \ . \]

\end{subsub}

\sepprop

\begin{subsub}\label{cha_can_fil} Characterisation of the canonical filtrations 
-- \rm (cf. lemma \ref{lem1g}). Let
\[0\subset E\subset\ke\]
be a filtration such that $E$ is a sheaf on $C$. Then we have \ \m{E\subset 
G^{(1)}(\ke)}, and \m{F=\ke/E} is a sheaf on $C$ if and only if \
\m{G_1(\ke)\subset E}. Suppose that \ \m{G_1(\ke)\subset E}. Then the canonical 
morphism \m{\tau_\ke} induces \ 
\m{\phi:F\ot L\to E}, and we have
\begin{enumerate}
\item[--] $E=G_1(\ke)$ if and only if $\phi$ is surjective.
\item[--] $E=G^{(1)}(\ke)$ if and only $\phi$ is injective.
\end{enumerate}

\end{subsub}

\sepprop

\begin{subsub}\label{inva} Invariants -- \rm The integer \
\m{R(M)=\rk({\rm Gr}(M))} \ is called the {\em generalised rank} of $M$.

The integer \ \m{R(\ke)=\rk({\rm Gr}(\ke))} \ is called the {\em generalised 
rank} of $\ke$. So we have \m{R(\ke)=R(\ke_P)} for every $P$ in a nonempty open
subset of $C$.

The integer \ \m{\Deg(\ke)=\deg({\rm Gr}(\ke))} \ is called the {\em 
generalised degree} of $\ke$. The generalised rank and degree of sheaves are 
deformation invariants.

If \m{R(\ke)>0}, let \ \m{\mu(\ke)=\Deg(\ke)/R(\ke)}. We call this rational 
number the {\em slope} of $\ke$.

Let \ \m{0\to M'\to M\to M''\to 0} \ be an exact sequence of 
\m{\ko_{Y,P}}-modules of finite type. Then we have \ \m{R(M)=R(M')+R(M'')} .

Let \ \m{0\to\ke\to\kf\to\kg\to 0} \ be an exact sequence of coherent sheaves 
on $Y$. Then we have \ \m{R(\kf)=R(\ke)+R(\kg)} , 
\m{\Deg(\kf)=\Deg(\ke)+\Deg(\kg)}.

If $\ke$ is a coherent sheaf on $Y$, we have
\[\chi(\ke) \ = \ \Deg(\ke)+R(\ke)(1-g) \ , \]
where $g$ is the genus of $C$ (Riemann-Roch theorem). It follows that if
\m{\ko_Y(1)} is an ample line bundle on $Y$, \m{\ko_C(1)=\ko_Y(1)_{|C}} \ and \
\m{\delta=\deg(\ko_C(1))}, the Hilbert polynomial of $\ke$ with respect to
\m{\ko_Y(1)} is given by
\[P_{\ke}(m) \ = \ R(\ke)\delta m+\Deg(\ke)+R(\ke)(1-g) \ . \]
\end{subsub}
\end{sub}

\sepsub

\Ssect{Quasi locally free sheaves -- Reflexive sheaves}{QLF-R}

\begin{subsub}\label{tors}Torsion -- \rm 
Let $M$ be a \m{\ko_{Y,P}}-module of finite type. Le {\em torsion sub-module
\m{T(M)} of $M$} consists of the elements $u$ of $M$ such that there exists an 
integer \m{p>0} such that \m{x^pu=0} (recall that \m{x\in\ko_{Y,P}} is such 
that its image in \m{\ko_{C,P}} is a generator of the maximal ideal). We say 
that $M$ is {\em torsion free} if 
\ \m{T(M)=\nsp}.

Let $\ke$ be a coherent sheaf on $Y$. The {\em torsion subsheaf \m{T(\ke)} of 
$\ke$} is the maximal subsheaf of $\ke$ with finite support. For every closed 
point $P$ of $C$ we have \ \m{T(\ke)_P=T(\ke_P)}. We say that $\ke$ is a {\em 
torsion sheaf} if \ \m{\ke=T(\ke)} , or equivalently, if its support is finite.
We say that $\ke$ is {\em torsion free} if \m{T(\ke)=0}.

The following conditions are equivalent:
\begin{enumerate}
\item[(i)] $\ke$ is torsion free.
\item[(ii)] $G^{(1)}(\ke)$ is locally free.
\end{enumerate}
Moreover, if $\ke$ is torsion free, $G^{(2)}(\ke)$ is also locally free.
\end{subsub}

\sepprop

\begin{subsub}\label{QLL}Quasi locally free sheaves -- \rm
Let $M$ be a \m{\ko_{Y,P}}-module of finite type. We say that $M$ is 
{\em 
quasi free} if there exist integers \m{m_1,m_2\geq 0} such that \ 
\m{M\simeq m_1\ko_{C,P}\oplus m_2\ko_{Y,P}}. These integers are uniquely 
determined. In this case we say that $M$ is {\em of type} \m{(m_1,m_2)}. We 
have \ \m{R(M)=m_1+2m_2}.

Let $\ke$ be a coherent sheaf $Y$. We say that $\ke$ is {\em quasi locally 
free} at a point $P$ of $C$ if there exists a neighbourhood $U$ of $P$ and 
integers \m{m_1,m_2\geq 0} such that for every \m{Q\in U}, \m{\ke_Q} is quasi 
free of type \m{m_1,m_2}. The integers \m{m_1,m_2} are uniquely
determined and depend only of $\ke$, and \m{(m_1,m_2)} is called the 
{\em type of} $\ke$.

We say that $\ke$ is {\em quasi locally free} if it is quasi locally free at 
every point of $C$.

The following conditions are equivalent: 
\begin{enumerate}
\item[(i)] $\ke_P$ is quasi locally free at $P$.
\item[(ii)] $G_0(\ke)$ and \m{G_1(\ke)} are free $\ko_{C,P}$-modules.
\item[(iii)] $\coker(\lambda_\ke)$ and $G_0(\ke)$ are free $\ko_{C,P}$-modules.
\end{enumerate}

The following conditions are equivalent: 
\begin{enumerate}
\item[(i)] $\ke$ is quasi locally free.
\item[(ii)] $G_0(\ke)$ and $G_1(\ke)$ are locally free on $C$.
\item[(iii)] $\coker(\lambda_\ke)$ and $G_0(\ke)$ are locally free on $C$.
\end{enumerate}

Assume that $\ke$ is torsion free. Then \m{\lambda_\ke} is a morphism of vector 
bundles, which is an injective morphism of sheaves, and $\ke$ is quasi locally 
free if and only if  \m{\lambda_\ke} is an injective morphism of vector bundles.
\end{subsub}

\sepprop

\begin{subsub}\label{refle} Reflexive sheaves -- \rm
Let \m{P\in C} and $M$ be a \m{\ko_{Y,P}}-module of finite type. Let \
\m{M^{\vee}} be the {\em dual} of $M$: \m{M^{\vee}=\Hom(M,\ko_{Y,P})}. If $N$ 
is \m{\ko_{C,P}}-module of finite type, we will denote \m{N^*} the dual of $N$: 
\m{N^*=\Hom(N,\ko_{C,P})}.

Let $\ke$ be a coherent sheaf on $Y$. Let \m{\ke^{\vee}} be the {\em
dual} of $\ke$: \ \m{\ke^{\vee}=\HHom(\ke,\ko_Y)}. If $E$ is a coherent sheaf 
on $C$, we will denote by \m{E^*} the dual of $E$: 
\m{E^*=\HHom(E,\ko_C)}. We have \ \m{E^\vee=E^*\ot L}.

We have \ \m{G^{(1)}(\ke^\vee)=G_0(\ke)^*\ot L}. If $\ke$ is quasi locally free 
we have
\[G_1(\ke^\vee) \ = \ G_1(\ke)^*\ot L^2 \ , \qquad G^{(2)}(\ke^\vee) \ = \
G_1(\ke)^*\ot L \ . \]

The following properties are equivalent (cf. \cite{dr2}, \cite{dr4}):
\begin{itemize}
\item[(i)] $\ke$ is reflexive.
\item[(ii)] $\ke$ is torsion free.
\item[(iii)] We have \ $\EExt^1_{\ko_Y}(\ke,\ko_Y)=0$ .
\end{itemize}
Moreover we have in this case \
\m{\EExt^i_{\ko_Y}(\ke,\ko_Y)=0} \ for every \m{i\geq 1}.

It follows that if \ \m{0\to\ke\to\kf\to\kg\to 0} is an exact sequence of 
coherent torsion free sheaves on $Y$, the dual sequence \
\m{0\to\kg^\vee\to\kf^\vee\to\ke^\vee\to 0} \ is also exact.
\end{subsub}

\sepprop

\begin{subsub}\label{mod} Moduli spaces of semi-stable reflexive sheaves -- \rm
A reflexive sheaf $\ke$ is semi-stable (resp. stable) (as per M.~Maruyama
\cite{ma2} or C.~Simpson \cite{si}) if for every proper subsheaf
\m{\kf\subset\ke} we have
\[\mu(\kf) \ \leq \ \mu(\ke) \quad \text{(resp. \ } < \text{)} \ . \]
As for smooth curves, the definition of semi-stability does not depend on the
choice of an ample line bundle on $Y$. If \m{R>0} and $D$ are integers, we will
denote by \m{\km(R,D)} the moduli space of semi-stable sheaves on $Y$ of rank
$R$ and degree $D$, and by \m{\km_s(R,D)} the open subset corresponding to
stable sheaves.
\end{subsub}

\sepprop

\begin{subsub}\label{rig_t} Quasi locally free sheaves of rigid type -- \rm
(cf. \cite{dr4}, \cite{dr5}) A quasi locally free sheaf $\ke$ on $Y$ is called
{\em of
rigid type} if it is locally free, or locally isomorphic to
\m{a\ko_Y\oplus\ko_C}, for some integer \m{a\geq 0}. In the latter case, we have
\ \m{\rk(\ke_{|C})=a+1} \ and \ \m{\rk(\ke_1)=a}. The deformations of quasi
locally free sheaves of rigid type are quasi locally free sheaves of rigid type,
and $a$, \m{\deg(\ke_1)} and \m{\deg(\ke_{|C})} are invariant by deformation. 

Let \m{R=2a+1} and \m{d_0}, \m{d_1}, $D$ integers such that \ \m{D=d_0+d_1}.
The stable sheaves $\ke$ of rigid type, of rank $R$ and such that
\m{\deg(\ke_1\ot L^*)=d_1} and \m{\deg(\ke_{|C})=d_0} form an open subset
\m{\kn(R,d_0,d_1)} of an irreducible component of \m{\km_s(R,D)}. If
\m{\kn(R,d_0,d_1)} is non empty, we have
\[\dim(\kn(R,d_0,d_1)) \ = \ 1-a(a+1)\deg(L)+(g-1)(2a^2+2a+1) \ . \]
In general \m{\kn(R,d_0,d_1)} is not reduced. The reduced subvariety
\m{\kn(R,d_0,d_1)_{\text{red}}} is smooth and its tangent sheaf at the point
corresponding to a sheaf $\kf$ is canonically isomorphic to \m{H^1(\EEnd(\kf))}.
We have also canonical isomorphisms
\[T\kn(R,d_0,d_1)_\kf/T\kn(R,d_0,d_1)_{\text{red},\kf} \ \simeq \
H^0(\EExt^1(\kf,\kf)) \ \simeq \ H^0(L^*)\]
(for the last isomorphism see \cite{dr2}, cor. 6.2.2). It is known that
\m{\kn(R,d_0,d_1)} is non empty if
\[\frac{\epsilon}{a+1} \ < \ \frac{\delta}{a} \ < \
\frac{\epsilon-\deg(L)}{a+1} \ . \]
\end{subsub}

\end{sub}

\sepsub

\Ssect{Properties and construction of coherent sheaves}{cstr}

If $E$ is a vector bundle on $C$, there exists a vector bundle $\E$ of $Y$ such 
that \ \m{\E_{|Y}\simeq E} (\cite{dr2}, th. 3.1.1). In particular there exists 
a line bundle $\L$ on $Y$ such that \ \m{\L_{|C}\simeq L}. It follows that we 
have a locally free resolution of $E$ :
\[\cdots\E\ot\L^3\lra\E\ot\L^2\lra\E\ot\L\lra\E\lra E\lra 0 \ . \]
The following results follows easily:

\sepprop

\begin{subsub}\label{lem1c}{\bf Lemma: } Let $F$ be a coherent sheaf on $C$ 
and $G$ a vector bundle on $C$. Then we have, for every integer \m{i\geq 0}
\[\Tor^i_{\ko_Y}(F,\ko_C) \ \simeq \ F\ot L^i \ , \qquad
\EExt^1_{\ko_Y}(G,F) \ \simeq \ \HHom(G\ot L^i,F) \ . \]
\end{subsub}

\sepprop

It follows that we have a canonical exact sequence
\xmat{0\ar[r] & \Ext^1_{\ko_C}(G,F)\ar[r]^-i\fleq[d] & 
\Ext^1_{\ko_Y}(G,F)\ar[r] & \Hom(G\ot L,F)\ar[r]\fleq[d] & 0 \ .\\ & 
H^1(\HHom(G,F)) & & H^0(\EExt^1_{\ko_Y}(G,F))}

\sepprop

\begin{subsub}\label{cst1b} \rm We can give an explicit construction of the 
morphism $i$ in the preceding exact sequence, using \v Cech cohomology: let 
\m{(U_i)} an open cover of $C$ and \ \m{\alpha\in\Ext^1_{\ko_C}(G,F)}, 
represented by a cocycle \m{(\alpha_{ij})}, where \ 
\m{\alpha_{ij}\in\Hom(G_{|U_{ij}},F_{|U_{ij}})}. Let \ 
\m{\sigma\in\Ext^1_{\ko_Y}(G,F)}, corresponding to an extension
\[0\lra F\lra\ke\lra G\lra 0 \ . \]
We can view \m{\alpha_{ij}} as an endomorphism of \m{\ke_{|U_{ij}}}. Then 
\m{\sigma+i(\alpha)} corresponds to the extension
\[0\lra F\lra\ke'\lra G\lra 0 \ , \]
where \m{\ke'} is the sheaf obtained by gluing the sheaves \m{\ke_{|U_{ij}}} 
using the automorphisms \Nligne \m{I+\alpha_{ij}\in\Aut(\ke_{|U_{ij}})}.
\end{subsub}

\sepsubsub

\begin{subsub}\label{lem1g}{\bf Lemma: } Suppose that $F$ is locally free. Let
\ \m{\sigma\in\Ext^1_{\ko_Y}(G,F)}, corresponding to an extension
\[0\lra F\lra\ke\lra G\lra 0 \ . \]
Then

{\bf 1 -- } We have \m{G=\ke_{|C}} and \m{F=\ke_1} if and only if 
\m{\theta(\sigma)} is surjective. In this case $\ke$ is quasi locally free.

{\bf 2 -- } We have \m{G=\ke^{(2)}} and \m{F=\ke^{(1)}} if and only 
\m{\theta(\sigma)} is injective. In this case $\ke$ is quasi locally free if 
and only if \m{\theta(\sigma)} is injective as a morphism of vector bundles.
\end{subsub}

\begin{proof}
{\bf 1-} follows from the factorisation of the canonical morphism \ 
\m{\ke_{|C}\ot L\to\ke_1}
\xmat{\ke_{|C}\ot L\flon[r] & G\ot L\flon[r]\ar@/^2pc/[rr]^-{\theta(\sigma)} & 
\ke_1\flinc[r] & F ,}
and {\bf 2-} from the factorisation of \m{\theta(\sigma)}
\xmat{G\ot L\flon[r] & \ke/\ke^{(1)}\ot L\fleq[r] & \ke_1\flinc[r] & F .}
\end{proof}

\sepprop

\begin{subsub}\label{lem1b}{\bf Lemma: } Let $\E$ be a vector bundle on $Y$, 
$E$ a coherent sheaf on $C$, \m{\phi:\E\to E} a surjective morphism and \ 
\m{\kn=\ker(\phi)}. Then we have an exact sequence
\xmat{0\ar[r] & E\ot L\ar[r] & \kn_{|C}\ar[r] & \E_{|C}\ar[r] & E \ar[r] & 0}
and a canonical isomorphism
\[\kn_{|C}/(E\ot L) \ \simeq \ \kn^{(2)} \ . \]
\end{subsub}

\begin{proof} The exact sequence follows from lemma \ref{lem1c}. Consider the 
exact sequence
\xmat{0\ar[r] & \kn\ar[r]^-i & \E\ar[r]^-p & E\ar[r] & 0 \ .}
We have a canonical surjective morphism \ 
\m{\phi:\kn_{|C}=\kn/\kn_1\to\kn^{(2)}=\kn/\kn^{(1)}}. Let \m{x\in C} and 
\m{z\in\ko_{Y,x}} and equation of $C$. Let \m{u\in\kn_{|C,x}}, and 
\m{\ov{u}\in\kn_x} over $u$. Then we have
\begin{eqnarray*}
i_{|C,x}(u)=0 & \Longleftrightarrow & i_x(\ov{u}) \ \text{is a multiple of} \ 
z\\
 & \Longleftrightarrow & z.i_x(\ov{u})=i_x(z\ov{u})=0, \ \ \text{because} \ \E 
 \ \text{is locally free}\\
 & \Longleftrightarrow & z\ov{u}=0\\
 & \Longleftrightarrow & \phi_x(u)=0.
\end{eqnarray*}
This proves the isomorphism of lemma \ref{lem1b}. \end{proof}

\sepprop

\begin{subsub}\label{pu_th}Punctual thickening of vector bundles on $C$. \rm
Let \ \m{D=\sigg_im_iP_i} \ be a divisor on $C$, where the \m{P_i} are distinct 
points and the \m{m_i} positive integers. Let \m{x_i\in\ko_{C_2,P_i}} be over a 
generator of the maximal ideal of \m{\ko_{C,P_i}}, and \m{z_i\in\ko_{C_2,P_i}} 
a local equation of $C$. Let \m{\ki[D]} be the ideal sheaf on \m{C_2} defined 
by: \m{\ki[D]_P=L_P} \ for $P$ distinct of all the \m{P_i}, and 
\m{\ki[D]_{P_i}=(x_i^{m_i}z_i)=x_i^{m_i}.L_{P_i}} \ . It depends only on $D$.

Let $E$ be a vector bundle on $C$. According to \cite{dr2}, thm. 3.1.1, $E$ can 
be extended to a vector bundle $\E$ on \m{C_2}. Let
\[E[D] \ = \ \E/(\E\ot\ki[D]) \ . \]
We have a canonical quotient  \  \m{E[D]\to E} \ which is an isomorphism ouside 
the \m{P_i}.
It is easy to see that the sheaf \m{E[D]} depends only of $E$ and $D$, i.e. if 
\m{\E'} is another vector bundle on \m{C_2} extending $E$, then there exists an 
isomorphism (not unique)
\[\E/(\E\ot\ki[D]) \ \simeq \ \E'/(\E'\ot\ki[D])\]
inducing the identity on $E$. We have \ \m{T(E[D])=\sigg T_i}, where \m{T_i} 
is the skyscraper sheaf concentrated at \m{P_i} with fibre \ 
\m{E_{P_i}/x_i^{m_i}E_{P_i}} \ at \m{P_i} (where \m{E_{P_i}} the fibre of the 
sheaf $E$ at this point). We have \ \m{E[D]/T(E[D])=E}.
\end{subsub}

\end{sub}

\sepsec

\section{Coherent sheaves on primitive reducible curves}\label{prim_red}

A {\em primitive reducible curve} is an algebraic curve $X$ such that
\begin{enumerate}
\item[--] $X$ is connected.
\item[--] The irreducible components of $X$ are smooth projective curves.
\item[--] Any three irreducible components of $X$  have no common point, and 
any two components are transverse.
\end{enumerate}

Let $X$ be a primitive reducible curve, $\bf C$ the set of components of $X$, 
and $I$ the set of intersection points of the components of $X$.

\sepsub

\Ssect{reflexive sheaves}{reflex}

The following result is an easy consequence of \cite{ses}, VII, VIII:

\sepprop

\begin{subsub}\label{prop1}{\bf Proposition: } Let $\ke$ be a coherent sheaf on
$X$. Then the following conditions are equivalent:
\begin{enumerate}
\item[(i)] $\ke$ is pure of dimension 1.
\item[(ii)] $\ke$ is of depth 1.
\item[(iii)] $\ke$ is locally free at every point of $X$ belonging to only one 
component, and if a point $x$ belongs to two components $D_1$, $D_2$, then 
there exist integers $a$, $a_1$, $a_2\geq 0$ and an isomorphism
\[\ke_x \ \simeq \ a\ko_{X,x}\oplus a_1\ko_{D_1,x}\oplus a_2\ko_{D_2,x} \ . \]
\item[(iv)] $\ke$ is torsion free, i.e. for every $x\in X$, every element of 
$\ko_{X,x}$ which is not a zero divisor in $\ko_{X,x}$ is not a zero divisor in 
$\ke_x$.
\item[(v)] $\ke$ is reflexive.
\end{enumerate}
\end{subsub}

\sepprop

In particular, every torsion free coherent sheaf on $X$ is reflexive.

\sepprop

\begin{subsub}\label{def_tau}{\bf Definition: }
Let $\ke$ be a torsion free sheaf on $X$. For every \m{D\in\bf C}, let \ 
\m{\ke_D=\ke_{|D}/T} (where $T$ is the torsion subsheaf of \m{\ke_{|D}}). The 
sequence of pairs \ \m{\tau(\ke)=(\rk(\ke_D),\deg(\ke_D))_{D\in\bf C}} \ is 
called the {\em type} of $\ke$.
\end{subsub}

\end{sub}

\sepsub

\Ssect{Structure and properties of torsion free sheaves}{struc}

\begin{subsub}\label{lin} Local structure at intersection points -- \rm
Let \m{x\in X}, belonging to two irreducible components \m{D_1}, \m{D_2}. The 
ring \m{\ko_{X,x}} can be identified with the ring of pairs 
\m{(\phi_1,\phi_2)\in\ko_{D_1,x}\times\ko_{D_2,x}} \ such that \ 
\m{\phi_1(x)=\phi_2(x)}. For \m{i=1,2}, let \m{t_i\in\ko_{D_i,x}} be a 
generator of the maximal ideal. We will also denote by \m{t_1} (resp. \m{t_2}) 
the element \m{(t_1,0)} (resp. \m{(0,t_2)}) of \m{\ko_{X,x}}. Let $M$ be the
\m{\ko_{X,x}}-module
\[M \ = \ (V\ot\ko_{X,x})\oplus (W_1\ot\ko_{D_1,x})\oplus (W_2\ot\ko_{D_2,x})\]
(where $V$, \m{W_1}, \m{W_2} are finite dimensional $\C$-vector spaces). Then
\m{\End(M)} consists of matrices
\[\begin{pmatrix}
\alpha & \lambda_1t_1 & \lambda_2t_2 \\ \beta_1 & \mu_1 & 0 \\ \beta_2 & 0 &
\mu_2
\end{pmatrix}
\]
with \m{\alpha\in\End(V)}, \m{\lambda_i\in 
L(W_i,V)\ot\ko_{D_i,x}}, \m{\beta_i\in L(V,W_i)\ot\ko_{D_i,x}},
\m{\mu_i\in\End(W_i)\ot\ko_{D_i,x}} \ for \m{i=1,2}. Such an endomorphism is 
invertible if and only if $\alpha$, \m{\mu_1} and \m{\mu_2} are invertible, if 
and only if their images in \m{\End(V)}, \m{\End(W_1)}, \m{\End(W_2)} 
respectively are invertible. It is easy to see that the sub-module
\[N \ = \ (V\ot{\mathfrak m}_x)\oplus (W_1\ot\ko_{D_1,x})\oplus 
(W_2\ot\ko_{D_2,x})\]
is invariant by all the endomorphisms of $M$. The module $M$ can be 
identified with the \m{\ko_{X,x}}-module of
\[(u_1,u_2,v_1,v_2) \ \in \ (V\ot\ko_{D_1,x})\oplus (V\ot\ko_{D_2,x})\oplus
(W_1\ot\ko_{D_1,x})\oplus (W_2\ot\ko_{D_2,x})\]
such that the images of \m{u_1,u_2} in $V$ are the same. 

It follows easily that if $U$ is a neighbourhood of $x$ which contains only one 
multiple point and meets only the components through $x$, then the reflexive 
sheaves $\ke$ on $U$ are obtained in the following way: take two vector 
bundles, \m{E_1} on \m{D_1\cap U}, \m{E_2} on \m{D_2\cap U} respectively, a 
$\C$-vector space $W$ and surjective linear maps
\[f_i:E_{i,x}\lra W \ , \qquad i=1,2 \ . \]
Then \m{\ke=E_1} on \m{D_1\backslash\{x\}}, \m{\ke=E_2} on 
\m{D_2\backslash\{x\}}, and
\[\ke_x \ = \ \{(\phi_1,\phi_2)\in E_1(x)\times E_2(x) \ ; \ f_1(\phi_1(x))=
f_2(\phi_2(x))\} \ , \]
(where \m{E_i(x)} is the module of sections of \m{E_i} defined in a 
neighbourhood on $x$). We have
\[E_i \ = \ \ke_{|D_i\cap U}/T_i \ , \]
where \m{T_i} is the torsion subsheaf of \m{\ke_{|D_i\cap U}}. We have
\[\ke_x \ \simeq \ (W\ot\ko_{X,x})\oplus(\ker(f_1)\ot\ko_{D_1,x})\oplus
(\ker(f_2)\ot\ko_{D_2,x}) \ . \]
The sheaf $\ke$ is locally free if and only if \m{f_1} and \m{f_2} are 
isomorphisms.

We have an exact sequence
\xmat{0\ar[r] & \ke\ar[r] & E_1\oplus E_2\ar[r]^-p & \ov{W}\ar[r] & 0,}
where \m{\ov{W}} denotes the skyscraper sheaf concentrated on $x$ with fibre 
$W$ at $x$, and $p$ is given at $x$ by
\[(f_1,-f_2):E_{1,x}\oplus E_{2,x}\lra W \ . \]
\end{subsub}

\sepprop

\begin{subsub}\label{oth_cnst} Another way of constructing sheaves -- \rm
We keep the notations of \ref{lin}. Let \Nligne \m{H\subset E_{1,x}\times 
E_{2,x}} \ be a linear subspace. We define a reflexive sheaf $\ke$ on $U$ by:
\begin{enumerate}
\item[--] $\ke=E_1$ on $D_1\backslash\{x\}$, and $\ke=E_2$ on 
$D_2\backslash\{x\}$,
\item[--] $\ke_x\ = \ \{(\phi_1,\phi_2)\in E_1(x)\times E_2(x) \ ; \ 
(\phi_1(x),\phi_2(x))\in H\}$ \ .
\end{enumerate}
Let \ \m{H_1\subset E_{1,x}} and \m{H_2\subset E_{2,x}} be the minimal linear 
subspaces such that \ \m{H\subset H_1\times H_2}, and \ \m{W=(H_1\times H_2)/H}.
Let \ \m{\tau_1:H_1=H_1\times\nsp\to W}, \m{\tau_2:H_2=\nsp\times H_2\to W} \ 
be the restrictions of the quotient map. Let \ \m{E'_1=(E_1)_{H_1}}, 
\m{E'_2=(E_2)_{H_2}} (cf. \ref{el_mod}). We have then canonical maps \ 
\m{\lambda_1:E'_{1,x}\to H_1}, \m{\lambda_2:E'_{2,x}\to H_2}. Let
\[f_1:\tau_1\circ\lambda_1:E'_{1,x}\to W \ , \quad 
f_2:\tau_2\circ\lambda_2:E'_{2,x}\to W \ . \]
Then the two linear maps \m{f_1}, \m{f_2} are surjective, and $\ke$ is the 
reflexive sheaf constructed from \m{E'_1}, \m{E'_2}, \m{f_1}, \m{f_2} by the 
method of \ref{lin}.
\end{subsub}

\sepprop

\begin{subsub}\label{oth_ex} Example: sheaves of homomorphisms -- \rm Let 
\m{F_i,G_i}, \m{i=1,2}, be vector bundles on \m{D_i\cap U}. Let $V$, $W$ be 
vector spaces, and
\[f_i:F_{i,x}\lra V \ , \quad g_i:G_{i,x}\lra W\]
surjective maps. Let $\kf$ (resp. $\kg$) be the reflexive sheaf on $U$ defined
by \m{F_1}, \m{F_2}, \m{f_1}, \m{f_2} (resp. \m{G_1}, \m{G_2}, \m{g_1},
\m{g_2}). We will describe the sheaf \m{\HHom(\kf,\kg)}. Of course we have \ 
\m{\HHom(\kf,\kg)=\HHom(F_i,G_i)} \ on \m{D_i\backslash\{x\}}. From the 
definitions of $\kf$, $\kg$ it is easy to see that
\[\HHom(\kf,\kg)_x \ = \ \{(\mu_1,\mu_2)\in\HHom(F_1,G_1)(x)\times
\HHom(F_2,G_2)(x) \ ; \ (\mu_1(x),\mu_2(x))\in H\} \ , \]
where \ \m{H\subset L(F_{1,x},G_{1,x})\times L(F_{2,x},G_{2,x})} \ is the 
linear subspace consisting of \m{(\alpha_1,\alpha_2)} such that for every \ 
\m{(u_1,u_2)\in F_{1,x}\times F_{2,x}} \ such that \m{f_1(u_1)=f_2(u_2)}, we 
have also \ \m{g_1(\alpha_1(u_1))=g_2(\alpha_2(u_2))}. If follows that 
\m{\HHom(\kf,\kg)} is reflexive.

We can suppose that we have decompositions into direct sums \ \m{F_i=(V\oplus 
A_i)\ot\ko_{D_i\cap U}}, \m{G_i=(W\oplus B_i)\ot\ko_{D_i\cap U}}, in such a way
that \m{f_i} (resp. \m{g_i}) is the projection to $V$  (resp. $W$). We can then
write linear maps \m{\alpha_i:F_{i,x}\to G_{i,x}} as matrices \m{\begin{pmatrix}
\phi_i & M_i\\ N_i & P_i\end{pmatrix}}. It is then easy to see that $H$ consists
of pairs \m{(\alpha_1,\alpha_2)} such that \m{M_1=0}, \m{M_2=0},
\m{\phi_1=\phi_2}. The minimal subspaces \m{H_i\subset L(F_{i,x},G_{i,x})} such
that \m{H\subset 
H_1\times H_2} are then clearly
\[H_i \ = \ \{\alpha_i\in L(F_{i,x},G_{i,x}) \ ; \ 
\alpha_i(\ker(f_i))\subset\ker(g_i)\} \ , \]
and \m{(H_1\times H_2)/H} is canonically isomorphic to \m{L(V,W)}. The linear 
maps \Nligne \m{\gamma_i:\HHom(F_i,G_i)_{H_i,x}\lra L(V,W)} \ are as follows:
\xmat{\gamma_i:\HHom(F_i,G_i)_{H_i,x}\flon[r] & H_i\ar[r]^-{\rho_i} & L(V,W)}
where \m{\rho_i} associates to \m{\alpha:F_{i,x}\to G_{i,x}} the induced map \
\m{V=F_{i,x}/\ker(f_i)\to W=G_{i,x}/\ker(g_i)}. The reflexive sheaf 
\m{\HHom(\kf,\kg)} is then defined (using the method of \ref{lin}) by the 
vector bundles \m{\HHom(F_i,G_i)_{H_i}} on \m{D_i}, and the linear maps 
\m{\gamma_i}.
\end{subsub}

\sepprop

\begin{subsub}\label{lem1}{\bf Lemma:} We have, for \m{i=1,2} and integers 
\m{j\geq 0}, \m{k\geq 1}
\begin{enumerate}
\item[--] $\Ext^{2j+1}_{\ko_{X,x}}(\ko_{D_1,x},\ko_{D_2,x})\simeq\C$, 
$\Ext^{2k}_{\ko_{X,x}}(\ko_{D_1,x},\ko_{D_2,x})=0$ .
\item[--] $\Ext^{2j+1}_{\ko_{X,x}}(\ko_{D_i,x},\ko_{D_i,x})=0$, 
$\Ext^{2k}_{\ko_{X,x}}(\ko_{D_i,x},\ko_{D_i,x})\simeq\C$ .
\item[--] $\Tor^{2j+1}_{\ko_{X,x}}(\ko_{D_i,x},\ko_{D_i,x})\simeq\C$, 
$\Tor^{2k}_{\ko_{X,x}}(\ko_{D_i,x},\ko_{D_i,x})=0$ .
\item[--] $\Tor^{2j+1}_{\ko_{X,x}}(\ko_{D_1,x},\ko_{D_2,x})=0$, 
$\Tor^{2k}_{\ko_{X,x}}(\ko_{D_1,x},\ko_{D_2,x})\simeq\C$ .
\end{enumerate}
\end{subsub}

\begin{proof}
We use a locally free resolution of \m{\ko_{D_1,x}}:
\xmat{\cdots\ko_{X,x}\ar[r]^-{t_1} & \ko_{X,x}\ar[r]^-{t_2} & \ko_{X,x}\ar[r] &
\ko_{D_1,x}\ar[r] & 0 \ .}
Hence \m{\Ext^1_{\ko_{X,x}}(\ko_{D_1,x},\ko_{D_2,x})} is isomorphic to the 
middle cohomology of the complex
\xmat{\ko_{D_2,x}\ar[r]^-{t_2} & \ko_{D_2,x}\ar[rr]^-{t_1=0} & & \ko_{D_2,x} \ 
.}
The result follows immediately. The other equalities are proved in the same way.
\end{proof}

\sepprop

\begin{subsub}\label{rel} Relative version of lemma \ref{lem1} -- \rm
A similar proof shows that if $S$ is an algebraic variety and \m{s\in S}, then 
we have, for \m{i=1,2} and integers \m{j\geq 0}, \m{k\geq 1}
\[\Ext^{2j+1}_{\ko_{X\times 
S,(x,s)}}(\ko_{D_1,x}\ot\ko_{S,s},\ko_{D_2,x}\ot\ko_{S,s}) \ \simeq \ \ko_{S,s} 
\ , \quad \Ext^{2k}_{\ko_{X\times 
S,(x,s)}}(\ko_{D_1,x}\ot\ko_{S,s},\ko_{D_2,x}\ot\ko_{S,s}) \ = \ 0 \ , \]
\[\Ext^{2j+1}_{\ko_{X\times 
S,(x,s)}}(\ko_{D_i,x}\ot\ko_{S,s},\ko_{D_i,x}\ot\ko_{S,s}) \ = \ 0 \ , \quad
\Ext^{2k}_{\ko_{X\times 
S,(x,s)}}(\ko_{D_i,x}\ot\ko_{S,s},\ko_{D_i,x}\ot\ko_{S,s}) \ \simeq \ \ko_{S,s} 
\ . \]
and the same as in lemma \ref{lem1} for the Tor.
\end{subsub}

\sepprop

\begin{subsub}\label{linked_s} Linked sheaves -- \rm We keep the notations of 
\ref{lin}. Let \ \m{r=\dim(W)}, \m{r_1=\rk(E_1)}, \m{r_2=\rk(E_2)}. Then we 
have \ \m{r\leq\inf(r_1,r_2)}, \m{\rk(\ke_x)=r_1+r_2-r}. The sheaf $\ke$ is 
called {\em linked at $x$} if \ \m{r=\inf(r_1,r_2)}, the maximal possible value 
(so one could say that the most unlinked sheaf is \m{E_1\oplus E_2}). Suppose 
that $\ke$ is linked at $x$. For \m{i=1,2}, we have \m{E_i=\ke_{|D_i}/T_i} 
(where \m{T_i} is the torsion subsheaf of \m{\ke_{|D_i}}). Suppose that 
\m{r_2\geq r_1}. Then \m{T_2=0} and \m{T_1\simeq\C_x\ot\C^{r_1-r_2}} (where 
\m{\C_x} is the torsion sheaf concentrated on $x$ with fibre $\C$ at that 
point).

A torsion free sheaf $\kf$ on $X$ is called {\em linked} if is linked at every 
intersection point of any two irreducible components of $X$. 
\end{subsub}

\sepprop

\begin{subsub}\label{ext_sh} Extensions of vector bundles on two components --
\rm
Let \m{V_1}, \m{V_2} be finite dimensional vector spaces. From lemma 
\ref{lem1}, we have a canonical isomorphism
\begin{equation}\label{equ12}\Ext^1_{\ko_{X,x}}(\ko_{D_1,x}\ot 
V_1,\ko_{D_2,x}\ot V_2) \ \simeq \ 
L(V_1,V_2) \ . \end{equation}
Let \ \m{\phi\in L(V_1,V_2)}. We now give an explicit description of the 
extension
\[0\lra\ko_{D_2,x}\ot V_2\lra M\lra\ko_{D_1,x}\ot V_1\lra 0\]
corresponding to $\phi$. We have an exact sequence
\xmat{0\ar[r] & \ko_{D_2,x}\ot V_1\ar[rr]^-{(t_2,\phi)} & & (\ko_{X,x}\ot V_1)
\oplus(\ko_{D_2,x}\ot V_2)\ar[r] & M\ar[r] & 0}
(cf. \cite{dr1b}, 4--), and
\[M \ = \ (\imm(\phi)\ot\ko_{X,x})\oplus (\ker(\phi)\ot\ko_{D_1,x})\oplus 
(\coker(\phi)\ot\ko_{D_2,x}) \ 
. \]
The module $M$ can be the fibre of some linked torsion free sheaf if and only 
$\phi$ has maximal rank, i.e. is injective or surjective.

For vector bundles on \m{D_1\cup D_2\subset X} we need more canonical 
isomorphisms. Let \ \m{D=D_1\cap D_2}. We use the obvious locally free 
resolution of \m{\ko_{D_1}}
\xmat{\cdots\kl_3\ar[r] & \kl_2\ar[r] & \ko_X\ar[r] &
\ko_{D_1}\ar[r] & 0 \ ,}
where \m{\kl_2} is a line bundle on \m{D_1\cup D_2} with 
\m{\kl_{2|D_2}=\ko_{D_2}(-D)}, \m{\kl_{2|D_1}=\ko_{D_1}}, and \m{\kl_3} is 
a line bundle on \m{D_1\cup D_2} with \m{\kl_{3|D_1}=\ko_{D_1}(-D)}, 
\m{\kl_{3|D_2}=\ko_{D_2}}. In this way we obtain, if \m{E_1}, \m{E_2} are 
vector bundles on \m{D_1}, \m{D_2} respectively, a canonical isomorphism
\begin{equation}\label{equ10}\Ext^1_{\ko_{D_1\cup D_2}}(E_1,E_2) \ = \ 
\bigoplus_{x\in 
D}L(E_1(-x)_x,E_{2,x}) \ . \end{equation}

Now we make the link between this description of extensions and the 
construction of \ref{lin}. Let \m{E_1}, \m{E_2} be vector bundles on \m{D_1}, 
\m{D_2} respectively. We want to describe torsion free sheaves $\ke$ on 
\m{D_1\cup D_2} such that \ \m{\ke_{|D_i}/T_i=E_i} \ for \m{i=1,2} (where 
\m{T_i} denotes the torsion subsheaf). As in in \ref{lin}, $\ke$ is defined by 
vector spaces \m{W^x}, \m{x\in D}, and surjective maps
\begin{equation}\label{equ14}f_1^x:E_{1,x}\lra W^x \ , \quad f_2^x:E_{2,x}\lra 
W^x \ . \end{equation}
Let \ \m{N_1^x=\ker(f_1^x)}, \m{N_2^x=\ker(f_2^x)} \ for \m{x\in D}, 
\m{i=1,2}. We have an exact sequence
\begin{equation}\label{equ11}0\lra\big(E_2\big)_{(N_2^x)_{x\in D_1\cap 
D_2}}\lra\ke\lra E_1\lra 0\end{equation}
(cf. \ref{el_mod} for the definition of \m{\big(E_2\big)_{(N_2^x)_{x\in D_1\cap 
D_2}}}). The corresponding linear map \m{(\ref{equ12})} for every 
\m{x\in D} is as follows:
\xmat{(E_1(-x))_x\ar[rr]^-{f_1^x} & &
W^x\ot\ko_X(-x)_x=E_2(-x)_x/N_2^x\ot\ko_X(-x)_x\flinc[r] & (E_2)_{(N_2^x),x}
}
(the last inclusion coming also from \ref{el_mod}).

We have also a canonical exact sequence
\xmat{0\ar[r] & \ke\ar[r] & E_1\oplus E_2\ar[r]^-p & \bigoplus_{x\in D
}\ov{W^x}\ar[r] & 0}
where \m{\ov{W_x}} denotes the skyscraper sheaf concentrated on $x$ with fibre 
\m{W_x} at $x$, and $p$ is given at $x$ by
\[(f_1^x,-f_2^x):E_{1,x}\oplus E_{2,x}\lra W^x \ . \]
This can of course be generalised to torsion free sheaves on the whole of $X$:
recall that $\bf C$ denotes the set of components of $X$ and $I$ the set of 
intersection points. Let $\kf$ be a torsion free sheaf on $X$, and for every 
\m{D\in{\bf C}}, \m{\kf_D=\kf_{|D}/T_D} (where \m{T_D} is the torsion 
subsheaf). Then we have an exact sequence 
\begin{equation}\label{equ17}
0\lra\kf\lra\bigoplus_{D\in\bf C}\kf_D\lra\bigoplus_{P\in I}(\C_P\ot W_P)\lra 0 
\ ,
\end{equation}
where for each \m{P\in I}, \m{W_P} is a vector space. If $\kf$ is linked, such 
a $P$ belongs to two components \m{X_1}, \m{X_2}, and \m{W_P=\kf_{X_i,x}}, 
where $i$ is such that \ \m{\rk(\kf_i)=\inf(\rk(\kf_1),\rk(\kf_2))}.
\end{subsub}

\sepprop

\begin{subsub}\label{dual_ext} Duality -- \rm To simplify the notations, we 
suppose in the rest of \ref{struc} that $X$ has two components \m{D_1}, \m{D_2} 
(it is easy to generalise the following results to the case where $X$ has more 
than two components). Recall that \m{D=D_1\cap D_2}. If $\ke$ is a coherent
sheaf on $X$, let \m{\ke^\vee} denote its dual. If $E_1$ is a vector bundle on
\m{D_1} we have
\[E_1^\vee \ \simeq \ E_1^*(-D) \ . \]

We keep the notations of \ref{ext_sh}. For \m{i=1,2}, let 
\m{Q_i^x=E_{i,x}/N^x_i}. From equation \m{(\ref{equ11})} we deduce the exact 
sequence
\[0\lra E_1^*(-D)\lra\ke^\vee\lra\left[\big(E_2\big)_{(N_2^x)_{x\in 
D}}\right]^*(-D)\lra 0 \]
(cf. \ref{el_mod}).
According to \ref{GEM_dual}, we have
\[\left[\big(E_2\big)_{(N_2^x)_{x\in D}}\right]^*(-D) \ = \
(E_2^*)_{({Q_2^x}^*)_{x\in D}} \ . \]
It follows that if \ \m{G_i=\ke^\vee_{|D_i}/T_i}, \m{i=1,2}, then we have
\[G_i \ = \ (E_i^*)_{({Q_i^x}^*)_{x\in D}} \ . \]
The corresponding maps \m{(\ref{equ14})} are of course
\[f_i^{x*}:G_{i,x}=\left[(E_i^*)_{({Q_i^x}^*)}\right]\lra 
(E_{i,x}/N_i^x)^*=W_x^* \ 
. \]
\end{subsub}

\sepprop

\begin{subsub}\label{tens} Tensor products -- \rm Let $\ke$, \m{\ke'} be 
torsion free sheaves on \m{X=D_1\cup D_2}. Let \m{E_i=\ke_{|D_i}/T_i}, 
\m{E'_i=\ke'_{|D_i}/T'_i} for \m{i=1,2}, where \m{T_i}, \m{T'_i} denote the 
torsion subsheaves, and \ \m{r_i=\rk(E_i)}, \m{r'_i=\rk(E'_i)}. For every 
\m{x\in D}, let
\[f_1^x:E_{1,x}\lra W^x \ , \quad f_2^x:E_{2,x}\lra W^x \ , \quad 
{f'_1}^x:E'_{1,x}\lra W'^x \ , \quad {f'_2}^x:E'_{2,x}\lra W'^x \]
be surjective linear maps that define $\ke$ and \m{\ke'} at $x$ (cf. 
\ref{ext_sh}).
\end{subsub}

\sepprop

\begin{subsub}\label{prop10}{\bf Proposition: } Suppose that $\ke$ and \m{\ke'} 
are linked and that \ \m{(r_1-r_2)(r'_1-r'_2)\geq 0}. Then \m{\ke\ot\ke'} is 
torsion free and linked, and we have \
\m{(\ke\ot\ke')_{|D_i}/\T_i\simeq E_i\ot E'_i} \ for \m{i=1,2} (where \m{\T_i} 
is the torsion subsheaf), and for every 
\m{x\in D}, \m{\ke\ot\ke'} is defined at $x$ by the linear maps
\[f_1^x\ot{f'_1}^x:E_{1,x}\ot E'_{1,x}\lra W_x\ot W'_x \ , \quad
f_2^x\ot{f'_2}^x:E_{2,x}\ot E'_{2,x}\lra W_x\ot W'_x \ . \]
\end{subsub}

\begin{proof} The condition \m{(r_1-r_2)(r'_1-r'_2)\geq 0} implies that for 
every \m{x\in D}
\begin{enumerate}
\item[--] if $r_1\geq r_2$, then there exist integers $a$, $b$, $a'$, $b'$ such 
that
$$\ke_x\simeq a\ko_{D_1,x}\oplus b\ko_{X,x} \ , \quad \ke'_x\simeq 
a'\ko_{D_1,x}\oplus b'\ko_{X,x} \ , $$
\item[--] if $r_1<r_2$, then there exist integers $a$, $b$, $a'$, $b'$ such 
that 
$$\ke_x\simeq a\ko_{D_2,x}\oplus b\ko_{X,x} \ , \quad \ke'_x\simeq 
a'\ko_{D_2,x}\oplus b'\ko_{X,x} \ . $$
\end{enumerate}
This condition is necessary, otherwise \m{\ke\ot\ke'} would not be torsion 
free, because \ \m{\ko_{D_1,x}\ot\ko_{D_2,x}} \ is the torsion sheaf 
concentrated on $x$ with fibre $\C$.

Let $\kf$ be the linked sheaf defined by \m{E_1\ot E'_1}, \m{E_2\ot E'_2} and 
the preceding maps. Then there is an obvious morphism \ \m{\ke\ot\ke'\to\kf}. 
Using local decompositions of $\ke$, \m{\ke'} in direct sums (as indicated in 
the beginning of the proof) it is easy to see that it is an isomorphism.
\end{proof}

\sepprop

\begin{subsub}\label{rem1} Computation of \ \ \m{\Ext^1_{\ko_X}(\ke,\ke)} -- 
\rm 
Let $\ke$ be a linked torsion free sheaf on \m{X=D_1\cup D_2}. \Nligne
According to \ref{dual_ext} and \ref{tens}, the 
sheaf \ \m{\ke\ot\ke^\vee} \ is defined by the vector bundles \Nligne
\m{E_i\ot(E_i^*)_{((E_{i,x}/N_i^x)^*)}}, and the maps \m{f_i^x\ot f_i^{x*}} 
(recall that \ \m{N_i^x=\ker(f_i^x)}). \Nligne
According to \ref{oth_ex}, the sheaf \m{\EEnd(\ke)} is defined by the vector 
bundles \m{\EEnd(E_i)_{(H_i^x)}} and the maps \ 
\m{\gamma_i^x:\EEnd(E_i)_{(H_i^x),x}\to\End(W^x)}, where for every 
\m{x\in D=D_1\cap D_2}, \Nligne
\m{H_i^x=\{\alpha\in\End(E_{i,x});\alpha(\ker(f_i^x))\subset
\ker(f_i^x)\}}.

Since $\ke$ is linked, we can assume that \m{r_1\leq r_2}. In this case we have
\m{N_1^x=\nsp} for every \m{x\in D}, and
\[E_1\ot(E_1^*)_{({Q_1^x}^*)} \ = \ \EEnd(E_1)_{(H_1^x)} \ = \
\EEnd(E_1) \ , \]
and only the vector bundles on \m{D_2} can differ. We have two canonical 
commutative diagrams with exact rows and columns
\xmat{& & & 0\ar[d]\\
& 0\ar[d] & & \dsp\bigoplus_{x\in D}\ov{\End(N_2^x)}\ar[d]\\
0\ar[r] & E_2\ot(E_2^*)_{({Q_2^x}^*)}\ar[r]\ar[d] & 
\EEnd(E_2)\fleq[d]\ar[r] & \dsp\bigoplus_{x\in D}
\ov{L(N_{2x},E_{2,x})}\ar[d]\ar[r] & 0\\
0\ar[r] & \EEnd(E_2)_{(H_2^x)}\ar[r]\ar[d] &  \EEnd(E_2)\ar[r] & 
\dsp\bigoplus_{x\in D}\ov{L(N_2^x,Q_2^x)}\ar[r]\ar[d] & 0\\
& \dsp\bigoplus_{x\in D}\ov{\End(N_2^x)}\ar[d] & & 0\\
& 0}
\xmat{& 0\ar[d] & 0\ar[d]\\
0\ar[r] & \ke^\vee\ot\ke\ar[r]\ar[d] & \EEnd(E_1)\oplus
(E_2\ot(E_2^*)_{({Q_2^x}^*)})\ar[r]\ar[d] &  \dsp\bigoplus_{x\in 
D}\ov{\End(E_{1,x})}\ar[r]\fleq[d] & 0\\
0\ar[r] & \EEnd(\ke)\ar[r]\ar[d] & \EEnd(E_1)\oplus
\EEnd(E_2)_{(H_2^x)}\ar[r]\ar[d] & \dsp\bigoplus_{x\in 
D}\ov{\End(E_{1,x})}\ar[r] & 0\\
& \dsp\bigoplus_{x\in D}\ov{\End(N_2^x)}\fleq[r]\ar[d] &
\dsp\bigoplus_{x\in D}\ov{\End(N_2^x)}\ar[d]\\
& 0 & 0}
\end{subsub}

\sepprop

It follows that

\sepprop

\begin{subsub}\label{prop11}{\bf Proposition: } Suppose that \m{E_1} and 
\m{E_2} are simple. Then there is a vector space \m{H_\ke} and canonical 
exact sequences
\[0\lra\dsp\big(\bigoplus_{x\in D}\End(E_{1,x})\big)/\C\lra 
H_\ke\lra\dsp\bigoplus_{x\in D}L(N_2^x,Q_2^x)\lra 0 \ , \]
\[0\lra H_\ke\lra\Ext^1_{\ko_X}(\ke,\ke)\lra\Ext^1_{\ko_{D_1}}(E_1,E_1)\oplus
\Ext^1_{\ko_{D_2}}(E_2,E_2)\lra 0 \ . \]
In particular \m{\dim(\Ext^1_{\ko_X}(\ke,\ke))} depends only on \m{r_1}, 
\m{r_2}, \m{\dim(\Ext^1_{\ko_{D_1}}(E_1,E_1))} and \Nligne
\m{\dim(\Ext^1_{\ko_{D_2}}(E_2,E_2))}.
\end{subsub}

\begin{proof} From lemma \ref{lem1} we have \ \m{\EExt^1_{\ko_X}(\ke,\ke)=0}, 
hence \ \m{\Ext^1_{\ko_X}(\ke,\ke)\simeq H^1(\EEnd(\ke))}.

The identity morphism \m{E_2\to E_2} is also a section of 
\m{\EEnd(E_2)_{(H_2^x)}}, hence it follows from the preceding diagram that we 
have an exact sequence
\[0\lra\big(\bigoplus_{x\in D}\End(E_{1,x})\big)/\C\lra H^1(\EEnd(\ke))\lra
\Ext^1_{\ko_{D_1}}(E_1,E_1)\oplus H^1(\EEnd(E_2)_{(H_2^x)})\lra 0 \ . \]
From the diagram before we deduce the exact sequence
\[0\lra\bigoplus_{x\in D}L(N_2^x,Q_2^x)\lra H^1(\EEnd(E_2)_{(H_2^x)})\lra
\Ext^1_{\ko_{D_2}}(E_2,E_2)\lra 0 \ . \]
The result follows immediately.
\end{proof}

\end{sub}

\sepsub

\Ssect{Deformations of linked sheaves}{def_ref}

\begin{subsub}\label{struc_fami} Local structure of families of linked sheaves 
-- \rm We keep the notations of \ref{struc}. Let $Z$ be an algebraic variety 
and \m{z_0} a closed point of $Z$. Let $\F$ be a coherent sheaf on \m{U\times 
Z}, flat on $Z$, and such that \m{\kf=\F_{z_0}} is a reflexive linked sheaf on 
$U$.
\end{subsub}

\sepprop

\begin{subsub}\label{lem1l}{\bf Lemma:} There exists a neighbourhood \m{Z_0} of 
\m{z_0} in $Z$ such that for every closed point \m{z\in Z_0}, \m{\F_z} is a 
linked torsion free sheaf on $X$.
\end{subsub}

\begin{proof} According to \cite{ma1}, prop. 2.1, the deformations of torsion 
free sheaves are torsion free. The fact that \m{\F_z} is linked for $z$ in a 
neighbourhood \m{Z_0} of \m{z_0} follows then easily from the semicontinuity of 
the rank of the fibres of coherent sheaves, since the rank of torsion free 
sheaves on $X$ at the intersection points is minimal precisely when the sheaves 
are linked.
\end{proof}

\sepprop

From now on, we assume that \m{Z_0=Z}. Let \m{r_i=\rk(\kf_{|D_i})} and suppose 
that \m{r_2\geq r_1}. Then \ \m{\F_2=\F_{|Z\times D_2}} \ is a vector bundle on 
\m{Z\times D_2\cap U} (because it is locally free of rank \m{r_2} on 
\m{D_2\cap(U\backslash\{x\})\times Z} and also on \m{\{x\}\times Z}). Similarly 
the kernel \m{\F'_1} of the restriction morphism \ \m{\F\to\F_{D_2\times Z}} \ 
is a vector bundle on \m{Z\times D_1}. According to lemma \ref{lem1} and 
\ref{rel}, for every \m{z\in Z} we have
\[\F_{z,x} \ \simeq \ r_1\ko_{X\times Z,(x,z)}\oplus(r_2-r_1)\ko_{D_2\times 
Z,(x,z)} \ . \]
Hence if \m{\T_1} is the torsion subsheaf of \m{\F_{|Z\times D_1}}, then \ 
\m{\F_1=\F_{|Z\times D_1}/\T_1} \ is a vector bundle on \m{Z\times D_1}, and we 
have a canonical isomorphism \ \m{\F'_1=\F_1\ot p_X^*(\ko_X(-x))} (where 
\m{p_X} is the projection \ \m{Z\times X\to X}).

\sepprop

\begin{subsub}\label{lem1m}{\bf Theorem: } {\bf 1 -- } Suppose that $Z$ is 
connected. Let $\U$ be a coherent sheaf on \m{X\times Z}, flat on $Z$, and such 
that for every closed point \m{z\in Z}, \m{\U_z} is a reflexive linked sheaf on 
$Z$. Let \m{D\in\bf C}. Let \m{\T_D} be the torsion subsheaf of \m{\U_{|Z\times 
D}} and \m{\U_D=\U_{|Z\times D}/\T_D}. Then \m{\U_D} is a vector bundle on 
\m{Z\times D}. If \m{z\in Z} is a closed point, then \ \m{\U_{D,z}=\U_{z|D}/T} 
(where $T$ is the torsion subsheaf of \m{\U_{z|D}}). In particular 
\m{\deg(\U_{D,z})} is independent of $z$. So the type of \m{\U_z} in 
independent of $z$.

{\bf 2 -- } Conversely, let $\ke$, $\kf$ be linked torsion free sheaves on $X$, 
such that \ \m{\tau(\ke)=\tau(\kf)} (cf. \ref{def_tau}).
Then there exist an integral variety $Y$, a flat family $\V$ of linked torsion 
free sheaves on $X$ parametrised by $Y$, and two closed points \m{s,t\in Y} 
such that \ \m{\V_s\simeq\ke} and \m{\V_t\simeq\kf}.
\end{subsub}

\begin{proof}
{\bf 1 } follows immediately from the discussion after lemma \ref{lem1l}. To 
prove {\bf 2} we remark first that we have exact sequences
\[0\lra\ke\lra\bigoplus_{D\in\bf C}\ke_D\lra\bigoplus_{P\in I}(\C_P\ot W_P)\lra 
0 \ ,\]
\[0\lra\kf\lra\bigoplus_{D\in\bf C}\kf_D\lra\bigoplus_{P\in I}(\C_P\ot W_P)\lra 
0 \ ,\]
(cf. $(\ref{equ17})$) with the same vector spaces \m{W_P}, where for every 
\m{D\in\bf C}, \m{\ke_D=\ke_{|D}/T_{\ke,D}} and \m{\kf_D=\kf_{|D}/T_{\kf,D}} 
(\m{T_{\ke,D}}, \m{T_{\kf,D}} are the torsion subsheaves). For every \m{D\in\bf 
C}, since \m{\ke_D} and \m{\kf_D} have the same rank and degree, there exist an 
integral variety \m{Y_D}, a vector bundle \m{\V_D} on \m{D\times Y_D} and two 
closed points \m{s_S,t_D\in Y_D} such that \m{\V_{D,s_D}\simeq\ke_D} and 
\m{\V_{D,t_D}\simeq\kf_D}. Let \ \m{{\bf Y}=\prod_{D\in\bf C}\Y_D}, and
\[\B \ = \ p_{X*}(\HHom(\bigoplus_{D\in\bf C}p_D^\#(\V_D),\bigoplus_{P\in 
I}p_P^*(\C_P\ot W_P)) \ , \]
(where \m{p_X}, \m{p_D}, \m{p_P} are the appropriate projections), which is a 
vector bundle on $\bf Y$. For every \m{y=(y_D)_{D\in\bf C}\in\bf Y}, 
\m{\dsp\B_y=\Hom(\bigoplus_{D\in\bf C}\V_{D,y_D},\bigoplus_{P\in I}\C_P\ot 
W_P)}. Now it suffices to take for $Y$ the open subset \m{\B_0} of $\B$ of 
surjective morphisms and for $\V$ the kernel of the obvious universal 
surjective morphism on \m{X\times\B_0}
\[\pi^\#(\bigoplus_{D\in\bf C}p_D^\#(\V_D))\lra\pi^\#(\bigoplus_{P\in 
I}p_P^*(\C_P\ot W_P))\]
(where $\pi$ is the projection \m{\B_0\to\bf Y}).
\end{proof}

\sepprop

\begin{subsub}\label{struc_fami2}\rm The description of the local structure of 
$\U$ implies also that the exact sequence $(\ref{equ17})$ can be globalised:
\[0\lra\U\lra\bigoplus_{D\in\bf C}\U_D\lra\bigoplus_{x\in I}p_{x*}(W_x)\lra 0 \ 
, \]
where $\bf C$ is the set of components of $X$, $I$ is the set of intersection
points, and for each \m{x\in I}, \m{W_x} is a vector bundle on $Z$ and \m{p_x} 
is the inclusion \ \m{\{x\}\times Z\to X\times Z}. More precisely, $x$ belongs 
to two components \m{D_1}, \m{D_2}, and \m{W_x=\U_{D_i|\{x\}\times Z}}, where 
$i$ is such that \ \m{\rk(\U_i)=\inf(\rk(\U_1),\rk(\U_2))}.
\end{subsub}

\sepprop

\begin{subsub}\label{loc_un} Construction of complete families -- \rm Let $\ke$ 
be a linked torsion free sheaf on $X$,\break \m{E_D=\ke_{|D}/T_D} for 
\m{D\in\bf C} (where \m{T_D} is the torsion subsheaf), \m{r_D=\rk(E_D)}.  
Then there exists a {\em complete deformation} \m{\E_D} of \m{E_D} 
parametrised by a smooth irreducible variety \m{Z_D}: it is defined by a vector 
bundle of \m{D\times Z_D} and a closed point \m{z_D\in Z_D} such that 
\m{E_D\simeq\E_{D,z_D}}, with the following local universal property: for every 
flat family $\kf$ of sheaves on $D$ parametrised by an algebraic variety $Y$, 
and \m{y\in Y} such that \m{\kf_y\simeq E_D}, there exists a neighbourhood 
\m{Y_0} of $y$ and a morphism \ \m{f:Y_0\to Z_D} \ such that \m{f(y)=z_D} and 
that there is an isomorphism \ \m{f^\#(\E_D)\simeq\kf_{|D\times Y_0}} 
compatible with the isomorphisms \m{\kf_y\simeq E_D}, \m{E_D\simeq\E_{D,z_D}}. 
This deformation can be constructed for example using Quot schemes. Now let 
\m{\boldsymbol{Z}=\prod_{D\in\C}Z_D} and \m{p_D} the projections 
\m{\boldsymbol{Z}\to Z_D}. For every intersection point $x$, let 
\m{D_1^x,D_2^x} be the components through $x$, such that \ 
\m{\rk(E_{D_1^x})\leq\rk(E_{D_2^x})}. Let $\F_x$ be the vector bundle on 
$\boldsymbol{Z}$ with fibre at \m{(z_D)} equal to 
\[\F_x \ = \ L(\E_{D_2^x,(x,z_{D_2^x})},\E_{D_1^x,(x,z_{D_1^x})}) \ , \]
and \ \m{\F=\bigoplus_{x\in I}\F_x}. Let \m{\F_0\subset\F} be the open subset 
corresponding to surjective morphisms, and \ \m{q:\F_0\to\boldsymbol{Z}} the 
projection. Let $\G$ be the torsion sheaf on \m{X\times \boldsymbol{Z}}
\[\G \ = \ \bigoplus_{x\in I}p_{D_1^x}^\#(\E_{D_1^x|\{x\}\times Z_{D_1^x}}) \ ,
\]
i.e. for every \m{(z_D)\in \boldsymbol{Z}}, \m{\G_{(z_D)}} is the torsion sheaf 
on $X$
\[\G_{(z_D)} \ = \ \bigoplus_{x\in I}\E_{D_1^x,(x,z_D)}\ot\C_x \ . \]
Then we have a canonical surjective morphism of sheaves on \m{X\times\F_0}
\[\Phi:\bigoplus_{D\in\bf C}q^\#(p_D^\#(\E_D))\lra q^\#(\G) \ . \]
Let \m{z=(z_D)_{D\in\bf C}\in\boldsymbol{Z}}, \m{\phi=(\phi_x)_{x\in I}\in\F_0} 
over $z$, 
so \m{\phi_x} is a surjective map \Nligne 
\m{\E_{D_2^x,(x,z_{D_2^x})}\to\E_{D_1^x,(x,z_{D_1^x})}}. For \m{x\in I}, 
\m{\Phi} at \m{(x,\phi)} is the linear map which is
\[(-I,\phi_x):\E_{D_1^x,(x,z_{D_1^x})}\oplus\E_{D_2^x,(x,z_{D_2^x})}\lra
\E_{D_1^x,(x,z_{D_1^x})}\]
on the two summands \m{\E_{D_1^x,(x,z_{D_1^x})}}, \m{\E_{D_2^x,(x,z_{D_2^x})}}, 
and zero on the others. Let \ \m{\boldsymbol{\ke}=\ker(\Phi)}. 

It follows then easily from \ref{struc_fami2} that $\boldsymbol{\ke}$ has the 
following {\em local universal property}: under the hypotheses of theorem 
\ref{lem1m}, if \m{z\in Z} is a closed point such that \m{\U_z\simeq\ke}, there 
exists a neighbourhood \m{Z_0} of $z$ and a morphism \m{f:Z_0\to\F_0} such that 
\ \m{\U_{X\times Z_0}\simeq f^\#(\boldsymbol{\ke})}.
\end{subsub}

\sepprop

\begin{subsub}\label{prop_st} Properties of stable linked sheaves -- \rm
Let $\ko_X(1)$ be an ample line bundle on $X$. For every polynomial $P$ in one 
variable and rational coefficients, let \m{\km_X(P)} denote the moduli space of 
sheaves on $X$, stable with respect to \m{\ko_X(1)} and with Hilbert 
polynomial $P$. 
\end{subsub}

\sepprop

\begin{subsub}\label{theo4}{\bf Theorem: } Let $\ke$ be a linked reflexive 
sheaf on $X$ and $P$ its Hilbert polynomial. Suppose that $\ke$ is stable with 
respect to \m{\ko_X(1)}, and that the restrictions of $\ke$ to the irreducible 
components of $X$ modulo torsion are simple vector bundles. Then \m{\km_X(P)} 
is smooth at the point corresponding to $\ke$.
\end{subsub}

\begin{proof}
We use the construction of \m{\km_X(P)} as a good quotient of an open subset 
\m{\bf R} of a \m{Quot} scheme (cd \cite{ma2}, \cite{si}) by a reductive group 
$G$. The quotient morphism \ \m{\pi:{\bf R}\to\km_X(P)} \ is associated to the 
universal sheaf $\ke$ on \m{X\times{\bf R}}. Let \m{{\bf R}_\ell} be the 
$G$-invariant open subset of points $y$ of $\bf R$ such that \m{\ke_y} is 
linked. According to \ref{loc_un}, $\pi$ can be locally factorised, at any 
point $y$ of \m{{\bf R}_\ell}, through a smooth variety: \m{\pi:{\bf 
R}\to\F_0\to\km_X(P)}. It follows that \m{\pi({\bf R}_\ell)} is contained in 
\m{\km_X(P)_\text{red}}, hence \m{\pi({\bf R}_\ell)} is integral.

At any point $z$ of \m{\km_X(P)} corresponding to a stable sheaf $\ke$, the 
tangent space of \m{\km_X(P)} at $z$ is canonically isomorphic to 
\m{\Ext^1_{\ko_X}(\ke,\ke)}. According to proposition \ref{prop11} and 
theorem \ref{lem1m}, if $\ke$ is linked, then the dimension of this space 
depends only on $P$, hence is constant. It follows that \m{\pi({\bf R}_\ell)} 
is smooth.
\end{proof}

\sepprop

\begin{subsub}\label{rem2}{\bf Remark: }\rm Let $\ke$ be a linked torsion free 
sheaf on $X$. Then we have \Nligne 
\m{H^0(\EExt^2_{\ko_X}(\ke,\ke))\subset\Ext^2_{\ko_X}(\ke,\ke)}. It follows 
from lemma \ref{lem1} that if $\ke$ is not locally free then \ 
\m{\Ext^2_{\ko_X}(\ke,\ke)\not=\nsp}.
\end{subsub}

\sepprop

\begin{subsub}\label{theo5}{\bf Theorem: } Every torsion free sheaf on $X$ can 
be deformed to linked torsion free sheaves.
\end{subsub}

\begin{proof} Let $\kf$ be a torsion free sheaf. We will use the exact sequence 
$(\ref{equ17})$:
\xmat{0\ar[r] & \kf\ar[r] & \bigoplus_{D\in\bf C}\kf_D\ar[r]^-\phi & 
\bigoplus_{P\in I}(\C_P\ot W_P)\ar[r] & 0 \ .}
For each intersection point \m{x\in I}, belonging to the components 
\m{D_1,D_2\in{\bf C}}, \m{\phi_x} is non zero only on 
\m{\kf_{D_1}\oplus\kf_{D_2}}, and comes from two surjective maps 
\m{\kf_{D_1,x}\to W_x}, \m{\kf_{D_2,x}\to W_x}. And $\kf$ is linked at $x$ if 
and only \m{\dim(W_x)=\inf(\rk(\kf_{D_1}),\rk(\kf_{D_2}))}.
Let \m{N_1\subset\kf_{D_1,x}} be a linear subspace, and 
\m{\ke=(\kf_{D_1})_{N_1}} (cf. \ref{el_mod}). Let \ \m{A_1=N_1\ot\ko_X(x)_x}.
Then we have an exact sequence
\[0\lra\kf_{D_1}\lra\ke(x)\lra\C_x\ot A_1\lra 0 \ , \]
associated to \ \m{\sigma\in\Ext^1_{\ko_X}(\C_x\ot A_1,\kf_{D_1})}. Let \ 
\m{\sigma'\in\Ext^1_{\ko_X}(\C_x\ot A_1,W_x)} \ be the image of $\sigma$ by 
$\phi$, and \ \m{0\to\C_x\ot W_x\to\km\lra\C_x\ot A_1\to 0} \ the associated 
extension, so that we have a commutative diagram with exact rows
\xmat{ & & 0\ar[d] & 0\ar[d]\\
0\ar[r] & \kg\ar[r]\fleq[d] & \kf_{D_1}\oplus\kf_{D_2}\ar[d]\ar[r]^-\phi & 
\C_x\ot W_x\ar[r]\ar[d] & 0\\
0\ar[r] & \kg\ar[r] & \ke(x)\oplus\kf_{D_2}\ar[r]\ar[d] & \km\ar[r]\ar[d] & 0\\
& & \C_x\ot A_1\fleq[r]\ar[d] & \C_x\ot A_1\ar[d]\\
& & 0 & 0}
where the sheaf $\kg$ is identical to $\kf$ in a neighbourhood of $x$. The 
sheaf $\km$ is concentrated on $x$ if and only if \ \m{\ko_X(-x)\km=0}, if and 
only if \ \m{\ko_X(-x)\ke(x)\subset\kg}. We have \ 
\m{\ko_X(-x)\ke(x)\subset\kf_{D_1}} \ and the image of \m{\ko_X(-x)\ke(x)} in 
\m{\kf_{D_1,x}} is precisely \m{N_1}. Hence $\km$ is concentrated on $x$ if and 
only if \ \m{N_1\subset\ker(\phi_{x|\kf_{D_1,x}})}. So it is possible to choose 
\m{N_1} such that $\km$ is concentrated on $x$ and of dimension 
\m{\inf(\rk(\kf_{D_1}),\rk(\kf_{D_2}))}. If we make such a modification at 
every \m{x\in I}, we obtain an exact sequence
\xmat{0\ar[r] & \kf\ar[r] & \bigoplus_{D\in\bf C} E_D\ar[r] & 
\bigoplus_{P\in I}(\C_P\ot V_P)\ar[r] & 0 \ ,}
where for every \m{D\in\bf C}, \m{E_D} is a vector bundle on $D$, and for every 
\m{x\in I}, \m{V_x} is a vector space of dimension 
\m{\inf(\rk(E_{D_1}),\rk(E_{D_2}))}. Now we consider the flat family of kernels 
of surjective morphisms
\[\bigoplus_{D\in\bf C}\E_D\lra\bigoplus_{P\in I}(\C_P\ot V_P) \ . \]
This family contains $\kf$ and its generic sheaf is linked.
\end{proof}

\sepprop

\begin{subsub}\label{rem3}{\bf Remark: }\rm In the preceding proof it would be 
possible to use more generally pairs of subspaces \m{N_1,N_2}, with 
\m{N_i\subset\kf_{D_i,x}}. In this case we would obtain vector bundles \m{\E_D} 
with the same ranks as the ones of the theorem but not necessarily the same
degrees.
\end{subsub}

\sepprop

\end{sub}

\sepsub

\Ssect{Moduli spaces of semi-stable sheaves}{mod_sh}

\begin{subsub}\label{polari} Polarisations -- \rm An ample line bundle 
\m{\ko_X(1)} on $X$ is defined by 
\begin{enumerate}
\item[--] for every $D\in\bf C$, a line bundle $\ko_D(1)$ of positive degree 
$\delta_D$ on $D$,
\item[--] for every intersection point $x\in I$, belonging to the components 
$D_1$ and $D_2$, an isomorphism \ $\ko_{D_1}(1)_x\simeq\ko_{D_2}(1)_x$.
\end{enumerate} 
Let $\kf$ be a torsion free sheaf. Consider the exact sequence $(\ref{equ17})$:
\xmat{0\ar[r] & \kf\ar[r] & \bigoplus_{D\in\bf C}\kf_D\ar[r]^-\phi & 
\bigoplus_{P\in I}(\C_P\ot W_P)\ar[r] & 0 \ , }
(recall that for every component \m{D\in\bf C}, \m{\kf_D=\kf_{|D}/T_D}, where 
\m{T_D} is the torsion subsheaf). For every  \m{D\in\bf C} and \m{P\in I}, let
\m{g_D} be the genus of $D$ and
\[r_D \ = \ \rk(\kf_D) \ , \quad d_D \ = \ \deg(\kf_D) \ , \quad h_P \ = \ 
\dim(W_P) \ . \]
Then the Hilbert polynomial of $\kf$ corresponding to \m{\ko_X(1)} is:
\begin{equation}\label{equ18} P_\kf(m) \ = \ \left(\sigg_{D\in\bf 
C}r_D\delta_D\right).m+\sigg_{D\in\bf C}(r_D(1-g_D)+d_D)-\sigg_{P\in I}h_P \ .
\end{equation}
\end{subsub}

\sepprop

\begin{subsub}\label{polari2}\rm
We will be mainly interested in the following case: $X$ has only two components
\m{D_1}, \m{D_2}, of the same genus $g$, intersecting in $d$ points. We will
also suppose that \ \m{\delta_{D_1}=\delta_{D_2}=\delta}. We have then
\[P_\kf(m) \ = \ (r_{D_1}+r_{D_2})\delta
m+(r_{D_1}+r_{D_2})(1-g)+d_{D_1}+d_{D_2}-\sigg_{P\in I}h_P \ . \]
If $\kf$ is a linked sheaf, we have \m{h_P=\inf(r_{D_1},r_{D_2})}
for every \m{P\in I}, hence
\[P_\kf(m) \ = \ (r_{D_1}+r_{D_2})\delta
m+(r_{D_1}+r_{D_2})(1-g)+d_{D_1}+d_{D_2}-d.\inf(r_{D_1},r_{D_2}) \ . \]
It follows that \m{P_\kf} depends only on
\m{r_{D_1}+r_{D_2}} \ and \ \m{d_{D_1}+d_{D_2}-d.\inf(r_{D_1},r_{D_2})}.
\end{subsub}

\sepprop

\begin{subsub}\label{comp_mod} Components of the moduli spaces of sheaves -- 
\rm Recall that \m{\km_X(\boldsymbol{P})} denotes the moduli space of sheaves
on $X$, stable with respect to \m{\ko_X(1)} and with Hilbert polynomial
$\boldsymbol{P}$. According to theorems \ref{lem1m} and \ref{theo5}
\begin{enumerate}
\item[--] Every component of $\km_X(\boldsymbol{P})$ contains points
corresponding to linked sheaves.
\item[--] If $\ke$ is a linked sheaf corresponding to a point in a component 
$\kn$ of $\km_X(\boldsymbol{P})$, for every linked sheaf $\kd$ corresponding to
a point of $\kn$, we have \ $\tau(\kf)=\tau(\ke)$.
\item[--] If $\ke$, $\kf$ are stable linked sheaves of Hilbert polynomial
$\boldsymbol{P}$ such that \ $\tau(\kf)=\tau(\ke)$, then their corresponding
points in $\km_X(\boldsymbol{P})$ belong to the same irreducible component.
\end{enumerate}
It follows that the components of \m{\km_X(\boldsymbol{P})} are indexed by the
types of the linked sheaves that they contain.

Using theorem \ref{theo4} and the same method as in the proof of theorem 
\ref{theo5}, it is easy to prove that every component of
\m{\km_X(\boldsymbol{P})} containing points corresponding to stable sheaves is 
generically smooth.

In the particular case of \ref{polari2}, $\boldsymbol{P}$ depends only on two
parameters $R$ and $D$: if $\ke$ is a linked semi-stable sheaf of Hilbert
polynomial $\boldsymbol{P}$, let \ \m{r_i=\rk(\ke_{|D_i}/T_i)} (where \m{T_i}
is the torsion subsheaf). Then \ \m{R=r_{D_1}+r_{D_2}},
\m{D=d_{D_1}+d_{D_2}-d.\inf(r_{D_1},r_{D_2})}.
We have seen that if \m{\km_X(\boldsymbol{P})} is non empty, then there exist
linked semi-stable sheaves of Hilbert polynomial $\boldsymbol{P}$, so we can use
the following notation: \ \m{\km_X(R,D)=\km_X(\boldsymbol{P})}.

If $R$ and $D$ are fixed, the components of \m{\km_X(R,D)} depend also on two
parameters: $r$ and $\epsilon$. If it contains the point corresponding to
$\ke$, take \m{r=r_{D_1}} and \m{\epsilon=d_{D_1}}. We will denote this
component by \m{\kn_X(R,D,r,\epsilon)}.
\end{subsub}

\end{sub}

\sepsec

\section{Coherent sheaves on blowing-ups of ribbons}\label{bl_up}

\Ssect{Blowing-ups of ribbons}{blp} 

(cf. \cite{ba_ei}, 1-, \cite{dr1}, 5.4).

Let $Y$ a primitive double curve with associated smooth curve $C$ and 
associated line bundle $L$. Let $D$ be a divisor on $C$ and
\[\phi:\widetilde{Y}\lra Y\]
the blowing-up of $Y$ along $D$. Then $\widetilde{Y}$ is a primitive double 
curve with associated smooth curve $C$ and associated line bundle \m{L(D)}.

If $\E$ is a vector bundle on $Y$, then \m{\phi^*(\E)} is a vector bundle on 
$\widetilde{Y}$, and we have \ \m{\phi^*(\E)_{|C}=\E_{|C}}. On the other hand, 
we have
\[\phi^*(\E)_1 \ \simeq \ \E_1(D) \ = \ \E_{|C}\ot L(D) \ . \]

Let $P$ be a point of $D$, with multiplicity $m$. Let \m{x\in\ko_{Y,P}} be over 
a generator of the maximal ideal of \m{\ko_{C,P}} and \m{z\in\ko_{Y,P}} an 
equation of $C$. Then there is a canonical isomorphism
\[\ko_{\widetilde{Y},P} \ \simeq \ \ko_{Y,P}\left[\frac{z}{x^m}\right] \ , \]
where \m{\dsp\ko_{Y,P}\left[\frac{z}{x^m}\right]} is the subring generated by 
\m{\ko_{Y,P}} and \m{\dsp\frac{z}{x^m}} in the localisation 
\m{\dsp\ko_{Y,P}\left[\frac{1}{x}\right]}. It follows easily that we have

\sepprop

\begin{subsub}\label{lem1i}{\bf Lemma: } Let $E$ be a vector bundle on $C$. 
Then we have
\[\phi^*(E) \ \simeq \ E[D] \ . \]
(cf. \ref{pu_th}). Let $\F$ be a vector bundle on $Y$, and \m{F=\F_{|C}}. Let 
\m{f:E\to\F} be a morphism, equivalent to a morphism \m{E\to F\ot L}. Then the 
morphism \m{E\to F\ot L(D)} induced by \m{\phi^*(f):\phi^*(E)\to\phi^*(\F)} is 
the composition
\[E\stackrel{f}{\lra} F\ot L\lra F\ot L(D) \ , \]
where the morphism on the right comes from the canonical section of 
\m{\ko_C(D)}.
\end{subsub}
\end{sub}

\sepsub

\Ssect{Inverse images of quasi-locally free sheaves}{blpsh}

Let $\ke$ be a quasi locally free sheaf on $Y$, and
\[0\lra K\lra\ke\stackrel{p}{\lra} E\lra 0\]
an exact sequence, where $K$ and $E$ are vector bundles on $C$. Then we have a 
canonical exact sequence
\[0\lra\ke_1\lra K\lra N\lra 0 \ , \]
where $N$ is the kernel of the canonical surjective morphism \ 
\m{\ke_{|C}\mapsto E}. Let $\sigma$ be the element of 
\m{\Ext^1_{\ko_C}(N,\ke_1)} associated to this exact sequence. Let
\[\widetilde{\ke} \ = \ \phi^*(\ke)/T(\phi^*(\ke)) \ . \]

\sepprop

\begin{subsub}\label{theo2} {\bf Theorem: } We have \ 
\m{(\widetilde{\ke})_{|C}\simeq\ke_{|C}},
\m{(\widetilde{\ke})_1\simeq\ke_1(D)}, and an exact sequence
\[0\lra\widetilde{K}\lra\widetilde{\ke}\stackrel{\widetilde{p}}{\lra} E\lra 0 \ 
, \]
where \m{\widetilde{K}} is the vector bundle on $C$ defined by the exact 
sequence
\[0\lra\ke_1(D)\lra\widetilde{K}\lra N\lra 0\]
associated to the element of \m{\Ext^1_{\ko_C}(N,\ke_1(D))} image of $\sigma$ 
by the morphism
\[\Ext^1_{\ko_C}(N,\ke_1)\lra\Ext^1_{\ko_C}(N,\ke_1(D))\]
induced by the canonical section of \m{\ko_C(D)}.
\end{subsub}

\begin{proof}
Let \m{U\subset C} be a non-empty open subset, \m{U_Y} (resp 
\m{U_{\widetilde{Y}}}) the corresponding open subset of $Y$ (resp. 
$\widetilde{Y}$), such that $\ke$ is split as
\[\ke_{|U_Y} \ \simeq \ \E\oplus F \ , \]
where $\E$ is locally free on \m{U_Y}, and $E$ is locally free on $U$. Then we 
have
\[\phi^*(\ke)_{|U_{\widetilde{Y}}} \ \simeq \ \phi^*(\E)\oplus F[D\cap U] \quad
\text{and} \quad \widetilde{\ke}_{|U_{\widetilde{Y}}} \ \simeq \ 
\phi^*(\E)\oplus F \ . \]
It follows that \ 
\m{(\widetilde{\ke}_{|U_{\widetilde{Y}}})_{|U}\simeq\ke_{|U}\simeq\E_{|U}\oplus 
F}, the first isomorphism being independent of the splitting \ 
\m{\ke_{|U_Y}\simeq\E\oplus F}. It follows that \ 
\m{(\widetilde{\ke})_{|C}\simeq\ke_{|C}}. We have also \ 
\m{(\widetilde{\ke})_{1|U_{\widetilde{Y}}}\simeq\ke_{1|U_Y}(D)=\E_1(D)}, the 
first isomorphism being independent of the splitting \ 
\m{\ke_{|U_Y}\simeq\E\oplus F}. It follows that \ 
\m{(\widetilde{\ke})_1\simeq\ke_1(D)}.

Let \ \m{p_C:\ke_{|C}\to E} \ be the morphism induced by $p$. The fact that $K$ 
is concentrated on $C$ is equivalent to \ \m{\ker(p_C)\subset p(\ke^{(1)})}. On 
$U$ this means that \ \m{N_{|U}=\ker(p_{C|U})\subset F}. The morphism \ 
\m{\widetilde{p}_C:\widetilde{\ke}_{|C}\to E} \ induced by $\widetilde{p}$ 
being the same as \m{p_C}, we see that \m{\ker(\widetilde{p})} is concentrated 
on $C$. We have in fact
\[K_{|U} \ = \ (\E_{|C}\ot L_{|U})\oplus N_{|U} \ , \quad \widetilde{K}_{|U} \ 
= \ (\E_{|C}\ot L(D)_{|U})\oplus N_{|U} \ . \]
It is easy to check that the obvious commutative diagram
\xmat{ & \ke_{1|U}\fleq[d] & K_{|U}\fleq[d] \\
0\ar[r] & \E_{|C}\ot L_{|U}\ar[r]\flinc[d] & (\E_{|C}\ot L_{|U})\oplus 
N_{|U}\ar[r]\flinc[d] & N_{|U}\ar[r]\fleq[d] & 0\\
0\ar[r] & \E_{|C}\ot L(D)_{|U}\ar[r]\fleq[d] & (\E_{|C}\ot L(D)_{|U})\oplus 
N_{|U}\ar[r]\fleq[d] & N_{|U}\ar[r] & 0\\
 & \ke_1(D)_{|U} & \widetilde{K}_{|U} 
}
is independent of the splitting \ \m{\ke_{|U_Y}\simeq\E\oplus F}. Hence all 
these diagrams can be glued to define this one on $C$
\xmat{0\ar[r] & \ke_1\ar[r]\flinc[d] & K\ar[r]\flinc[d] & N\ar[r]\fleq[d] & 0\\
0\ar[r] & \ke_1(D)\ar[r] & \widetilde{K}\ar[r] & N\ar[r] & 0}
This implies the last statement of the theorem by \cite{dr1b}, prop. 4.3.2. 
\end{proof}

\sepprop

Recall that we have an exact sequence
\xmat{0\ar[r] & \Ext^1_{\ko_C}(E,K)\ar[r]^-i &  
\Ext^1_{\ko_Y}(E,K)\ar[r]^-\theta & \Hom(E\ot L,K)\lra 0}
(cf. \ref{cstr}). On $\widetilde{Y}$ we have an exact sequence
\xmat{0\ar[r] & \Ext^1_{\ko_C}(E,\widetilde{K})\ar[r]^-{\widetilde{i}} &  
\Ext^1_{\ko_{\widetilde{Y}}}(E,\widetilde{K})\ar[r]^-{\widetilde{\theta}} & 
\Hom(E\ot L(D),\widetilde{K})\lra 0 .}
Let \ \m{0\to K\to\ke\to E\to 0} \ be an exact sequence on $Y$, and \
\m{\eta\in\Ext^1_{\ko_Y}(E,K)} \ corresponding to it. The morphism 
\m{\theta(\eta)} can be written as a composition
\xmat{E\ot L\ar[r]^-\tau & \ke_1\flinc[r] & K .}
Let \ \m{\widetilde{\eta}\in\Ext^1_{\ko_{\widetilde{Y}}}(E,\widetilde{K})} \ 
corresponding to the first exact sequence of theorem \ref{theo2}. Then 
\m{\widetilde{\theta}(\widetilde{\eta})} is the composition
\xmat{E\ot L(D)\ar[rr]^-{\tau\ot I_{\ko_C(D)}} & & \ke_1(D)\flinc[r] & 
\widetilde{K} .}

\sepsubsub

\begin{subsub}\label{cas_1can} The case of the first canonical filtration -- \rm
We consider the exact sequence on $Y$
\[0\lra\ke_1\lra\ke\lra\ke_{|C}\lra 0 \ , \]
(i.e. we suppose that \m{K=\ke_1} and \m{E=\ke_{|C}}) and \ 
\m{\eta\in\Ext^1_{\ko_Y}(E,K)} \ corresponding to it. In this case we have 
\m{N=0}, hence \ \m{\widetilde{K}=\ke_1(D)=K(D)}, we have an exact sequence on 
$\widetilde{Y}$
\[0\lra K(D)\lra\widetilde{\ke}\lra\widetilde{\ke}_{|C}=\ke_{|C}\lra 0\]
and \ \m{\widetilde{\theta}(\widetilde{\eta})=\theta(\eta)\ot I_{\ko_C(D)}}.
\end{subsub}

\sepprop

Now consider extensions
\[0\lra K\lra\ke\stackrel{p}{\lra} E\lra 0\]
corresponding to \ \m{\eta\in\Ext^1_{\ko_Y}(E,K)} \ such that \ 
\m{\theta(\eta):E\ot L\to K} \ is surjective. According to lemma \ref{lem1g} 
this condition is equivalent to the fact that \m{K=\ke_1} and \m{E=\ke_{|C}}.
Let \m{\Surj(E\ot L,K)} denote the set of surjective morphisms \m{E\ot 
L\to K}.

Let \ \m{0\to K\to\ke\to E\to 0} \ be an extension, corresponding to \ 
\m{\eta\in\Ext^1_{\ko_Y}(E,K)} \ such that \m{\theta(\eta)} is surjective. Let 
\m{E'\subset E} be a subbundle such that the restriction \m{E'\ot L\to K} of 
\m{\theta(\eta)} is surjective. Let \ \m{\eta'\in\Ext^1_{\ko_Y}(E',K)} \ be 
induced by $\eta$, and \ \m{0\to K\to\ke'\to E'\to 0} \ the corresponding exact 
sequence. Then the associated morphism \m{E'\ot L\to K} is the restriction of 
\m{\theta(\eta)}, and we have a commutative diagram on $Y$
\xmat{0\ar[r] & K\ar[r]\fleq[d] & \ke'\ar[r]\flinc[d] & E'\ar[r]\flinc[d] & 0\\
0\ar[r] & K\ar[r] & \ke\ar[r] & E\ar[r] & 0}
Similarly let \m{K\twoheadrightarrow K''} be a quotient bundle of $K$. Let \ 
\m{\eta''\in\Ext^1_{\ko_Y}(E,K'')} \ be induced by $\eta$, and \ \m{0\to 
K''\to\ke''\to E\to 0} \ the corresponding exact sequence. The associated 
morphism \m{E\ot L\to K''} is the composition
\xmat{E\ot L\ar[r]^-{\theta(\eta)} & K\flon[r] & K'' ,}
hence it is surjective, and we have a commutative diagram
\xmat{0\ar[r] & K\ar[r]\flon[d] & \ke\ar[r]\flon[d] & E\ar[r]\fleq[d] & 0\\
0\ar[r] & K''\ar[r] & \ke''\ar[r] & E\ar[r] & 0}

\sepprop

\begin{subsub}\label{lem1h}{\bf Lemma: } {\bf 1 -- } We have a commutative 
diagram on $\widetilde{Y}$
\xmat{0\ar[r] & K(D)\ar[r]\fleq[d] & \widetilde{\ke'}\ar[r]\flinc[d] &
E'\ar[r]\flinc[d] & 0\\
0\ar[r] & K(D)\ar[r] & \widetilde{\ke}\ar[r] & E\ar[r] & 0}
(where the two horizontal sequences are implied by theorem \ref{theo2}).

{\bf 2 -- } We have a commutative diagram on $\widetilde{Y}$
\xmat{0\ar[r] & K(D)\ar[r]\flon[d] & \widetilde{\ke}\ar[r]\flon[d] &
E\ar[r]\fleq[d] & 0\\
0\ar[r] & K''(D)\ar[r] & \widetilde{\ke''}\ar[r] & E\ar[r] & 0}
(where the two horizontal sequences are implied by theorem \ref{theo2}).
\end{subsub}

\begin{proof} We will prove only {\bf 1-} ({\bf 2-} is analogous). The exact 
sequences on $Y$ are coming from the first canonical filtrations, and the exact 
sequences of lemma \ref{lem1h} too. Hence the existence of the diagram comes 
from the functoriality of the canonical filtrations. The fact that the vertical 
morphisms are equality on the left and inclusions on the middle and the right 
can be seen by taking local descriptions of \m{\widetilde{\ke'}} and 
$\widetilde{\ke}$ as in theorem \ref{theo2}.
\end{proof}

\sepprop

\begin{subsub}\label{theo3} {\bf Theorem: } Suppose that \ \m{\Surj(E\ot 
L,K)\not=\emptyset}. Then there exists a linear map
\[\nu:\Ext^1_{\ko_Y}(E,K)\lra\Ext^1_{\ko_{\widetilde{Y}}}(E,K(D))\]
such that

{\bf 1 -- } For every \m{\eta\in\Ext^1_{\ko_Y}(E,K)} such that \m{\theta(\eta)} 
is surjective, if \ \m{0\to K\to\ke\to E\to 0} \ is the exact sequence 
corresponding to $\eta$, then \ \m{0\to K(D)\to\widetilde{\ke}\to E\to 0} \ is 
the exact sequence corresponding to \m{\nu(\eta)}.

{\bf 2 -- } We have a commutative diagram
\xmat{0\ar[r] & \Ext^1_{\ko_C}(E,K)\ar[r]^-i\ar[d]^\alpha &  
\Ext^1_{\ko_Y}(E,K)\ar[r]^-\theta\ar[d]^\nu & \Hom(E\ot L,K)\fleq[d]\lra 0\\
0\ar[r] & \Ext^1_{\ko_C}(E,K(D))\ar[r]^-{\widetilde{i}} &  
\Ext^1_{\ko_{\widetilde{Y}}}(E,K(D))\ar[r]^-{\widetilde{\theta}} & 
\Hom(E\ot L(D),K(D))\lra 0}
(where $\alpha$ is induced by the canonical section of \m{\ko_C(D)}).
\end{subsub}

\begin{proof} We can already define $\nu$ on the open subset of 
\m{\Ext^1_{\ko_Y}(E,K)} of elements $\eta$ such that \m{\theta(\eta)} is 
surjective. We will prove that it can be extended to a linear map.

Let \m{n\geq 2} be an integer and 
\m{\eta_1,\ldots,\eta_n\in\Ext^1_{\ko_Y}(E,K)} such that 
\m{\theta(\eta_1),\ldots,\theta(\eta_n)} are surjective, as well as 
\m{\theta(\eta_1+\cdots+\eta_n)}. For \m{1\leq i\leq n} let
\[0\lra K\lra\ke_i\lra E\lra 0\]
be the exact sequence corresponding to \m{\eta_i}.
We have a commutative diagram on $Y$ with exact horizontal sequences
\xmat{0\ar[r] & K\oplus\cdots\oplus K\ar[r]\flon[d]^f & 
\ke_1\oplus\cdots\oplus\ke_n\ar[r]\flon[d] & E\oplus\cdots\oplus 
E\ar[r]\fleq[d] & 0\\
0\ar[r] & K\ar[r]\fleq[d] & \kf\ar[r] & E\oplus\cdots\oplus E\ar[r] & 0\\
0\ar[r] & K\ar[r] & \ke\ar[r]\flinc[u] & E\ar[r]\flinc[u]_j & 0}
We can view \m{\Ext^1_{\ko_Y}(E\oplus\cdots\oplus E,K\oplus\cdots\oplus K)} as 
the space of \m{n\times n}-matrices with coefficients in 
\m{\Ext^1_{\ko_Y}(E,K)}. The exact sequence at the top is associated to the 
matrix 
with zero terms outside the diagonal, and with $i$-th diagonal term \m{\eta_i}. 
It is also the direct sum of all the exact sequences \ \m{0\to K\to\ke_i\to 
E\to 0}. The middle exact sequence is associated to \ 
\m{(\eta_1,\ldots,\eta_n)\in\Ext^1_{\ko_Y}(E\oplus\cdots\oplus E,K)}. The exact 
sequence down corresponds to \m{\eta_1+\cdots+\eta_n}. The morphisms $f$ and 
$j$ are the identity on each factor.

We will how see how this diagram is transformed on $\widetilde{Y}$ (using 
theorem \ref{theo2}). The first transformed exact sequence will obviously be 
the direct sum of the exact sequences \m{0\to K(D)\to\widetilde{\ke_i}\to E\to 
0}. According to lemma \ref{lem1h} the whole transformed diagram is as follows
\xmat{0\ar[r] & K(D)\oplus\cdots\oplus K(D)\ar[r]\flon[d]^{\widetilde{f}} & 
\widetilde{\ke_1}\oplus\cdots\oplus\widetilde{\ke_n}\ar[r]\flon[d] & 
E\oplus\cdots\oplus E\ar[r]\fleq[d] & 0\\
0\ar[r] & K(D)\ar[r]\fleq[d] & \widetilde{\kf}\ar[r] & E\oplus\cdots\oplus 
E\ar[r] & 0\\
0\ar[r] & K(D)\ar[r] & \widetilde{\ke}\ar[r]\flinc[u] & 
E\ar[r]\flinc[u]_{\widetilde{j}} & 0}
The morphisms $\widetilde{f}$ and $\widetilde{j}$ are the identity on each 
factor. It follows from this diagram that the exact sequence \m{0\to 
K(D)\to\widetilde{\ke}\to F\to 0} corresponds to 
\m{\nu(\eta_1)+\cdots\nu(\eta_n)}. Hence we have
\[\nu(\eta_1+\cdots+\eta_n) \ = \ \nu(\eta_1)+\cdots+\nu(\eta_n) \ . \]
Let \m{\eta\in\Ext^1_{\ko_Y}(E,K)} such that \m{\theta(\eta)} is surjective, and
\[0\lra K\lra\ke\stackrel{p}{\lra}E\lra 0\]
the associated exact sequence. Let \m{\lambda\in\C^*}. Then \m{\lambda\eta} 
corresponds to the exact sequence
\[0\lra K\lra\ke\stackrel{\frac{1}{\lambda}p}{\lra}E\lra 0 \ . \]
It follows easily that \ \m{\nu(\lambda\eta)=\lambda\nu(\eta)}.

Now let \m{(\eta_1,\ldots,\eta_n)} be a basis of \m{\Ext^1_{\ko_Y}(E,K)} such 
that for \m{1\leq i\leq n}, \m{\theta(\eta_i)} is surjective. We can now define 
coherently $\nu$ on the whole of \m{\Ext^1_{\ko_Y}(E,K)} by
\[\nu(\lambda_1\eta_1+\cdots+\lambda_n\eta_n) \ = \ \lambda_1\nu(\eta_1)+\cdots
+\lambda_n\nu(\eta_n) \ . \]
This proves {\bf 1-}.

It is clear that the right square of the diagram of {\bf 2-} is commutative, 
and we have to show that the left square is commutative. We will use 
\ref{cst1b}. Let \m{\eta\in\Ext^1_{\ko_Y}(E,K)} such that \m{\theta(\eta)} is 
surjective, and \ \m{0\to K\to\ke\to E\to 0} \ the corresponding extension.

The sheaf $\ke$ can be constructed using trivialisations as in the 
proof of \ref{theo2}: for every \m{P\in C} there exists a neighbourhood $U$ of 
$P$ such that on the corresponding open subset \m{U_Y} of $Y$ we have \ 
\m{\ke_{|U_{Y}}\simeq\E\oplus F}, where $\E$ is a vector bundle on $Y$ and $F$ 
a vector bundle on $C$. So $\ke$ can be constructed by gluing sheaves of type 
\m{\E\oplus F}. The gluing are isomorphisms
 \[f:\E\oplus F\stackrel{\simeq}{\lra}\E'\oplus F' \ , \]
where $\E$, \m{\E'} are vector bundles on \m{U_Y}, $F$, \m{F'} are vector 
bundles on $U$. Note that we have \ \m{\E_1=\E'_1=K_{U}}, hence 
\m{\E_{|U}=\E'_{|U}=(K\ot L^*)_{|U}}, and \ \m{\ke_U=\E_{|U}\oplus 
F=\E'_{|U}\oplus F'}. Suppose that on \m{U_Y}, $U$ is defined by an equation 
\m{z\in\ko_{Y}(U_Y)}. The isomorphism $f$ is defined by a matrix
\m{\begin{pmatrix}1+za & zb\\ c & d\end{pmatrix}}
where \m{a\in\End(K_{|U})}, \m{b\in\Hom(F,K_{|U})}, \m{c\in\Hom(K_{|U},F')},
\m{d\in\Hom(F,F')}. Let \m{u\in\Ext^1_{\ko_C}(E,K)}. According to \ref{cst1b}, 
the sheaf \m{\ke'} defined by \m{\eta+i(u)} can be constructed by the same 
trivialisations as $\ke$, but with gluings defined by matrices of type
\m{\begin{pmatrix}1+z(a+a') & z(b+b')\\ c & d\end{pmatrix}}
with \m{a'\in\End(K_{|U})}, \m{b'\in\Hom(F,K_{|U})}, \m{(a',b')} representing 
an element of \m{\Hom(E_{|U},K_{|U})}. The pairs \m{(a',b')} (for various open 
subsets $U$ covering $C$) make a cocycle representing $u$.

Now, according to lemma \ref{lem1i}, the sheaf $\widetilde{\ke}$ is obtained by 
gluing the sheaves \m{\phi^*(\E)\oplus F} on \m{U_{\widetilde{Y}}} (the open 
subset of $\widetilde{Y}$ corresponding to $U$), using the isomorphisms 
$\widetilde{f}$ defined the matrices
\m{\begin{pmatrix}1+zsa & zsb\\ c & d\end{pmatrix}},
where $s$ is the canonical section of \m{\ko_C(D)}. Similarly 
\m{\widetilde{\ke'}} is defined by the matrices
\m{\begin{pmatrix}1+zs(a+a') & zs(b+b')\\ c & d\end{pmatrix}}. It follows that 
\ \m{\nu(\eta+i(u))=\nu(\eta)+u'}, where \m{u'} is defined by a cocycle of 
pairs of type \m{(sa',sb')}, i.e. we have \m{u'=\alpha(u)}.
\end{proof}
\end{sub}

\sepsub

\Ssect{The dual construction}{dualx}

If instead of \ref{cas_1can} we use the second canonical filtration of a quasi 
locally free sheaf $\ke$:
\[0\lra\ke^{(1)}\lra\ke\lra\ke/\ke^{(1)}=\ke_1\ot L^*\lra 0 \ , \]
to describe the sheaf $\widetilde{\ke}$ using theorem \ref{theo2}, we find that 
the sheaf $\widetilde{K}$ does not depend only on \m{\ke_1\ot L^*} and 
\m{\ke^{(1)}}. We have an exact sequence \ \m{0\to\ke_1(D)\to\widetilde{K}\to 
N\to 0}, where $N$ and the extension may vary even if \m{\ke_1\ot L^*} and 
\m{\ke^{(1)}} remain constant.

To use the second canonical filtration, one can consider the following sheaf on 
$\widetilde{Y}$ associated to $\ke$:
\[\widehat{\ke} \ = \ \left(\widetilde{\ke^\vee}\right)^\vee\]
(cf. \ref{refle}). We have the exact sequence
\[0\lra(\ke_1\ot L^*)^\vee=\ke_1^*\ot 
L^2\lra\ke^\vee\lra\ke^{(1)\vee}=\ke^{(1)*}\ot L\lra 0 \ , \]
corresponding to the first canonical filtration of \m{\ke^\vee}. We have then 
the exact sequences on $\widetilde{Y}$:
\[0\lra\ke_1^*\ot L^2(D)\lra\widetilde{\ke^\vee}\lra\ke^{(1)*}\ot L\lra 0 \ , \]
corresponding to the first canonical filtration of $\widetilde{\ke^\vee}$, and
\[0\lra\ke^{(1)}(D)\lra\widehat{\ke}\lra\ke_1\ot L^*\lra 0 \ , \] 
corresponding to the second canonical filtration of $\widehat{\ke}$.

Now consider extensions
\[0\lra K\lra\ke\stackrel{p}{\lra} E\lra 0\]
corresponding to \ \m{\eta\in\Ext^1_{\ko_Y}(E,K)} \ such that \ 
\m{\theta(\eta):E\ot L\to K} \ is injective (as a morphism of vector bundles). 
According to lemma \ref{lem1g} this condition is equivalent to the fact that 
$\ke$ is quasi locally free, \m{K=\ke^{(1)}} and \m{E=\ke_1\ot L^*}. Let 
\m{\Inj(E\ot L,K)} denote the set of injective morphisms of vector bundles 
\m{E\ot L\to K}. The following result is similar to theorem \ref{theo3} and can 
be proved in the same way:

\sepprop

\begin{subsub}\label{theo3b} {\bf Theorem: } Suppose that \ \m{\Inj(E\ot 
L,K)\not=\emptyset}. Then there exists a linear map
\[\nu':\Ext^1_{\ko_Y}(E,K)\lra\Ext^1_{\ko_{\widetilde{Y}}}(E,K(D))\]
such that

{\bf 1 -- } For every \m{\eta\in\Ext^1_{\ko_Y}(E,K)} such that the morphism of 
vector bundles \m{\theta(\eta)} is \Nligne injective, if \ \m{0\to K\to\ke\to 
E\to 0} \ is the exact sequence corresponding to $\eta$, then \Nligne \m{0\to 
K(D)\to\widehat{\ke}\to E\to 0} \ is the exact sequence corresponding to 
\m{\nu(\eta)}.

{\bf 2 -- } We have a commutative diagram
\xmat{0\ar[r] & \Ext^1_{\ko_C}(E,K)\ar[r]^-i\ar[d]^\alpha &  
\Ext^1_{\ko_Y}(E,K)\ar[r]^-\theta\ar[d]^{\nu'} & \Hom(E\ot L,K)\fleq[d]\lra 0\\
0\ar[r] & \Ext^1_{\ko_C}(E,K(D))\ar[r]^-{\widetilde{i}} &  
\Ext^1_{\ko_{\widetilde{Y}}}(E,K(D))\ar[r]^-{\widetilde{\theta}} & 
\Hom(E\ot L(D),K(D))\lra 0}
(where $\alpha$ is induced by the canonical section of \m{\ko_C(D)}).
\end{subsub}

\sepprop

If \m{K=E\ot L}, then of course \m{\nu=\nu'}, and in this case the injective 
(or surjective) morphisms correspond to locally free $\ke$.

\end{sub}

\sepsec

\section{Coherent sheaves on reducible deformations of primitive double 
curves}\label{prim_reddef}

\Ssect{Preliminaries}{prel2}

(cf. \cite{dr7}, \cite{dr8})

Let $C$ be a projective irreducible smooth curve and \m{Y=C_2} a primitive 
double curve, with underlying smooth curve $C$, and associated line bundle $L$ 
on $C$. Let $S$ be a smooth curve, \m{P\in S} and \ \m{\pi:\kc\to S} \ a {\em 
maximal reducible deformation of \m{Y}} (cf. \cite{dr7}). This means that
\begin{enumerate}
\item[(i)] $\kc$ is a reduced algebraic variety with two irreducible components 
$\kc_1,\kc_2$.
\item[(ii)] We have \ $\pi^{-1}(P)=Y$. So we can view $C$ as a curve in
$\kc$.
\item[(iii)] For $i=1,2$, let \ $\pi_i:\kc_i\to S$ \ be the restriction
of $\pi$. Then \ $\pi_i^{-1}(P)=C$ \ and $\pi_i$ is a flat family of smooth
irreducible projective curves.
\item[(iv)] For every $z\in S\backslash\{P\}$, the components
$\kc_{1,z},\kc_{2,z}$ of $\kc_z$ meet transversally.
\end{enumerate}

For every \m{z\in S\backslash\{P\}}, \m{\kc_{1,z}} and \m{\kc_{2,z}} meet in 
exactly \m{-\deg(L)} points. If \ \m{\deg(L)=0}, then $\pi$ (or $\kc$) is 
called a {\em fragmented deformation}.

Let \m{\kz\subset\kc} be the closure in $\kc$ of the locus of the intersection
points of the components of \m{\pi^{-1}(z)}, \m{z\not=P}. Since $S$ is a curve,
$\kz$ is a curve of \m{\kc_1} and \m{\kc_2}. It intersects $C$ in a finite
number of points. If \m{x\in C}, let \m{r_x} be the number of
branches of $\kz$ at $x$ and \m{s_x} the sum of the multiplicities of the
intersections of these branches with $C$.

Let \m{t\in\ko_{S,P}} be a generator of the maximal ideal. We will also denote 
\m{\pi^*t} by $\pi$ and \m{\pi_i^*t} by \m{\pi_i}. So we have 
\m{\pi=(\pi_1,\pi_2)\in\ko_\kc(-C)}.

\sepprop

\begin{subsub} {\bf Theorem: }\label{theo1} (cf. \cite{dr8}, th. 3.1.1) Let 
\m{x\in C}. Then

{\bf 1 -- } There exists an unique integer \m{p>0} such that 
\m{\ki_{C,x}/\span{(\pi_1,\pi_2)}} is generated by the image of 
\m{(\pi_1^p\lambda_1,0)}, for some \m{\lambda_1\in\ko_{\kc_1,x}} not divisible 
by \m{\pi_1}. This integer does not depend on $x$.

{\bf 2 -- } $\lambda_1$ is unique up to multiplication by an invertible element
of \m{\ko_{\kc_1,x}}, and \m{(\pi_1^p\lambda_1,0)} is a generator of the ideal
\m{\ki_{\kc_1,\kc,x}} of \m{\kc_1} in $\kc$.

{\bf 3 -- } Let \m{m_x} be the multiplicity of \m{\lambda_{1\mid
C}\in\ko_{C,x}}. Then we have \ \m{m_x>0} \ if and only if \ \m{x\in\kz\cap
C}, and in this case we have \ \m{m_x=r_x=s_x}, and the branches of $\kz$ at
$x$ intersect transversally with $C$. Moreover
\[L \ \simeq \ \ko_C(-\sigg_{x\in \kz\cap C}r_xx) \ \simeq \ \ki_{\kz\cap C,C}
\ . \]

{\bf 4 -- } There exists \m{\lambda_2\in\ko_{\kc_2,x}} not divisible by 
\m{\pi_2}, such that \m{\ki_{C,x}/\span{(\pi_1,\pi_2)}} is generated by the 
image of \m{(0,\pi_2^p\lambda_2)}, and \m{\lambda_2} unique up to 
multiplication by an invertible element of \m{\ko_{\kc_2,x}}. 
\end{subsub}

\sepprop

In the preceding theorem, it is even possible to choose \m{\lambda_1}, 
\m{\lambda_2} such that \ \m{(\lambda_1,\lambda_2)\in\ko_{\kc,x}} (\cite{dr8}, 
corollary 3.1.2). The ideal sheaf of $\kz$ in \m{\kc_1} (resp. \m{\kc_2}) at 
$x$ is generated by \m{\lambda_1} (resp. \m{\lambda_2}).

Let
\[\kz_0 \ = \ \kc_1\cap\kc_2 \ \subset \kc \ . \]
We have then \ \m{(\kz_0)_{red}=\kz\cup C}. The ideal sheaf \ 
\m{\L_1=\ki_{Z_0,\kc_1}} (resp. \m{\L_2=\ki_{Z_0,\kc_2}}) of \m{\kz_0} in 
\m{\kc_1} (resp. \m{\kc_2}) at $x$ is generated by \m{\lambda_1\pi_1^p} (resp. 
\m{\lambda_2\pi_2^p}). Hence \m{\L_1} (resp. \m{\L_2}) is a line bundle on 
\m{\kc_1} (resp. \m{\kc_2}). The ideal sheaf  \ \m{\ki_{Z,\kc_1}} (resp. 
\m{\ki_{Z,\kc_2}}) of $\kz$ in \m{\kc_1} (resp. \m{\kc_2}) is canonically 
isomorphic to \m{\L_1} (resp. \m{\L_2}). The $p$-th infinitesimal 
neighbourhoods of $C$ in \m{\kc_1}, \m{\kc_2} (generated respectively by 
\m{\pi_1^p} and \m{\pi_2^p}) are canonically isomorphic, we will denote them by 
\m{C^{(p)}}. We have also a canonical isomorphism \ 
\m{\L_{1|C^{(p)}}\simeq\L_{2|C^{(p)}}}, and \ \m{\L_{1|C}\simeq\L_{2|C}\simeq 
L}. It is also possible, by replacing $S$ with a smaller neighbourhood of $P$, 
to assume that \ \m{\L_{1|\kz_0}\simeq\L_{2|\kz_0}}. Let \ 
\m{\L=\L_{1|\kz_0}=\L_{2|\kz_0}}.

We have \ \m{\ki_{\kc_1,\kc}=\L_2} \ and \ \m{\ki_{\kc_2,\kc}=\L_1}.

Let \m{\alpha\in\ko_{\kc_1,x}}. Then there exists \m{\beta\in\ko_{\kc_2,x}} 
such that \m{(\alpha,\beta)\in\ko_{\kc,x}}. By associating the class of $\beta$ 
to that of $\alpha$ we obtain an isomorphism of rings
\[\Phi_x:\ko_{\kc_1,x}/(\lambda_1\pi_1^p)\lra\ko_{\kc_2,x}/(\lambda_1\pi_2^p)\]
(the two are of course isomorphic to \m{\ko_{Z_0,x}}), and we have \ 
\m{\Phi_x(\pi_1)=\pi_2}.

\sepprop

\begin{subsub}\label{assoc} The associated fragmented deformation -- \rm 
The preceding isomorphism \m{\Phi_x} induces an isomorphism
\[\Psi_x:\ko_{\kc_1,x}/(\pi_1^p)\lra\ko_{\kc_2,x}/(\pi_2^p) \ . \]
These isomorphisms define a fragmented deformation \ \m{\rho:\kd\to S} \ with 
the same components \m{\kc_1}, \m{\kc_2} (which meet only along $C$ in $\kd$), 
and we have a canonical surjective morphism \ \m{\eta:\kd\to\kc}, induced by 
the obvious inclusion \m{\ko_\kc\subset\ko_\kd}. So we have
\[\ko_{\kc,x} \ = \ \left\{(\alpha,\beta)\in\ko_{\kc_1,x}\times\ko_{\kc_2,x}\ ; 
\ \alpha\equiv\beta\ \text{mod} \ \lambda_1\pi_1^p\right\}\]
and
\[\ko_{\kd,x} \ = \ \left\{(\alpha,\beta)\in\ko_{\kc_1,x}\times\ko_{\kc_2,x}\ ; 
\ \alpha\equiv\beta\ \text{mod} \ \pi_1^p\right\}\]
\end{subsub}

\sepprop

\begin{subsub}\label{lem3b}{\bf Lemma: } The primitive double curve 
\m{\rho^{-1}(P)} is the blowing up of \Nligne \m{\kz\cap C=\sigg_{x\in\kz\cap 
C}r_xx\subset Y}. 
\end{subsub}

\begin{proof} Let \ \m{C'_2=\rho^{-1}(P)}. It is clear that the morphism 
induced by $\rho$, \m{C'_2\to Y}, is an isomorphism between the 
open subsets corresponding to \m{C\backslash C\cap\kz}. Let \m{x\in C\backslash 
C\cap\kz}. Let \m{\lambda_1\in\ko_{\kc_1,x}} be a generator of the ideal of 
$\kz$ in \m{\kc_1}. Then the maximal ideal of \m{\ko_{C'_2,x}} (resp. 
\m{\ko_{Y,x}}) is generated by the image \m{z'} of \m{(\pi_1^p\lambda_1,0)}
(resp. by the image $z$ of \m{(\pi_1^p,0)}). If \m{t\in\ko_{\kc_1,x}} is over a 
generator of the maximal ideal of \m{\ko_{C,x}}, then we have \ 
\m{\lambda_1=\alpha t^{r_x}}, for some invertible \m{\alpha\in\ko_{\kc_1,x}}.
It follows easily that there is a canonical isomorphism
\[\ko_{C'_2,x} \ \simeq \ \ko_{Y,x}\left[\frac{z}{t^m}\right] \ , \]
where \m{\dsp\ko_{Y,x}\left[\frac{z}{z^m}\right]} is the subring generated by 
\m{\ko_{Y,z}} and \m{\dsp\frac{z}{t^m}} in the localisation 
\m{\dsp\ko_{Y,x}\left[\frac{1}{t}\right]}. This proves the lemma (cf. 
\ref{blp}).
\end{proof}

\sepprop

\begin{subsub}\label{prop4}{\bf Proposition: } The morphism \m{\eta:\kd\to\kc} 
is the blowing-up of $\kz$.
\end{subsub}

\begin{proof} It suffices to prove that $\eta$ satisfies the universal property 
of blowing-ups: if $X$ is an algebraic variety and \m{f:X\to\kc} is a morphism 
such that the ideal sheaf \m{f^{-1}(\ki_{\kz,\kc})\ko_X} is invertible, then 
there is a unique morphism \m{g:X\to\kd} such that \ \m{f=\eta\circ g}. The 
problem is local, we need only to prove that for every \m{x\in X}, there exists 
a neighbourhood $U$ of $x$ such that \m{f_{|U}} can be factorised through 
$\kd$, and this factorisation is unique. This is obvious if 
\m{f(x)\in\kc\backslash\kz}, and well known if \m{f(x)\in\kz\backslash C}. 
Suppose that \m{f(x)\in\kz\cap C}. For \m{i=1,2}, let 
\m{\lambda_i\in\ko_{\kc_i,x}} be a generator of the ideal of $\kz$ in 
\m{\kc_i}, such that \m{(\lambda_1,\lambda_2)\in\ko_{\kc,x}}. Then 
\m{\ki_{\kz,\kc,x}} is generated by \m{(\lambda_1,\lambda_2)} and 
\m{(\pi_1^p\lambda_1,0)}. Let \m{\gamma\in(f^{-1}(\ki_{\kz,\kc})\ko_X)_x} be a 
generator, it is not a zero divisor. Then there exist \m{a,b\in\ko_{X,x}} such 
that
\[\gamma \ = \ a.f^*(\lambda_1,\lambda_2)+b.f^*(\pi_1^p\lambda_1,0) \ . \]
There exist \m{\beta_1,\beta_2\in\ko_{X,x}} such that
\[f^*(\pi_1^p\lambda_1,0) \ = \beta_1\gamma \ , \quad
f^*(0,\pi_2^p\lambda_2) \ = \ \beta_2\gamma \ . \]
Since \ \m{(\pi_1^p\lambda_1,0)(0,\pi_2^p\lambda_2)=0}, we have \ 
\m{\beta_1\beta_2\gamma^2=0}, hence \ \m{\beta_1\beta_2=0}. Hence 
\m{\beta_1(x)=0} or \m{\beta_2(x)=0}. We can suppose that \m{\beta_1(x)=0}. We 
have
\[(1-\beta_1b)f^*(\pi_1^p\lambda_1,0) \ = \ \beta_1a.f^*(\lambda_1,\lambda_2) \ 
. \]
But\  \m{1-\beta_1b} \ is invertible, let \ 
\m{\dsp\rho=\frac{\beta_1a}{1-\beta_1b}}. We have then
\begin{equation}\label{equ5}f^*(\pi_1^p\lambda_1,0) \ = \ 
\rho.f^*(\lambda_1,\lambda_2)\end{equation}
and we can suppose that \ \m{\gamma=f^*(\lambda_1,\lambda_2)}, i.e. \m{a=1} and 
\m{b=0}. Now we describe the factorisation of \m{f^*} at $x$. According to 
\ref{assoc}, \m{\ko_{\kd,f(x)}} is generated by \m{\ko_{\kc,f(x)}} and 
\m{(\pi_1^p,0)}. Since \ 
\m{(\pi_1^p\lambda_1,0)=(\pi_1^p,0)(\lambda_1,\lambda_2)}, the equation 
(\ref{equ5}) suggests that we must take \ \m{f^*(\pi_1^p,0)=\rho}. This will 
define uniquely the extension of \m{f^*} provided that we show that there is no 
ambiguity in this definition, that is if $Q$ is a polynomial in one variable 
with coefficients in \m{\ko_{\kc,f(x)}} such that \ \m{Q((\pi_1^p,0))=0} \ in 
\m{\ko_{\kd,f(x)}}, then \ \m{Q(\rho)=0} in \m{\ko_{X,x}}. Let $m$ be the 
degree of $Q$. Then \ \m{(\lambda_1,\lambda_2)^mQ((\pi_1^p,0))\in\ko_{\kc,x}}.
We have
\[f^*\left((\lambda_1,\lambda_2)^mQ((\pi_1^p,0))\right) \ = 0 \ = \ 
f^*(\lambda_1,\lambda_2)^mQ(\rho) \ , \]
and since \m{f^*(\lambda_1,\lambda_2)} is not a zero divisor, we have 
\m{Q(\rho)=0}.
\end{proof}

\sepprop

\begin{subsub}\label{fami} The relative case -- \rm Let \m{S'} be a smooth 
connected curve and \m{P'\in S'}. Let \ \m{f:S'\to S} be a non constant 
morphism such that \m{f(P')=P}, and suppose that \ \m{f^{-1}(P)=\{P'\}} (which 
can be done by shrinking \m{S'}). Then \ \m{f^*(\kc)=\kc\times_S S'\to S'} \ is 
again a maximal deformation of $Y$. The integer $p$ remains the same if and 
only if the tangent map of $f$ at $P$ is not zero.
\end{subsub}

\end{sub}

\sepsub

\Ssect{Coherent sheaves on reducible deformations}{co_max_def}

We use the notations of \ref{prel2}, and we suppose that $L$ can be written as
\[L \ = \ \ko_C(-P_1-\cdots-P_d) \ , \]
where \m{P_1,\cdots,P_d} are distinct points of $C$, or that \m{L=\ko_C}.

\sepprop

\begin{subsub}\label{lo_fr_re} Locally free resolutions -- \rm For every \
\m{x\in\kz\cap C}, let \ 
\m{(\lambda_{1x},\lambda_{2x})\in\ko_{\kc,x}} \ such that 
\m{\ki_{C,x}/\span{(\pi_1,\pi_2)}} is generated by the image of 
\m{(\pi_1^p\lambda_1,0)}, and also by the image of \m{(0,\pi_2^p\lambda_2)}.
If \ \m{x\in C\backslash(\kz\cap C)}, let \ 
\m{(\lambda_{1x},\lambda_{2x})=(1,1)}.

Let $\D$ the ideal sheaf on $\kc$ such that
\begin{enumerate}
\item[--] $\D=\ko_\kc$ on \m{\kc\backslash(\kz\cap C)},
\item[--] For every $x\in\kz\cap C$, $\D_x=((\lambda_{1x},\lambda_{2x}))$.
\end{enumerate}
The sheaf $\D$ is well defined (because \m{\kz\cap C} is finite) and locally 
free. But it depends of the choice of \m{(\lambda_{1x},\lambda_{2x})} 
(\m{x\in\kz\cap C}).

We have a canonical morphism
\[r_1:\D\lra\ko_\kc\]
such that \ \m{\imm(r_1)=\ki_{\kc_2}}, which, for every \m{x\in C}, sends 
\m{(\alpha,\beta)\in\D_x} to \m{(\pi_1^p\alpha,0)}. And similarly, we have a 
canonical morphism
\[r_2:\D\lra\ko_\kc \ , \]
with image \m{\ki_{\kc_1}}.
\end{subsub}

It is easy to see that we have exact sequences
\begin{equation}\label{equ1}\xymatrix{\ldots\D^3\ar[rr]^-{r_1\ot I_{\D^2}} & & 
\D^2\ar[rr]^-{r_2\ot I_\D} & & \D\ar[r]^-{r_1} & \ko_\kc\ar[r] & 
\ko_{\kc_2}\ar[r] & 0 \ , }\end{equation} 
\begin{equation}\label{equ2}\xymatrix{\ldots\D^3\ar[rr]^-{r_2\ot I_{\D^2}} & & 
\D^2\ar[rr]^-{r_1\ot I_\D} & & \D\ar[r]^-{r_2} & \ko_\kc\ar[r] & 
\ko_{\kc_1}\ar[r] & 0 \ , }\end{equation}
i.e. locally free resolutions of \m{\ko_{\kc_1}} and \m{\ko_{\kc_2}}. It 
follows immediately that

\sepprop

\begin{subsub}\label{coro1}{\bf Corollary: } For \m{i=1,2}, we have \ 
\m{\EExt^1_{\ko_\kc}(\ko_{\kc_i},\ko_\kc)=0} \ and \
\m{\Tor^1_{\ko_\kc}(\ko_{\kc_1},\ko_{\kc_2})=0} .
\end{subsub}

\sepprop

\begin{subsub}\label{lem2}{\bf Lemma: } Let \m{x\in C}, and $M$ a torsion free
\m{\ko_{\kc_i,x}}-module (\m{i=1} or $2$). Then for \m{j\geq 1}, we have \ 
\m{\Tor^j_{\ko_{\kc,x}}(\ko_{Y,x},M)=0}.
\end{subsub}

\begin{proof}
This follows immediately from the free resolution of \m{\ko_Y}
\xmat{0\ar[r] & \ko_\kc\ar[r]^-\pi & \ko_\kc\ar[r] & \ko_Y\ar[r] & 0 \ .}
\end{proof}

\sepsubsub

\begin{subsub}\label{fl_sh} Extensions of vector bundles -- \rm
Let \m{E_i}, \m{i=1,2}, be a vector bundle on \m{\kc_i}. We are interested in 
extensions on $\kc$
\begin{equation}\label{equ3}0\lra E_1\lra\ke\lra E_2\lra 0\end{equation}
(of course the case of extensions \ \m{0\to E_2\to\ke\to E_1\to 0} \ is 
analogous). The sheaf $\ke$ is then flat on $S$.
\end{subsub}

\sepprop

\begin{subsub}\label{prop2}{\bf Proposition: } There are canonical isomorphisms
\[\EExt^1_{\ko_\kc}(E_2,E_1) \ \simeq \ \HHom(E_{2|\kz_0},\L_{1|\kz_0}^*\ot 
E_{1|\kz_0}) \ , \quad \Ext^1_{\ko_\kc}(E_2,E_1) \ \simeq \ 
\Hom(E_{2|\kz_0},\L^*\ot E_{1|\kz_0}) \ . \]
\end{subsub}

\begin{proof}
The first equality is an easy consequence of (\ref{equ1}). The second follows 
from the Ext spectral sequence (\cite{go}, 7.3) and the fact that \ 
\m{\HHom(E_2,E_1)=0}.
\end{proof}

\sepprop

Consider an extension (\ref{equ3}), corresponding to \ 
\m{\sigma\in\Hom(E_{2|\kz_0},\L_{1|\kz_0}^*\ot E_{1|\kz_0})}. Then according to 
lemma \ref{lem2} we have an exact sequence
\[0\lra E_{1|C}\lra\ke_{|Y}\lra E_{2|C}\lra 0\]
Recall that we have a canonical exact sequence
\xmat{0\ar[r] & \Ext^1_{\ko_C}(E_{2|C},E_{1|C})\ar[r]^-i & 
\Ext^1_{\ko_Y}(E_{2|C},E_{1|C})\ar[r]^-\theta & \Hom(E_{2|C}\ot L,E_{1|C})\ar[r]
& 0 . }

\sepprop

\begin{subsub}\label{lem1j}{\bf Lemma: } There are canonical morphisms
\[\tau:\Ext^1_{\ko_\kc}(E_2,E_1)\lra\Ext^1_{\ko_Y}(E_{2|C},E_{1|C}) \ , \quad
\ov{\tau}:\EExt^1_{\ko_\kc}(E_2,E_1)\lra\EExt^1_{\ko_Y}(E_{2|C},E_{1|C}) \ ,  \]
$\tau$ associating to the extension \ \m{0\to E_1\to\ke\to E_2\to 0} \ the 
extension 
\Nligne \m{0\to E_{1|C}\to\ke_{|Y}\to E_{2|C}\to 0}.
\end{subsub}

\begin{proof} Let \m{\ko_\kc(1)} be an ample line bundle on $\kc$. There exists 
a vector bundle \m{\F_0} on $\kc$ such that \ 
\m{\Ext^1_{\ko_\kc}(\F_0,E_1)=\nsp}, \m{\Ext^1_{\ko_Y}(\F_{0|Y},E_{1|C})=\nsp} 
\ and such that there exists a surjective morphism \m{f_0:\F_0\to E_2} (we can 
take \m{\F_0} of the form \ 
\m{\F_0=\ko_\kc(-n)\ot\C^k}, for suitable $k$ and $n$). Let \ 
\m{N_0=\ker(f_0)}, so we have an exact sequence
\[0\lra N_0\lra\F_0\lra E_2\lra 0 \ . \]
It follows from lemma \ref{lem2} that this sequence restricted to $Y$ is exact:
\[0\lra N_{0|Y}\lra\F_{0|Y}\lra E_{2|C}\lra 0 \ . \]
Hence we have a commutative diagram with exact rows:
\xmat{\Hom(\F_0,E_1)\ar[r]\ar[d] & \Hom(N_0,E_1)\ar[r]\ar[d] & 
\Ext^1_{\ko_\kc}(E_2,E_1)\ar[r] & 0\\
\Hom(\F_{0|Y},E_{1|C})\ar[r] & \Hom(N_{0|Y},E_{1|C})\ar[r] & 
\Ext^1_{\ko_Y}(E_{2|C},E_{1|C})\ar[r] & 0}
where the vertical arrows are the canonical morphisms. The morphism $\tau$ is 
then deduced immediately, and it is easy to see that it does not depend on 
\m{f_0}. The existence of $\ov{\tau}$ can be proved in the same way.
\end{proof}

\sepprop

We have a commutative diagram
\begin{equation}\label{equ4}
\xymatrix{\Ext^1_{\ko_\kc}(E_2,E_1)\fleq[r]\ar[d]^\tau & 
H^0(\EExt^1_{\ko_Y}(E_{2|C},E_{1|C}))
=\Hom(E_{2|\kz_0},\L_{1|\kz_0}^*\ot E_{1|\kz_0})\ar[d]^r\\
\Ext^1_{\ko_Y}(E_{2|C},E_{1|C})\ar[r]^-\theta & 
H^0(\EExt^1_{\ko_Y}(E_{2|C},E_{1|C}))=\Hom(E_{2|C}\ot L,E_{1|C})}
\end{equation}
where \ \m{r:\Hom(E_{2|\kz_0},\L_{1|\kz_0}^*\ot E_{1|\kz_0})\to
\Hom(E_{2|C}\ot L,E_{1|C})} \ is the restriction to $C$.

\sepsubsub

\begin{subsub}\label{dual_Y} Duality -- \rm If $\kf$ is a coherent sheaf on 
$\kc$, we will denote its dual \m{\Hom(\kf,\ko_\kc)} by \m{\kf^\vee}, and if 
$F$ is a coherent sheaf on \m{\kc_i}, its dual on \m{\kc_i} will be denoted by 
\m{F^*}. Note that, when we consider $F$ as a sheaf on $\kc$, we have
\[F^\vee \ = \ F^*\ot\L_i \ . \]
\end{subsub}

According to corollary \ref{coro1}, we have an exact sequence on $\kc$
\[0\lra E^*_2\ot\L_2\lra\ke^\vee\lra E^*_1\ot\L_1\lra 0 \ . \]
It follows from proposition \ref{prop2} that we have a canonical isomorphism
\[\Ext^1_{\ko_\kc}(E_1^*\ot\L_1,E_2^*\ot\L_2) \ \simeq \ 
\Ext^1_{\ko_\kc}(E_2,E_1) \ . \]
Of course, the elements of \m{\Ext^1_{\ko_\kc}(E_2,E_1)} corresponding to the 
exact sequence \m{(\ref{equ3})} and its dual are the same.

\sepprop

\begin{subsub}\label{polari3} Polarisations -- \rm Let \m{\ko_\kc(1)} be an
ample line bundle on $\kc$. As usually, if \m{Z\subset\kc} is a subvariety, we
will denote by \m{\ko_Z(1)} the restriction of \m{\ko_\kc(1)} to $Z$. For every
\m{s\in S\backslash\{P\}}, let \m{D_1}, \m{D_2} be the irreducible components
of \m{\kc_s}. Let \ \m{{\bf d}=\deg(\ko_C(1))}. Then since $C$ belongs to the
two families of curves \m{\kc_1}, \m{\kc_2}, we have 
\[\deg(\ko_{D_1}(1)) \ = \ {\bf d} \ = \ \deg(\ko_{D_2}(1)) \ . \]
Let $\ke$ be a coherent sheaf on $\kc$, flat on $S$. Suppose that for
every \m{s\in S}, \m{\ke_s} is torsion free. For \m{i=1,2}, let
\m{\ke_{s,i}=\ke_{s|D_i}/T_i}, where \m{T_i} is the torsion subsheaf. Then
\m{\ke_{s,i}} is a vector bundle on \m{D_i}. Let \m{r_i=\rk(\ke_{s,i})} and
\m{d_i=\deg(\ke_{s,i})}. From \ref{struc} we have an exact sequence
\[0\lra\ke_s\lra\ke_{s,1}\oplus\ke_{s,2}\lra\bigoplus_{x\in
D_1\cap D_2}\ko_{\kc_s,x}\ot W_x\lra 0 \ , \]
where for every \m{x\in D_1\cap D_2}, \m{W_x} is a finite dimensional vector
space. Then we have (cf. \ref{polari2})
\[P_{\ke_s}(m) \ = \ (r_1+r_2){\bf d}m +
(r_1+r_2)(1-g)+d_1+d_2-\sigg_{x\in D_1\cap D_2}\dim(W_x) \ . \]
On the other hand, we have (cf. \ref{inva})
\[P_{\ke_{|Y}}(m) \ = \ R(\ke_{|Y}){\bf d}m +R(\ke_{|Y})(1-g)+\Deg(\ke_{|Y}) \
. \]
Since \ \m{P_{\ke_s}=P_{\ke_{|Y}}}, we have
\[R(\ke_{|Y}) \ = \ r_1+r_2 \ , \quad \Deg(\ke_{|Y}) \ = \ d_1+d_2-\sigg_{x\in
D_1\cap D_2}\dim(W_x) \ . \]
\end{subsub}

\end{sub}

\sepsub

\Ssect{Coherent sheaves on fragmented deformations}{sh_frag}

We keep the notations of \ref{co_max_def}. If \ \m{\deg(L)=0}, we say that \ 
\m{\pi:\kc\to C} \ is a {\em fragmented deformation} of \m{Y}. Hence for 
every \m{s\in S\backslash\{P\}}, \m{\kc_s=\pi^{-1}(s)} is the disjoint union on 
\m{\kc_{1,s}} and \m{\kc_{2,s}}. In this case we have \ \m{L=\ko_C}.

Suppose that $\pi$ is a fragmented deformation and that \m{p=1}. In this case 
$\kc$ is simple, it is the gluing of \m{\kc_1} and \m{\kc_2} along $C$, i.e. 
for every \m{x\in C} we have
\[\ko_{\kc,x} \ = \ \left\{(\phi_1,\phi_2)\in\ko_{\kc_1,x}\times\ko_{\kc_2,x} \ 
; \ \phi_{1|C}=\phi_{2|C}\right\} \ . \]
We have \m{\kz_0=C}, hence the diagram \m{(\ref{equ4})}is
\xmat{\Ext^1_{\ko_\kc}(E_2,E_1)\fleq[r]\ar[d]^\tau & 
\Hom(E_{2|C},E_{1|C})\fleq[d]\\
\Ext^1_{\ko_Y}(E_{2|C},E_{1|C})\ar[r]^-\theta & \Hom(E_{2|C},E_{1|C})}
It follows that

\sepprop

\begin{subsub}\label{lem1k}{\bf Proposition: } Given \ \m{f:E_{2|C}\to 
E_{1|C}}, there exists only one extension on $Y$
\[0\lra E_{1|C}\lra\kf\lra E_{2|C}\lra 0\]
corresponding to \ \m{\sigma\in\Ext^1_{\ko_Y}(E_{2|C},E_{1|C})} such that 
\m{\theta(\sigma)=f}, and which is the restriction to $Y$ on an extension 
\m{(\ref{equ3})} on $\kc$.
\end{subsub}

\sepprop

\begin{subsub}\label{geom_r} The case of vector bundles -- \rm The sheaf $\ke$ 
is locally free if and only $f$ is an isomorphism. In this case $\ke$ is 
obtained by gluing \m{E_1} and \m{E_2} along \m{E_{1|C}\simeq E_{2|C}}, so it 
is unique.
\end{subsub}

\end{sub}

\sepsub

\Ssect{Regular sheaves on reducible deformations}{reg_sh}

\begin{subsub}\label{defin_reg}{\bf Definition:} A coherent sheaf $\ke$ on 
$\kc$ is called {\em regular} if it is locally free on \m{\kc\backslash\kz_0}, 
and if for every \m{x\in\kz_0} there exists a neighbourhood of $x$ in $\kc$, a 
vector bundle $\E$ on $U$, \m{i\in\{1,2\}}, and a vector bundle $F$ on 
\m{U\cap\kc_i}, such that \ \m{\ke_{|U}\simeq\E\oplus F}.
\end{subsub}

\sepprop

\begin{subsub}\label{lem3} {\bf Lemma: } In the preceding definition, $i$, and 
the ranks of $\E$ and $F$ are unique.
\end{subsub}

\begin{proof} Suppose that \ \m{\ke_{|U}\simeq\E'\oplus F'}, with \m{\E'} 
locally free on $U$ and \m{F'} a vector bundle on \m{U\cap\kc_j}. The induced 
isomorphism \ \m{\E\oplus F\to\E'\oplus F'} \ is defined by a matrix \ 
\m{\begin{pmatrix}A & B\\C & D\end{pmatrix}}. We have \m{B_{|\kz_0}=0}, hence 
\m{A_{|\kz_0}} and \m{D_{|\kz_0}} are isomorphisms. The result follows 
immediately.
\end{proof}

\sepprop

\begin{subsub}\label{prop3}{\bf Proposition: } Let $\ke$ be a coherent sheaf on 
$\kc$. Then the following assertions are equivalent:

{\em (i)} $\ke$ is regular (with $i=1$ in definition \ref{defin_reg}).

{\em (ii)} There exists an exact sequence \ \m{0\to E_2\to\ke\to E_1\to 0}, 
where for \m{j=1,2}, \m{E_j} is a vector bundle on \m{\kc_j}, such that the 
associated morphism \ \m{E_{1|\kz_0}\to\L^*\ot E_{2|\kz_0}} \ is surjective on
a neighbourhood of $C$.

{\em (iii)} There exists an exact sequence \ \m{0\to E_1\to\ke\to E_2\to 0}, 
where for \m{j=1,2}, \m{E_j} is a vector bundle on \m{\kc_j}, such that the 
associated morphism \ \m{E_{2|\kz_0}\to\L^*\ot E_{1|\kz_0}} \ is injective (as 
a morphism of vector bundles) on a neighbourhood of $C$.

We have a similar result by taking \m{i=2} in definition
\ref{defin_reg}.\end{subsub}

\begin{proof} Suppose that (i) is true. Let \m{x\in\kz_0} and $U$ a 
neighbourhood of $x$ such that \ \m{\ke_{|U}\simeq\E\oplus F}, with $\E$ 
locally free on $U$ and $F$ a vector bundle on \m{U\cap\kc_1}. Then we have \ 
\m{\ke_{|U\cap\kc_1}\simeq\E_{|U\cap\kc_1}\oplus F}. Hence \ 
\m{E_1=\ke_{|\kc_1}} \ is a vector bundle on $\kc_1$. The kernel of the 
restriction morphism \m{\ke\to\ke_{|\kc_1}} is \m{\E\ot\L_2}, which is a vector 
bundle on  \m{U\cap\kc_2}. Hence the kernel \m{E_2} of the restriction morphism 
\m{\ke\to E_1} is a vector bundle on \m{\kc_2}. On $U$, the associated morphism 
\ \m{E_{1|\kz_0}\to\L^*\ot E_{2|\kz_0}} \ is the projection
\[\E_{|U\cap\kz_0}\oplus F_{|U\cap\kz_0}\lra\E_{|U\cap\kz_0} \ , \]
hence it is surjective. This proves (ii). The proof of (iii) is similar, in 
this case \m{E_2} is \m{\ke_{\kc_2}} quotiented by its torsion subsheaf.

Suppose that (ii) is true. Let Let \m{x\in\kz_0} and $U$ a 
neighbourhood of $x$ such that \m{E_{1|U\cap\kc_1}}, \m{E_{2|U\cap\kc_2}} and 
\m{\L_{|\kz_0\cap U}} are trivial, say \ 
\m{E_{1|U\cap\kc_1}\simeq\ko_{U\cap\kc_1}\ot\C^{r_1}}, 
\m{E_{2|U\cap\kc_2}\simeq\ko_{U\cap\kc_2}\ot\C^{r_2}}. We can also assume that 
the associated surjective morphism  \ \m{E_{1|\kz_0\cap U}\to\L^*\ot 
E_{2|\kz_0\cap U}} \ is the projection
\[\ko_{U\cap\kz_0}\ot\C^{r_1}=(\ko_{U\cap\kz_0}\ot\C^{r_1-r_2})\oplus
(\ko_{U\cap\kz_0}\ot\C^{r_2})\lra\ko_{U\cap\kz_0}\ot\C^{r_2} \ . \]
The extension corresponding to this morphism (cf. prop. \ref{prop2}) is the
direct sum of
\[0\lra\ko_{U\cap\kc_2}\ot\C^{r_2}\lra\ko_U\ot\C^{r_2}\lra
\ko_{U\cap\kc_1}\ot\C^{r_2}\lra 0\]
and
\[0\lra 0\lra\ko_{U\cap\kc_1}\ot\C^{r_1-r_2}\lra\ko_{U\cap\kc_1}\ot\C^{r_1-r_2}
\lra 0 \ , \]
hence \ \m{\ke_{|U}\simeq(\ko_U\ot\C^{r_2})\oplus(\ko_{U\cap\kc_1}
\ot\C^{r_1-r_2})}, and (i) is proved. The proof that (iii) implies (i) is 
similar.
\end{proof}

\sepprop

\begin{subsub}\label{geom_co} Geometric construction of regular sheaves on 
fragmented deformations -- \rm We suppose now that $\kc$ is a fragmented 
deformation as in \ref{sh_frag}, with \m{p=1}. Let \m{E_i}, \m{i=1,2}, a vector
bundle on \m{\kc_i}, and \ \m{\phi:E_{1|C}^*\to E_{2|C}^*} \ an injective
morphism of vector bundles. Viewing \m{E_1^*} and \m{E_2^*} as algebraic
varieties, let $Y$ be the variety obtained by gluing \m{E_{1|C}^*\subset E_1^*}
and the corresponding subbundle of \m{E_{2|C}^*\subset E_2^*}. The existence of
$Y$ is easily seen for example by using th\'eor\`eme 5.4 of \cite{ferr}. We have
a canonical projection \ \m{\rho:Y\to\kc}, whose fibres are vector spaces. Let
$\ke$ be the sheaf of \m{\ko_\kc}-modules defined by: for every open subset
\m{U\subset\kc}, \m{\ke(U)} is the subset of \m{\ko_Y(\rho^{-1}(U))} of regular
maps which are linear on the fibres of $\rho$. By considering open subsets on
which \m{E_1} and \m{E_2} are trivial, it is easy to see that $\ke$ is a regular
sheaf, and that we have an exact sequence
\[0\lra E_1\lra\ke\lra E_2\lra 0\]
associated with the surjective morphism \ \m{^t\phi:E_{2|C}\to E_{1|C}}.
\end{subsub}

\sepprop
%\newpage

\begin{subsub}\label{prop14}{\bf Proposition: } Let $\ke$ be a coherent sheaf
on $\kc$, flat on $S$. Suppose that for every \m{s\in S}, \m{\ke_s} is torsion
free, and that \m{\ke_{|Y}} is quasi locally free of rigid type (cf.
\ref{rig_t}). Then $\ke$ is regular.
\end{subsub}

\begin{proof} We will suppose that \m{\ke_{|Y}} is locally isomorphic to \
\m{a\ko_Y\oplus\ko_C} \ for some integer \m{a\geq 0} (the case where
\m{\ke_{|Y}} is locally free is similar and easier). Let \m{s\in
S\backslash\{P\}}, \m{D_1}, \m{D_2} the irreducible components of \m{\kc_s},
and for \m{i=1,2}, \m{F_i=\ke_{s|D_i}/T_i} (where \m{T_i} is the torsion
subsheaf), it is a vector bundle on \m{D_i}. Let \ \m{r_i=\rk(F_i)}. Then from
\ref{polari3} we have
\[R(\ke_{|Y}) \ = \ 2a+1 \ = \ r_1+r_2 \ . \]
For every \m{x\in C}, we have \ \m{\rk(\ke_{|Y,x})=a+1}. By semi-continuity of
the rank, we have \m{r_1\leq a+1} and \m{r_2\leq a+1}. If follows that
\m{r_1=a}, \m{r_2=a+1} or \m{r_1=a+1}, \m{r_2=a}. We can suppose that we are in
the first case.

Let \m{Z_1\subset\kz} be an irreducible component. It is a smooth
irreducible curve meeting $C$ in one point $z$. We have \
\m{\rk(\ke_{Y,z})=a+1}, hence for a general point \m{z'\in Z_1}, we have \
\m{\rk(\ke_{|Z_1,z'})\leq a+1}. Suppose that \m{\pi(z')=s\not=P}. Then from
\ref{prim_red}-, there exist integers \m{a_1\geq 0}, \m{a_2\geq 0}, \m{b\geq 0}
such that
\[\ke_{s,z'} \ \simeq \ a_1\ko_{D_1,z'}\oplus a_2\ko_{D_2,z'}\oplus
b\ko_{\kc_s,z'} \ . \]
Hence we have \ \m{a_1+a_2+b=\rk(\ke_{|Z_1,z'})\leq a+1}, \m{a_1+b=r_1=a},
\m{a_2+b=r_2=a+1}. It follows that \m{a_1=0}. Hence \m{\ke_s} is linked at
\m{z'}. So we can suppose that there exists a neighbourhood $U$ of $P$ in $S$
such that for every \m{s\in S\backslash\{P\}}, \m{\ke_s} is linked. From now
on, we suppose that \m{U=S}.

The sheaf \m{E_2=\ke_{\kc_2}} is of rank \m{a+1} at every point of \m{\kc_2},
hence it is a vector bundle on \m{\kc_2}. Let \m{E_1} be the kernel of the
restriction morphism \ \m{\ke\to E_2}. It is a torsion free sheaf concentrated
on \m{\kc_1}. From lemma \ref{lem2}, we have an exact sequence
\[0\lra E_{1|C}\lra\ke_{|Y}\lra E_{2|C}\lra 0 \ . \]
Hence we have \ \m{E_{2|C}=\ke_C} \ and \m{E_{1|C}=(\ke_{|Y})_1}. In particular
\m{E_{1|C}} is of rank $a$ on $C$, as on \m{\kc_1\backslash C}. It follows that
\m{E_1} is a vector bundle on \m{\kc_1}. The morphism \
\m{E_{2|\kz_0}\to\L^*\ot E_{1|\kz_0}} \ corresponding to the exact sequence \
\m{0\to E_1\to\ke\to E_2\to 0} \ is surjective, hence $\ke$ is regular by
proposition \ref{prop3}.
\end{proof}

\sepprop

\end{sub}

\sepprop

\Ssect{Regular sheaves and blowing-ups}{regbl}

We keep the notations of \ref{co_max_def}, and we suppose that \m{p=1}.

We now consider the associated fragmented deformation \m{\rho:\kd\to S} (cf. 
\ref{assoc}). Recall that the canonical morphism \m{\eta:\kd\to\kc} is the 
blowing-up of $\kz$. If \m{C'_2=\rho^{-1}(P)}, the induced morphism 
\m{C'_2\to Y} is the blowing-up of \m{\kz\cap C}.

Let $\ke$ be a regular sheaf on $\kc$, \m{T(\eta^*(\ke))} the torsion subsheaf 
of \m{\eta^*(\ke)}, and
\[\widetilde{\ke} \ = \ \eta^*(\ke)/T(\eta^*(\ke)) \ . \]

\sepprop

\begin{subsub}\label{prop5}{\bf Proposition: } Let $\ke$ be a regular sheaf on 
$\kc$, and \ \m{0\to E_2\to\ke\to E_1\to 0} \ an exact sequence, where \m{E_i} 
is a vector bundle on \m{\kc_i}, for \m{i=1,2}, such that the associated 
morphism \ \m{\lambda:E_{1|\kz_0}\to\L^*\ot E_{2|\kz_0}} \ is surjective. Then 
we have an exact sequence on $\kd$
\[0\lra E_2\ot\L_2^*\lra\widetilde{\ke}\lra E_1\lra 0\]
such that the associated morphism \ \m{E_{1|C}\to E_{2|C}\ot L^*} \ is the 
restriction of $\lambda$.
\end{subsub}

\begin{proof} The proof is easy, using the description of $\rho$ in \ref{assoc} 
and the local description of $\ke$ given in the proof of proposition 
\ref{prop3}.
\end{proof}

\sepprop

It is clear that the dual sheaf \m{\ke^\vee} is also regular. Let 
\[\widehat{\ke} \ = \ \left(\widetilde{\ke^\vee}\right)^\vee \ . \]
The proof of the following result is similar to that of proposition \ref{prop5}:

\sepprop

\begin{subsub}\label{prop6}{\bf Proposition: } Let $\ke$ be a regular sheaf on 
$\kc$, and \ \m{0\to E_2\to\ke\to E_1\to 0} \ an exact sequence, where \m{E_i} 
is a vector bundle on \m{\kc_i}, for \m{i=1,2}, such that the associated 
morphism of vector bundles \ \m{\lambda:E_{1|\kz_0}\to\L^*\ot E_{2|\kz_0}} \ is 
injective. Then we have an exact sequence on $\kd$
\[0\lra E_2\ot\L_2^*\lra\widehat{\ke}\lra E_1\lra 0\]
such that the associated morphism \ \m{E_{1|C}\to E_{2|C}\ot L^*} \ is the 
restriction of $\lambda$.
\end{subsub}

\end{sub}

\sepsub

\Ssect{Deformations of regular sheaves on fragmented deformations}{def_reg_f}

We keep the notations of \ref{regbl}. Let $\ke$, \m{\ke'} be regular sheaves on 
$\kd$. We suppose that there are exact sequences
\[0\lra F\lra\ke\lra E\lra 0 \ , \quad 0\lra F\lra\ke'\lra E'\lra 0 \ , \]
where $F$ is a vector bundle on \m{\ke_2}, $E$, \m{E'} are vector bundles on 
\m{\kc_1} such that \ \m{E_{|C}=E'_{|C}=E_0}. We can then define the Koda\"\i 
ra-Spencer morphism
\[\omega_{E,E'}:TS_P\lra\Ext^1_{\ko_C}(E_0,E_0)\]
(cf. \ref{Koda}). We suppose also that the two morphisms \m{E_{|C}\to F_{|C}},
\m{E'_{|C}\to F_{|C}} associated to these exact sequences are the same and 
surjective, there are denoted by \ \m{\phi:E_0\to F_{|C}}. 

Let
\[0\lra F_{|C}\lra\ke_{|C'_2}\lra E_0\lra 0 \ , \quad 0\lra 
F_{|C}\lra\ke'_{|C'_2}\lra E_0\lra 0\]
be the restrictions of the preceding exact sequences, and \ \m{\sigma,\sigma'\in
\Ext^1_{\ko_{C'_2}}(E_0,F_{|C})} \ the corresponding elements. 
Let \m{\tau\in TS_P} be the element corresponding to the isomorphism \ 
\m{\ki_{C,C'_2}\simeq\ko_C}. Let
\xmat{0\ar[r] & \Ext^1_{\ko_C}(E_0,F_{|C})\ar[r]^-\iota & 
\Ext^1_{\ko_{C'_2}}(E_0,F_{|C})\ar[r]^-\theta & \Hom(E_0,F_{|C})\ar[r] & 0}
be the canonical exact sequence (cf. \ref{cstr}). Let
\[\ov{\phi}:\Ext^1_{\ko_C}(E_0,E_0)\lra\Ext^1_{\ko_C}(E_0,F_{|C})\]
be the map induced by $\phi$. 

\sepprop
%\newpage

\begin{subsub}\label{prop7}{\bf Proposition: } We have \ \m{\dsp\sigma'-
\sigma=\iota\left(\ov{\phi}(\omega_{E,E'}(\tau))\right)}.
\end{subsub}

\begin{proof} We will use \v Cech cohomology.
According to \ref{reg_sh}, there exists an open cover \m{(U_i)_{i\in I}} of a 
neighbourhood of $C$ in $\kd$ such that for every \m{i\in I} we have 
isomorphisms
\xmat{\ke_{|U_i}\ar[rrd]^{\lambda_i}\\
& & \V_i\oplus W_i\\\ke'_{|U_i}\ar[rru]^{\lambda'_i}}
where \m{\V_i}, \m{W_i} are trivial vector bundles on \m{U_i}, \m{U_i\cap\kc_1} 
respectively. We can suppose that the induced isomorphisms \ \m{(\V_i\oplus 
W_i)_{|C}\simeq E_C} \ and \ \m{(\V_i\oplus W_i)_{|C}\simeq E'_C} \ are the 
same, as well as the two \ \m{F_{i|U_i}\simeq\V_i\ot\ki_{\kc_1}}. Now let
\[\theta_{ij}=\lambda_j\circ\lambda_i^{-1}, \
\theta'_{ij}=\lambda'_j\circ\lambda^{'-1}_i:
(\V_i\oplus W_i)_{|U_{ij}}\to(\V_j\oplus W_j)_{|U_{ij}} \ ,\]
represented respectively by matrices
\[M_{ij} \ = \ \begin{pmatrix}A_{ij} & B_{ij}\\ C_{ij} & D_{ij}\end{pmatrix} 
\quad \text{and} \quad
M'_{ij} \ = \ \begin{pmatrix}A'_{ij} & B'_{ij}\\ C'_{ij} & D'_{ij}\end{pmatrix} 
\ . \]
We can write
\[A'_{ij} \ = \ A_{ij}+\alpha_{ij} \ , \]
where
\[\alpha_{ij}:\V_{i|U_{ij}\cap\kc_1}\to(\pi_1,0)\V_{j|U_{ij}}
=(\pi_1)\V_{j|U_{ij}\cap\kc_1} \ , \]
and the two matrices restricted to $C$ are the same. Now we have induced 
trivialisations of $E$ and $E'$: \m{\mu_i:E_{U_i}\to\V_{i|\kc_1}\oplus F_i},  
\m{\mu'_i:\ke'_{U_i}\to\V_{i|\kc_1}\oplus F_i}, such that
\[\nu_{ij}=\mu_j\circ\mu_i^{-1}, \ \nu'_{ij}=\mu'_j\circ\mu^{'-1}_i:
\V_{i|U_{ij}\cap\kc_1}\oplus W_{i|U_{ij}}\to
\V_{j|U_{ij}\cap\kc_1}\oplus W_{j|U_{ij}}\]
are represented by the matrices \m{M_{ij|\kc_1}} and  \m{M'_{ij|\kc_1}} 
respectively. It follows that \m{C_{ij}}, \m{C'_{ij}} are multiples of 
\m{\pi_1}, and that \m{D_{ij}}, \m{D'_{ij}} are of the form \ 
\m{D_{ij}=I+\pi_1\delta_{ij}}, \m{D'_{ij}=I+\pi_1\delta'_{ij}}. We have then on 
\m{C'_2}
\[(M_{ij}-M'_{ij})_{|C'_2} \ = \ \begin{pmatrix}\alpha_{ij|C'_2} & 
(B_{ij}-B'_{ij})_{|C'_2}\\ 0 & 0\end{pmatrix}\]
and a similar formula for \m{(M_{ij}-M'_{ij})_{|C}}.

The result follows then easily from the interpretation in \ref{Koda_4} of 
the Koda\"\i ra-Spencer morphism and from \ref{cst1b}. \end{proof}

\sepprop

We suppose that there are exact sequences
\[0\lra F\lra\ke\lra E\lra 0 \ , \quad 0\lra F'\lra\ke'\lra E\lra 0 \ , \]
where $E$ is a vector bundle on \m{\ke_1}, $F$, \m{F'} are vector bundles on 
\m{\kc_2} such that \ \m{F_{|C}=F'_{|C}=F_0}. We can then define the Koda\"\i 
ra-Spencer morphism
\[\omega_{E,E'}:TS_P\lra\Ext^1_{\ko_C}(E,E) \ . \]
We suppose also that the two morphisms \m{E_{|C}\to F_{|C}}, \m{E_{|C}\to 
F'_{|C}} associated to these exact sequences are the same and surjective, there 
are denoted by \ \m{\psi:E_{|C}\to F_0}. Let
\[0\lra F_0\lra\ke_{|C'_2}\lra E_{|C}\lra 0 \ , \quad 0\lra 
F_0\lra\ke'_{|C'_2}\lra E_{|C}\lra 0\]
be the restrictions of the preceding exact sequences, and \ \m{\rho,\rho'\in
\Ext^1_{\ko_{C'_2}}(E_{|C},F_0)} \ the corresponding elements. Let
\xmat{0\ar[r] & \Ext^1_{\ko_C}(E_{|C},F_0)\ar[r]^-\kappa & 
\Ext^1_{\ko_{C'_2}}(E_{|C},F_0)\ar[r]^-\zeta & \Hom(E_{|C},F_0)\ar[r] & 0}
be the canonical exact sequence (cf. \ref{cstr}). Let
\[\widetilde{\psi}:\Ext^1_{\ko_C}(F_0,F_0)\lra\Ext^1_{\ko_C}(E_{|C},F_0)\]
be the map induced by $\psi$. The following result is similar to proposition 
\ref{prop7}

\sepprop

\begin{subsub}\label{prop8}{\bf Proposition: } We have \ \m{\dsp\rho'-
\rho=\kappa\left(\widetilde{\psi}(\omega_{F,F'}(\tau))\right)}.
\end{subsub}

\sepprop

\begin{subsub}\label{coro2}{\bf Corollary:} Let $\E$ be a quasi locally free 
sheaf on $Y$. Then there exist a smooth curve $T$, \m{t_0\in T}, a non
constant morphism \m{\alpha:T\to S} such that \m{\alpha(t_0)=P}, and a regular
sheaf $\ke$ on \m{\alpha^*(\kd)} such that \m{\ke_{|Y}\simeq\E}. It is possible
to choose $\ke$ satisfying conditions (ii) or (iii) of proposition \ref{prop3}.
\end{subsub}

\begin{proof} We will construct $\ke$ satisfying (ii) (the case of (iii) is 
similar). We consider the first canonical filtration of $\ke$ (cf. \ref{cspdc}):
\[0\lra\ke_1\lra\ke\lra\ke_0\lra 0 \ , \]
and the associated \ \m{\sigma\in\Ext^1_{\ko_Y}(\ke_0,\ke_1)} \ and surjective
morphism \m{f:\ke_0\to\ke_1}.
According to \ref{Koda_2} there exist smooth curves \m{T^0}, \m{T^1}, and for
\m{i=0,1}, \m{t_0^i\in T^i}, morphisms \m{\phi_i:T^i\to S} such that
\m{\phi_i(t_0^i)=P}, and vector bundles \m{\kf_i} on \m{\phi_i^*(\kc_{i+1})}
such that \m{\kf_{i,t_0^i}\simeq\ke_i}. Let \ \m{Z=T^0\times_ST^1},
\m{p_0,p_1,q} the projections \m{Z\to T^0}, \m{Z\to T^1}, \m{Z\to S}
respectively.
According to proposition \ref{lem1k}, there exists a unique extension on
\m{q^*(\kd)}
\[0\lra p_1^*(\kf_1)\lra\kf\lra p_0^*(\kf_0)\lra 0\]
corresponding to $f$. We have an exact sequence
\[0\lra\ke_1\lra\kf_{|Y}\lra\ke_0\lra 0 \ , \]
associated to \ \m{\sigma'\in\Ext^1_{\ko_Y}(\ke_0,\ke_1)} (here $Y$ is the
fibre of \m{q^*(\kd)\to Z} over \m{(t_0^0,t_0^1)}). Recall the exact sequence
\xmat{0\ar[r] & \Ext^1_{\ko_C}(\ke_0,\ke_1)\ar[r]^-\iota &
\Ext^1_{\ko_Y}(\ke_0,\ke_1)\ar[r]^-\theta & \Hom(\ke_0,\ke_1)\ar[r] & 0}
(cf. \ref{cstr}). We have \ \m{\theta(\sigma)=\theta(\sigma')}, hence \
\m{\sigma''=\sigma-\sigma'\in\Ext^1_{\ko_C}(\ke_0,\ke_1)}.

Since $f$ is surjective, the induced map
\[\rho:\Ext^1_{\ko_C}(\ke_0,\ke_0)\lra\Ext^1_{\ko_C}(\ke_0,\ke_1)\]
is surjective. Let \ \m{\mu\in\Ext^1_{\ko_C}(\ke_0,\ke_0)} \ such that \
\m{\rho(\mu)=\sigma''}.

According to proposition \ref{prop1b}, there exists a smooth curve $T$,
\m{t_0\in T}, a morphism \m{\beta:T\to Z} such that
\m{\beta(t_0)=(t_0^0,t_0^1)}, and a vector bundle $\ku$ on
\m{\beta^*(q^*(\kc_2))} such that
\m{\ku_{|C}\simeq\ke_0} and \m{\omega_{\beta^\#(p_0^*(\kf_0)),\ku}=\mu}. Let
\m{\alpha=q\circ\beta:T\to S}. We now consider the extension on
\m{\alpha^*(\kd)}
\[0\lra\beta^*(p_1^*(\kf_1))\lra\kv\lra\ku\lra 0 \ . \]
It follows from proposition \ref{prop7} that its restriction to $Y$
\[0\lra\ke_1\lra\ku_{|Y}\lra\ke_0\lra 0\]
is associated to $\sigma$.
\end{proof}

\end{sub}

\sepsec

\section{Maximal reducible deformations and limit sheaves}\label{reach}

Let $C$ be a projective irreducible smooth curve and \m{Y=C_2} a primitive 
double curve, with underlying smooth curve $C$, and associated line bundle $L$ 
on $C$. Let $S$ be a smooth affine curve, \m{P\in S} and \ \m{\pi:\kc\to S} \ a 
maximal reducible deformation of \m{Y}. We keep the notations of 
\ref{prim_reddef}-.

\sepsub

\Ssect{Reachable sheaves}{reach_def}

\begin{subsub}\label{defin_rea}{\bf Definition: } A coherent sheaf $\ke$ on $Y$ 
is called {\em reachable} (with respect to $\pi$) if there exists a smooth 
curve \m{S'}, \m{P'\in S'}, a non constant morphism \m{f:S'\to S} such that 
\m{f(P')=P} and a coherent sheaf $\E$ on $f^*(\kc)$, flat on \m{S'}, such that 
\ \m{\E_{P'}\simeq\ke}.
\end{subsub}

\sepprop

The rest of section 7 is devoted to the proof of the

\sepprop

\begin{subsub}\label{theo6}{\bf Theorem: } Every quasi locally free sheaf on 
$Y$ is reachable.
\end{subsub}

\sepprop

Note that by corollary \ref{coro2} the theorem is true on $\kd$.

\end{sub}

\sepsub

\Ssect{Sheaves concentrated on $C$}{conce_C}

Let \m{E_i}, \m{i=1,2}, be a vector bundle on \m{\kc_i}. If \ \m{0\to 
E_1\to\ke\to E_2\to 0} \ is an exact sequence, the restriction to $Y$
\begin{equation}\label{equ19}0\lra E_{1|Y}=E_{1|C}\lra\ke_{|Y}\lra 
E_{2|Y}=E_{2|C}\lra 0\end{equation}
is also exact.

Recall that \m{L^*} has a canonical section $\bf s$ (defined by theorem 
\ref{theo1}, {\bf 3}). Let
\[\lambda_L:\Ext^1_{\ko_C}(E_{2|C},E_{1|C})\lra\Ext^1_{\ko_C}(E_{2|C}\ot
L,E_{1|C}) \ \]
be the induced linear map. The rest of \ref{conce_C} is dedicated to the proof 
of

\sepprop

\begin{subsub}\label{prop12}{\bf Proposition:} Let
\begin{equation}\label{equ20}0\lra E_{1|C}\lra E\lra E_{2|C}\lra 
0\end{equation}
be an extension on $C$, associated to \ 
\m{\sigma\in\Ext^1_{\ko_C}(E_{2|C},E_{1|C})}. Then there exists an extension \ 
\m{0\to E_1\to\ke\to E_2\to 0} \ on $\kc$ such that \m{\ke_{|Y}} is 
concentrated on $C$ and that $(\ref{equ19})$ is isomorphic to $(\ref{equ20})$ 
if and only if \ \m{\lambda_L(\sigma)=0}.
\end{subsub}

\sepprop

Let \m{\ko_\kc(1)} be an ample line bundle on $\kc$. It \m{X\subset\kc} is a 
closed subvariety, let \ \m{\ko_X(1)=\ko_\kc(1)_{|X}}.

Let \m{n_0} be an integer such that for every \m{n\geq n_0}, \m{E_2(n)} is 
generated by its global sections and 
\m{H^1(E_1(n))=H^1(E_{1|C}(n))=H^1(E_1(n)(-C))=\nsp}. 
Suppose that \m{n\geq n_0}. Let \ \m{k=h^0(E_2(n))}. Then the morphism (of 
sheaves on $\kc$)
\[\tau:k\ko_\kc(-n)=\ko_\kc(-n)\ot H^0(E_2(n))\lra E_2\]
induced by the evaluation morphism is surjective. Let \ \m{\kn=\ker(\tau)}. 
Using the fact that
\[\Hom(E_2,E_1) \ = \ \Ext^1_{\ko_\kc}(k\ko_\kc(-n),E_1) \ = \ \nsp \ , \]
we deduce from the exact sequence
\begin{equation}\label{equ21}
0\lra\kn\lra k\ko_\kc(-n)\lra E_2\lra 0
\end{equation}
this one
\[0\lra\Hom(k\ko_\kc(-n),E_1)\lra\Hom(\kn,E_1)\lra
\Ext^1_{\ko_\kc}(E_2,E_1)\lra 0 \ . \]
According to lemma \ref{lem2}, the restriction of $(\ref{equ21})$ to $Y$
\[0\lra\kn_{|Y}\stackrel{\epsilon}{\lra}k\ko_Y(-n)\lra E_{2|C}\lra 0\]
is also exact, and induces a surjective map
\[\gamma:\Hom(\kn_{|Y},E_{1|C})\lra\Ext^1_{\ko_Y}(E_{2|C},E_{1|C}) \ . \]
Let \ \m{\sigma\in\Ext^1_{\ko_Y}(E_{2|C},E_{1|C})}, 
\m{\phi\in\Hom(\kn_{|Y},E_{1|C})} \ such that \m{\gamma(\phi)=\sigma}, and
\[0\lra E_{1|C}\lra\ke\lra E_{2|C}\lra 0\]
the extension (on $Y$) corresponding to $\sigma$. Then $\ke$ is isomorphic to 
the cokernel of the injective morphism
\[(\epsilon,\phi):\kn_{|Y}\lra k\ko_Y(-n)\oplus E_{1|C} \ . \]

\sepprop

\begin{subsub}\label{lem4}{\bf Lemma: } The sheaf $\ke$ is concentrated on $C$ 
if and only if $\phi$ vanishes on \m{(\kn_{|Y})^{(1)}}.
\end{subsub}

\begin{proof} Let \m{x\in C}, \m{z\in\ko_{Y,x}} an equation of $C$. Then we 
have \ \m{z.k\ko_Y(-n)_x\subset\kn_{|Y,x}}, and \m{\ke_x} is a 
\m{\ko_{C,x}}-module if and only if \ \m{z\ke_x=0}, if and only if for every \ 
\m{(u,e)\in k\ko_Y(-n)_x\times E_{1|C,x}}, we have \ 
\m{z(u,e)=(zu,0)\in\imm(\epsilon_x,\phi_x)}, i.e. if and only if there exists \ 
\m{\nu\in\kn_{|Y,x}} \ such that \ \m{\epsilon_x(\nu)=zu} \ and \ 
\m{\phi_x(\nu)=0}. We have
\begin{eqnarray*}
\epsilon_x(\nu) \ \text{is a multiple of} \ z & \Longleftrightarrow &
z\epsilon_x(\nu)=0\\
&  \Longleftrightarrow & \epsilon_x(z\nu)=0\\
&  \Longleftrightarrow & z\nu=0\\
&  \Longleftrightarrow & \nu\in(\kn_{|Y})^{(1)}_x \ .
\end{eqnarray*}
Lemma \ref{lem4} follows immediately.
\end{proof}

\sepprop

According to lemma \ref{lem1b}, we have an exact sequence
\xmat{0\ar[r] & \kn_{|Y}/(\kn_{|Y})^{(1)}\ar[r]\fleq[d] & 
k\ko_C(-n)\ar[r]^-\rho & E_{2|C}\ar[r] & 0 .\\
& \kn_{|C}/\big((\kn_{|Y})^{(1)}/(\kn_{|Y})_1)
}
Let
\[U \ = \ \kn_{|Y}/(\kn_{|Y})^{(1)}=\kn_{|C}/\big((\kn_{|Y})^{(1)}/(\kn_{|Y})_1)
\ = \ \ker(\rho) \ . \]
Then it follows from lemma \ref{lem4} that:

\begin{subsub}\label{fact1} $\ke$ is concentrated on $C$ if 
and only if $\phi$ can be factorised through $U$. \rm\Nligne
We have a canonical surjection \ 
\m{\beta:\Hom(U,E_{1|C})\to\Ext^1_{\ko_C}(E_{1|C},E_{2|C})}, and for every 
\Nligne \m{\phi\in\Hom(U,E_{1|C})}, if \ \m{0\to E_{1|C}\to E\to E_{2|C}\to 0} 
\ is the extension on $C$ associated to \m{\beta(\phi)}, then $E$ is isomorphic 
to the cokernel of the injective morphism
\[(i,\phi):U\lra n\ko_C(-n)\oplus E_{1|C}\]
(where $i$ is the inclusion).\end{subsub}

Let \m{\kn_1} the kernel of the surjective composite morphism
\xmat{k\ko_{\kc_1}(-n)\ar[r]^-{\tau_1} & 
E_{2|\kc_1}=E_{2|\kz_0}\ar[r] &
E_{2|\kz}}
(cf. \ref{prel2}), with \m{\tau_1=\tau_{|\kc_1}}. Since \ 
\m{\Tor^1_{\ko_{\kc_1}}(\ko_C,E_{2|\kz})=\nsp}, we have an exact sequence
\[0\lra\kn_{1|C}=\kn_{|C}\lra k\ko_C(-n)\lra E_{2|\kz\cap C}\lra 0 \ , \]
and a commutative diagram with exact rows and columns
\begin{equation}\label{equ24}
\xymatrix{& & & 0\ar[d]\\ & 0\ar[d] & & E_{2|C}\ot L\ar[d]\\
0\ar[r] & U\ar[r]\ar[d] & k\ko_C(-n)\ar[r]\fleq[d] & E_{2|C}\ar[r]\ar[d] & 0\\
0\ar[r] & \kn_{1|C}\ar[r]\ar[d] & k\ko_C(-n)\ar[r] & E_{2|\kz\cap 
C}\ar[r]\ar[d] & 0\\
& E_{2|C}\ot L\ar[d] & & 0\\ & 0}
\end{equation}
Let 
\begin{equation}\label{equ22} 0\lra E_1\lra\ke\lra E_2\lra 0
\end{equation}
be an extension on $\kc$, associated to \ 
\m{\sigma\in\Ext^1_{\ko_\kc}(E_2,E_1)}. Recall that
\[\Ext^1_{\ko_\kc}(E_2,E_1) \ \simeq \ \Ext^1_{\ko_{\kc_1}}(E_{2|\kz_0},E_1) \
\simeq \ \Hom(E_{2|\kz_0},E_1(\kz_0)_{|\kz_0})\]
(cf. \ref{prel2}). The exact sequence $(\ref{equ22})$ restricts to an exact 
sequence on $Y$
\[0\lra E_{1|Y}\lra\ke_{|Y}\lra E_{2|Y}\lra 0 \ . \]
Let \ \m{\sigma_0\in\Ext^1(E_{2|C},E_{1|C})} \ be the associated element. 
Recall that we have a canonical exact sequence
\[0\lra\Ext^1_{\ko_C}(E_{2|C},E_{1|C})\stackrel{i}{\lra} \Ext^1(E_{2|C},E_{1|C})
\stackrel{\theta}{\lra} \Hom(E_{2|C}\ot L,E_{1|C})\lra 0\]
(cf. \ref{cstr}). We have a canonical obvious map
\[\theta_1:\Hom(E_{2|\kz_0},E_1(\kz_0)_{|kz_0})\lra\Hom(E_{2|C}\ot L,E_{1|C}) \ 
, \]
and \ \m{\theta_1(\sigma)=\theta(\sigma_0)}.

We have a commutative diagram with exact rows
\begin{equation}\label{equ23}
\xymatrix{0\ar[r] & \kn_{|\kc_1}\ar[r]\flinc[d] & 
k\ko_{\kc_1}(-n)\ar[r]\fleq[d] &
E_{2|\kz_0}\ar[r]\flon[d] & 0\\
0\ar[r] & \kn_1\ar[r] & k\ko_{\kc_1}(-n)\ar[r] & E_{2|\kz}\ar[r] & 0}
\end{equation}
(the exactness of the first row comes from corollary \ref{coro1}), 

From $(\ref{equ23})$ and lemma \ref{lem1f} we deduce the commutative diagram 
with exact rows
\xmat{\Hom(\kn_1,E_1)\ar[r]\flinc[d] & \Ext^1_{\ko_{\kc_1}}(E_{2|\kz},E_1)
=\Hom(E_{2|\kz},E_1(\kz)_{|\kz})\ar[r]\flinc[d]_\Psi & 0\\
\Hom(\kn_{|\kc_1},E_1)\ar[r] & \Ext^1_{\ko_{\kc_1}}(E_{2|\kz_0},E_1)
=\Hom(E_{2|\kz_0},E_1(\kz_0)_{|\kz_0})\ar[r] & 0}
From lemma \ref{lem1d}, the image of $\Psi$ is exactly the space of morphisms 
of vector bundles vanishing on $C$. It follows that:

\begin{subsub}\label{fact2} The extensions $(\ref{equ17})$ such that \m{\ke_Y} 
is concentrated on $C$ arise from the morphisms \ \m{\kn\to E_1} \ that can be 
factorised through \m{\kn_1}. \rm\Nligne
Let \ \m{\phi:\kn\to\kn_1\stackrel{\alpha}{\lra}E_1} \ be such a morphism. Then 
the morphism \ \m{U\to E_{1|C}} \ restriction of $\alpha$ can be used to define
the extension \ \m{0\to E_{1|C}\to\m{\ke_Y}\to E_{2|C}\to 0} (cf. \ref{fact1}).
\Nligne
We have a commutative diagram
\xmat{\Hom(k\ko_C(-n),E_{1|C})\ar[r]\flinc[d] & \Hom(U,E_{1|C})\flon[r]\fleq[d] 
& \Ext^1_{\ko_C}(E_{2|C},E_{1|C})\ar[d]_{\lambda_L}\\
\Hom(\kn_{1|C},E_{1|C})\ar[r] & \Hom(U,E_{1|C})\ar[r]^-\mu & 
\Ext^1_{\ko_C}(E_{2|C}\ot L,E_{1|C}),}
where the two rows come from $(\ref{equ24})$ and are exact. The commutativity 
of the right square follows from lemma \ref{lem1e}. It follows that: {\em a 
morphism \ \m{\alpha:U\to E_{1|C}} \ can be extended to \m{\kn_{1|C}} if and 
only if \ \m{\gamma(\mu(\alpha))=0}.
}\end{subsub}

\sepprop

\begin{subsub}\label{lem5} {\bf Lemma: } The restriction morphism
\[\lambda:\Hom(\kn_1,E_1)\lra\Hom(\kn_{1|C},E_{1|C})\]
is surjective.
\end{subsub}

\begin{proof} Let \m{T=\kz\cap C}. From the exact sequences
\[0\to\kn_1\to k\ko_{\kc_1}(-n)\to E_{2|\kz}\to 0 \ , \quad
0\to\kn_{1|C}\to k\ko_C(-n)\to E_{2|T}\to 0\]
We deduce the commutative diagram with exact rows
\xmat{& & & \Ext^1_{\ko_{\kc_1}}(E_{2|\kz},E_1)\fleq[d]\\
0\ar[r] & \Hom(k\ko_{\kc_1}(-n),E_1)\ar[r]\ar[d]^f & 
\Hom(\kn_1,E_1)\ar[r]\ar[d]^g & \Hom(E_{2|\kz},E_1(\kz)_{|\kz})\ar[r]\ar[d]^h & 
0\\
0\ar[r] & \Hom(k\ko_C(-n),E_{1|C})\ar[r] & \Hom(\kn_{1|C},E_{1|C})\ar[r] &
\Hom(E_{2|T},(E_{1|C})(T)_{|T})\ar[r]\fleq[d] & 0\\
& & & \Ext^1_{\ko_C}(E_{2|\kz\cap C},E_{1|C})
}
The map $f$ is surjective because \ \m{H^1(E_1(n)(-C))=\nsp}, and $h$ is also 
surjective because $\kz$ is made of disjoint affine curves. Hence $g$ is 
surjective.
\end{proof}

\sepprop

{\em Proof of proposition \ref{prop12}.} Suppose that the extension \ \m{0\to 
E_1\to\ke\to E_2\to 0} \ on $\kc$ exists. Then according to \ref{fact2} we have
 \ \m{\lambda_L(\sigma)=0}. This follows also from proposition \ref{lem1k}.

Conversely, suppose that \ \m{\lambda_L(\sigma)=0}. Let \ 
\m{\phi\in\Hom(U,E_{1|C})} \ over $\sigma$ (cf. \ref{fact1}). Then according to 
\ref{fact2}, $\phi$ can be extended to \ \m{\ov{\phi}:\kn_{1|C}\to U}. From 
lemma \ref{lem5}, there exists \ \m{\Phi\in\Hom(\kn_1,E_1)} \ restriction to 
$\ov{\phi}$ on $C$. This $\ov{\phi}$ defines the extension \ \m{0\to 
E_1\to\ke\to E_2\to 0} \ on $\kc$. $\Box$

\end{sub}

\sepsub

\Ssect{Proof of theorem \ref{theo6}}{pr_t6}

Let $\ke$ be a quasi locally free sheaf on $Y$ and
\begin{equation}\label{equ32}0\lra\ke_1\lra\ke\lra\ke_0\lra 0\end{equation}
the exact sequence coming from its first canonical filtration, corresponding
to\Nligne \ \m{\sigma\in\Ext^1_{\ko_Y}(\ke_0,\ke_1)}. Let \
\m{f:\ke_0\ot L\to\ke_1} \ be the associated surjective morphism. Let
\m{\rho:\kd\to\kc} be the blowing-up of $\kz$, inducing \
\m{\phi:\widehat{Y}=C'_2\to Y}, the blowing-up of \m{\kz\cap C} in $Y$ (cf.
\ref{regbl}).
We have then an exact sequence on \m{\widehat{Y}}
\begin{equation}\label{equ25}0\lra\ke_1\ot L^*\lra\widetilde{\ke}
\lra\ke_0\lra 0\end{equation}
corresponding to the first canonical filtration of \m{\widetilde{\ke}} (cf.
\ref{bl_up}). Using corollary \ref{coro2} it is possible to find a smooth curve
$T$, \m{t_0\in T}, a morphism \m{\alpha:T\to C} such that \m{\alpha(t_0)=P}, and
vector bundles \m{F_i} on \m{\alpha^*(\kc_i)}, \m{i=1,2}, such that
\[\F_{2|C} \ \simeq \ \ke_1\ot L^* \ , \quad \F_{1|C} \ \simeq \ \ke_0 \ , \]
and an exact sequence
\begin{equation}\label{equ28}0\lra F_2\lra\E\lra F_1\lra 0\end{equation}
on \m{\alpha^*(\kd)} restricting to $(\ref{equ25})$ on $\widehat{Y}$. It is
possible to extend $f$ to a surjective morphism \m{\ov{f}:F_{1|\kz_0}\to
F_{2|\kz_0}}, corresponding to an extension
\begin{equation}\label{equ29}0\lra F_2\ot\L\lra\ku\lra F_1\lra 0\end{equation}
on \m{\alpha^*(\kc)} (cf. proposition \ref{prop2}). Then by proposition
\ref{prop5}, $(\ref{equ28})$ is the exact sequence
\[0\lra F_2\lra\widetilde{\ku}\lra F_1\lra 0\]
arising from $(\ref{equ29})$ by blowing-up (cf. \ref{regbl}). The restriction
of $(\ref{equ29})$ to $Y$ is an exact sequence
\[0\lra\ke_1\lra\ke'\lra\ke_0\lra 0\]
corresponding to \ \m{\sigma'\in\Ext^1_{\ko_Y}(\ke_0,\ke_1)}, and such that the
associated morphism \m{\ke_0\ot L\to\ke_1} is $f$. We have
\[\sigma-\sigma' \ \in \ \Ext^1_{\ko_C}(\ke_0,\ke_1) \ . \]
Let
\[\lambda:\Ext^1_{\ko_C}(\ke_0,\ke_1)\lra\Ext^1_{\ko_C}(\ke_0,\ke_1\ot L^*)\]
be the map induced by the inclusion \m{L\subset\ko_C}. Then we have \
\m{\sigma-\sigma'\in\ker(\lambda)}, from theorem \ref{theo3}.
From proposition \ref{prop12}, it is possible to find an extension on
$\alpha^*(\kc)$
\[0\lra F_2\lra\widetilde{\kv}\lra F_1\lra0\]
such that its restriction to $Y$
\[0\lra\ke_1\lra\kv_{|Y}\lra\ke_0\lra 0\]
is in fact concentrated on $C$ and corresponds to \m{\sigma-\sigma'}. Let
\m{g:F_{1|\kz_0}\to F_{2|\kz_0}} corresponding to $g$. 
Then the extension
\[0\lra F_2\lra\widetilde{\kv}\lra F_1\lra 0\]
restricts to $(\ref{equ32})$ on $Y$.

\end{sub}

\vskip 1.5cm

\end{document}